\documentclass[11pt]{amsart}
\usepackage{amssymb, latexsym, mathrsfs,   color, amsmath }
\usepackage{setspace}

\usepackage[all]{xy}
\usepackage{mathrsfs}
\usepackage[colorlinks=true, pdfstartview=FitV,linkcolor=blue,citecolor=blue,urlcolor=blue]{hyperref}
\usepackage{chngcntr}
\usepackage{arydshln}
\usepackage{xcolor}       % 提供 \textcolor 命令（也可用 color 宏包）
\usepackage{forest}
\counterwithin{figure}{section}
\advance \topmargin by -\headheight
\advance \topmargin by -\headsep

\evensidemargin \oddsidemargin
\newtheorem {theorem}    {Theorem}[section]

\newtheorem {lemma}      [theorem]    {Lemma}
\newtheorem {corollary}  [theorem]    {Corollary}
\newtheorem {proposition}[theorem]    {Proposition}

\renewcommand{\rm}{\mathrm}

\newcommand{\GGL}{\mathrm{GL}}

\theoremstyle{definition}

\numberwithin{equation}{section}

\begin{document}

	\title{The Fourier kernel and invariant trace formula of $\rm{GL}(3)$}
	
	\date{\today}
	
	\author{Haoyang Wang, Xinghua Cui, Zhifeng Peng*}
	
	\address{School of Mathematical Science, Soochow University, Suzhou 310027, Jiangsu, P.R. China}
	
	\maketitle
	
	\begin{abstract}
		Trace formula is an important method to study the Langlands program.  Arthur obtains the existence of stable trace formula for connected reductive group. In this paper, we will give the explicit invariant trace formula of $\rm{GL}(3)$. For general reductive groups, Arthur gave an inductive proof of the existence of the invariant trace formula and provided explicit examples for $\GGL(2)$. For GL(3), we successively compute the invariant distributions, give all the Fourier kernels for the $\GGL(2)$ cases, and finally obtain the concrete invariant trace formula.
		
	\end{abstract}

\makeatletter

% 一级目录：加粗，间距 3pt
\def\l@section{\@tocline{1}{3pt}{1pc}{}{\bfseries}}

% 二级目录：间距 0pt
\def\l@subsection{\@tocline{2}{0pt}{1pc}{5pc}{}}

% 调整目录格式，保留页码和点线
\def\@tocline#1#2#3#4#5#6#7{%
  \relax
  \ifnum #1>\c@tocdepth
  \else
    \par\addpenalty\@secpenalty\addvspace{#2}%
    \begingroup
    \hyphenpenalty\@M
    \@ifempty{#4}{%
      \@tempdima\csname r@tocindent\number#1\endcsname\relax
    }{%
      \@tempdima#4\relax
    }%
    \parindent\z@
    \leftskip#3\relax
    \advance\leftskip\@tempdima\relax
    \rightskip\@pnumwidth plus4em
    \parfillskip-\@pnumwidth
    #5\leavevmode\hskip-\@tempdima
    #6\nobreak\relax
    \ifnum#1=1
      \hfill
    \else
      \dotfill
    \fi
    \hbox to\@pnumwidth{\@tocpagenum{#7}}\par
    \nobreak
    \endgroup
  \fi
}

\makeatother

    \setcounter{tocdepth}{2}
    \begingroup
\setstretch{1.0}
\setlength{\parskip}{0pt}
\tableofcontents
\endgroup
	
	\section{introduction}
	Arthur’s invariant trace formula relates geometric distributions associated with conjugacy classes to spectral distributions arising from automorphic representations. It provides a fundamental tool for the study of automorphic forms and the Langlands program (\cite{Art88(1)},\cite{Art88(2)}). Its explicit use requires a detailed understanding of the invariant distributions on both sides. The classical rank-one examples of \(\mathrm{GL}(2)\) and \(\mathrm{SL}(2)\) provide basic models for such calculations (\cite{Art05}, \cite{AHS86}), but the passage to \(\mathrm{GL}(3)\) requires additional local harmonic analysis.
	
	In this paper, we give an explicit treatment of the invariant distributions occurring in the trace formula for \(\mathrm{GL}(3)\). A key local input is the calculation of the \(p\)-adic Fourier kernels associated with regular semisimple elements of \(\mathrm{GL}(2)\). These kernels enter the correction terms used to make the geometric distributions invariant. We determine them by choosing suitable auxiliary test functions and applying Arthur’s reciprocity relation, which reduces the problem to the calculation of weighted characters of tempered representations. We then use these local formulas in the inductive construction of the invariant distributions for \(\mathrm{GL}(3)\).
	
	The explicit determination of nonarchimedean Fourier kernels presents a substantial analytic difficulty. In the real case, Hoffmann’s calculations use the differential equations satisfied by the Fourier transforms, together with asymptotic conditions and relations at singular sets (\cite{Hof01},\cite{Hof08}). This approach leads, in particular, to explicit formulas for \(\mathrm{GL}(2,\mathbb R)\) and \(\mathrm{GL}(3,\mathbb R)\). The corresponding local differential equations are unavailable in the \(p\)-adic setting (\cite{Art94}), so the method underlying those calculations cannot be used for the nonarchimedean problem considered here. To the best of our knowledge, explicit formulas for the collection of \(p\)-adic \(\mathrm{GL}(2)\) kernels required in our argument have not previously been available.
	
	Our calculation rests on Arthur’s theorem establishing the existence of Fourier kernels for the invariant distributions associated with weighted orbital integrals (\cite{Art94}, Theorem 4.1), together with the inductive relation between noninvariant and invariant distributions. In the notation of this construction, the latter relation is
    \begin{align}\label{descent1}
		J^L(f)=\sum_{M\in\mathcal L^L(M_0)}
		c(M)c(L)^{-1}\,
		\widehat I^M\!\left(\phi_M^L(f)\right),
\end{align}
	see (\cite{Art81}).
	
	Identity (\ref{descent1}) organizes the correction terms needed to recover the invariant distributions from the noninvariant ones. We apply it in the case of \(\mathrm{GL}(3)\) to reduce the Fourier-kernel calculations required by these correction terms to the corresponding calculations for \(\mathrm{GL}(2)\). Arthur’s existence theorem provides the kernel representations that make this reduction possible. The explicit local problem is therefore to determine the relevant \(\mathrm{GL}(2)\) kernels.
	
	For regular elliptic elements of any Levi subgroup over a \(p\)-adic field, we address this problem by a suitable choice of auxiliary test functions (Lemma \ref{local0} and Lemma \ref{simil}). Arthur’s reciprocity relation (\cite{Art94}, Theorem 4.5) connects the geometric Fourier kernels with the kernels of the dual invariant distributions constructed from weighted characters. Applied with our chosen test functions, it reduces the desired kernel calculation to weighted characters attached to tempered representations. We compute these weighted characters and thereby obtain explicit expressions for the required Fourier kernels. Since the Fourier kernel has nothing to do with the selection of the function, the special test functions serve to determine the kernels, which can then be inserted into (\ref{descent1}) to evaluate the invariant correction terms. This gives a method of computing the needed nonarchimedean kernels without differential equations.
	
	Having determined the regular semisimple terms, we apply Arthur’s limiting formulas
		\[J_M(\gamma, f) = \lim_{a \to 1} \sum_{L \in \mathcal{L}(M)} r_M^L(\gamma, a) J_L(a\gamma, f),\]
	 to obtain expressions for the invariant distributions at arbitrary nonregular semisimple elements. The kernel calculations and the limiting process together provide the geometric invariant distributions needed for the \(\mathrm{GL}(3)\) formula.

     Then we can get the main theorems:
     \begin{theorem}(Theorem \ref{main1})
		Suppose \(f\in\mathcal{H}(G(F_v))\). For \(\gamma_G\in G_{\rm{ell}}\), \[I_G(\gamma_G,f)=|D_G(\gamma_G)|^{\frac{1}{2}}\int_{G_{\gamma_G}(F)\backslash G(F)}f(x^{-1}\gamma_G x)dx.\] 
		For \(\gamma_{21}\in M_{21}(F)_{\rm{ell}}\), the invariant distribution \(I_{M_{21}}(\gamma_{21},f)\) equals 
        \begin{align*}
            &J_{M_{21}}(\gamma_{21},f)\\
-&
|D^{M_{21}}(\gamma_{21})|^{1/2}
\int_{\Pi_2(M_{21}(F))}
\Theta_\tau(\gamma_{21}^{-1})
J_{M_{21}}(\tau,f)\,d\tau\\
-&
\frac12
\int_{\Pi_{\mathrm{unit}}(M_0(F))}
I_{M_{21},M_0}^{M_{21}}(\gamma_{21},\chi)
J_{M_{21}}(\chi^{M_{21}},f)\,d\chi.
        \end{align*}
		For \(\gamma\in M_0(F)\),
		\(I_{M_0}(\gamma,f)\) equals 
		\begin{align*}
			&J_{M_0}(\gamma,f)-
\int_{\Pi_{\mathrm{unit}}(M_0(F))}
\tau(\gamma^{-1})J_{M_0}(\tau,f)\,d\tau\\
-&
\frac12\sum_{w\in\mathscr W}
\int_{\Pi_{\mathrm{unit}}(M_0(F))}
I_{M_0,M_0}^{M_{21}}(\gamma_w,\tau)\,
J_{M_{21}}(\tau^{M_{21}},f)\,d\tau\\
-&
\sum_{w\in\mathscr W}
\int_{T_{\mathrm{disc}}(M_{21})}
I_{M_0,M_{21}}^{M_{21}}(\gamma_w,\tau)\,
J_{M_{21}}(\tau,f)\,d\tau.
		\end{align*}
	\end{theorem}
    \begin{theorem}(Theorem \ref{main2})
		Suppose \(f\in\mathcal{H}(G(\mathbb{A}))\). The geometric side of invariant trace formula of \(\GGL(3)\) is 
		\[
		I_{\rm{geom}}(f) = \lim_{S} \sum_{M \in \mathcal{L}} |W_0^M| |W_0^G|^{-1} \sum_{\gamma \in \Gamma(M)_S} a^M(\gamma,S) I_M(\gamma, f)
		\]
	\end{theorem}
     On the spectral side, we start from Arthur’s invariant spectral expansion
	\begin{align*}
		&I(f)=\sum_{t>0}I_t(f) = \sum_{t>0}I_t^G(f) \\
        &=\sum_{t>0}\sum_{M \in \mathcal{L}(T)} \sum_{M_1 \subset M}\frac{|W_0^{M_1}|}{|W_0^G|} \sum_{\sigma \in \Pi_{\mathrm{disc}}(M_1, t)} a_{\mathrm{disc}}^{M_1}(\sigma)\int_{i\mathfrak{a}_{M_1}^* / i\mathfrak{a}_M^*} r_{M_1}^M(\sigma, \lambda) I_M \left( \mathrm{Ind}_{P_1\cap M}^M(\sigma_\lambda), f \right) d\lambda.
	\end{align*}
	This expansion involves invariant distributions attached to induced representations, Arthur’s discrete spectral coefficients, and factors formed from the normalizing factors of intertwining operators. Combining it with the classification of the relevant representations, we determine the spectral invariant distributions explicitly, including the nontempered cases that occur in the expansion. This spectral calculation, together with the geometric analysis above, yields an explicit description of the invariant distributions on both sides of the trace formula for \(\mathrm{GL}(3)\).
    
    \begin{theorem}(Theorem \ref{spectr})
        For \(f\in\mathcal{H}(G(F_v))\), the invariant distribution \(I_M \left( \mathscr{I}_{P_1\cap M}^M(\sigma_\lambda), f \right)\) equals
\[ \rm{tr}({\pi}(f)),\]if \(M=G,\ \pi\in \Pi_{\rm{disc}}(G,t)\),\[-\sqrt{\frac{3}{2}}\int_{i\mathfrak{a}_{M_{21}}^*/i\mathfrak{a}_G^*}h_{\xi,\eta}(z)\rm{tr}(\mathscr{I}_{P_{21}}(\sigma_{\kappa,z},f))dz,\]if \(M=M_{21},\ \sigma\in\Sigma(M_{21},t)\).
    \end{theorem}
    And then, we have 
    \begin{theorem}\label{main3}
        For any \(f\in \mathcal{H}(G(\mathbb{A}))\), the spectral side of invariant trace formula of \(\GGL(3)\) is 
        \begin{align*}
		I(f) = & \sum_{t>0}\sum_{M \in \mathcal{L}(T)} \sum_{M_1 \subset M}\frac{|W_0^{M_1}|}{|W_0^G|} \sum_{\sigma \in \Pi_{\mathrm{disc}}(M_1, t)} a_{\mathrm{disc}}^{M_1}(\sigma)\int_{i\mathfrak{a}_{M_1}^* / i\mathfrak{a}_M^*} r_{M_1}^M(\sigma, \lambda) I_M \left( \mathrm{Ind}_{P_1\cap M}^M(\sigma_\lambda), f \right) d\lambda.
	\end{align*}
    \end{theorem}
	
	The principal local contribution is thus the explicit calculation of the required \(p\)-adic \(\mathrm{GL}(2)\) Fourier kernels through special test functions, reciprocity, and tempered weighted characters. Their use in the invariant construction shows how local calculations for \(\mathrm{GL}(2)\) can be used to determine invariant distributions for \(\mathrm{GL}(3)\), and provides concrete formulas for further study of the relation between the geometric and spectral sides.

	Section 3 introduces the Fourier kernel. While we rely on Hoffmann's results for the Archimedean case, we compute the Fourier kernel for regular elliptic elements in the 
p-adic case by selecting specific test functions. Section 4 employs Arthur's limit formula to determine the Fourier kernels for all elements. Section 5 is devoted to the geometric side of the invariant trace formula, where we derive each of its terms following Arthur's structural results for arbitrary reductive groups. In Section 6, using the classification of representations and the results of Mœglin and Waldspurger, we derive each term on the spectral side, which ultimately allows us to write out the complete invariant trace formula for \(\GGL(3)\).

	\section{Preliminaries}
	Let \(F\) be a number field of characteristic zero and let \(G\) be the connected
	reductive group \(\GGL(3)\) defined over \(F\). We use the same notation \(G\) for the group
	of \(F\)-rational points \(G(F)\). Fix a minimal Levi subgroup \(M_0\) of \(G\).
	We denote by
	\[
	\mathcal{L}=\mathcal{L}(G)
	\]
	the set of Levi subgroups of \(G\) containing \(M_0\). If \(M\in\mathcal{L}\),
	we write \(\mathcal{L}(M)\) for the set of Levi subgroups of \(G\) containing
	\(M\). Let \( K \) be a fixed maximal compact subgroup of \( G(F) \).

	Every parabolic subgroup \( P \) of \( G \) defined over \( F \) and containing \( M_0 \) has a unique Levi component \( M_P \) containing \( M_0 \). Both \( M_P \) and the unipotent radical \( N_P \) of \(P\) are defined over \( F \). We write \[ \mathcal{F}(M) = \mathcal{F}^G(M) \] for the set of parabolic subgroups \( P \) of \( G \) over \( F \) such that \( M_P \supset M \), and \[ \mathcal{P}(M) = \mathcal{P}^G(M) \] for those \( P \in \mathcal{F}(M) \) with \( M_P = M \). For any \( M \in \mathcal{L} \), the intersection \( K_M = K \cap M(F) \) is a maximal compact subgroup of \( M(F) \).

	Let \(\sigma \in \Pi_2(M(F))\) be a discrete series representation of \(M(F)\).
	The principal series representation \(\pi = \mathcal{I}_P(\sigma_\lambda)\) is 
	induced from the parabolic subgroup \(P\) with Levi component \(M\).

	For \(M\in\mathcal{L}\), let \(A_M\) denote the maximal \(F\)-split torus in
	the center of \(M\), and put
	\[
	\mathfrak a_M
	= \operatorname{Hom}_{\mathbb Z}
	\bigl(X^*(M)_F,\mathbb R\bigr).
	\]
	If \(M\subset L\) are Levi subgroups, we write
	\[
	\mathfrak a_M^L=\mathfrak a_M/\mathfrak a_L.
	\]
	
	Let \(W_0^G\) be the Weyl group of \(G\) relative to \(M_0\), and for
	\(L\in\mathcal{L}\) let \(W_0^L\) denote the corresponding Weyl group of \(L\).
	We write \(|W_0^L|\) for its order. Thus the factor
	\(
	|W_0^L|/|W_0^G|
	\)
	is the Weyl coefficient occurring in Arthur's Fourier expansion.

	Define \[
	\tau_Q^R(H)
	=
	\begin{cases}
		1,
		&
		\alpha(H)>0
		\quad
		\text{for all }\alpha\in\Delta_Q^R,
		\\
		0,
		&
		\text{otherwise},
	\end{cases}
	\] and \[
	\widehat{\tau}_Q^R(H)
	=
	\begin{cases}
		1,
		&
		\varpi(H)>0
		\quad
		\text{for all }\varpi\in\widehat{\Delta}_Q^R,
		\\
		0,
		&
		\text{otherwise}.
	\end{cases}
	\]

	We denote by \(G_{\mathrm{reg}}(F)\) the set of regular semisimple elements of \(G(F)\). For \(M\in\mathcal{L}\),
	let \(\Gamma(M(F))\) denote the space of \(M(F)\)-conjugacy classes in \(M(F)\). Thus
	\[
	\gamma\in \Gamma(M(F))\cap G_{\mathrm{reg}}(F)
	\]
	will always denote a regular conjugacy class represented by an element of
	\(M(F)\).
	
	We write
	\(I_{\mathrm{ell}}(M(F))\) for the elliptic part of the geometric parameter space of \(M(F)\). In
	particular, the variable
	\[
	\gamma\in I_{\mathrm{ell}}(M(F))\cap G_{\mathrm{reg}}(F)
	\]
	in the geometric Fourier expansion is an elliptic regular geometric
	parameter.
	
	Let \(T(L)\) denote Arthur's tempered parameter space associated with a Levi
	subgroup \(L\) (\cite{Art93}), and let
	\[
	T_{\mathrm{disc}}(L)\subset T(L)
	\]
	be its discrete part. If \(\tau\in T(L)\), we denote by
	\(\tau^G\) the representation, or virtual representation, of \(G(F)\) obtained from
	\(\tau\) by normalized parabolic induction.
	
	For \(f\in\mathcal C(G(F))\), let
	\[
	\widehat f(\pi)=f_G(\pi)=\operatorname{tr}\pi(f),
	\qquad
	\pi\in\Pi_{\mathrm{temp}}(G),
	\]
	be its invariant Fourier transform. If \(\tau\in T(L)\), we put
	\[
	\widehat f_L(\tau)
	=\widehat f(\tau^G)=f_L(\tau).
	\]
	Thus \(\widehat f_L\) denotes the restriction of the Fourier transform of
	\(f\) to the spectral component induced from \(L\).
	
	For
	\[
	\gamma\in\Gamma(M(F))\cap G_{\mathrm{reg}}(F),
	\]
	Arthur's invariant distribution associated with \(M\) will be denoted by
	\[
	I_M(\gamma,f)=I_M^G(\gamma,f).
	\]
	Its Fourier kernel is denoted by
	\[
	I_M(\gamma,\tau)
	=I_M^G(\gamma,\tau),
	\qquad
	\tau\in T_{\mathrm{disc}}(L).
	\]
	We denote by
	\(J_M(\gamma,f)\) the weighted orbital integral associated with \(M\). Observe that when \( M = G \), \( J_M(\gamma, f) \) is just the invariant orbital integral over the \( G(F) \)-conjugacy class of \( \gamma \). In this case we shall write
	\[
	f_G(\gamma) = J_G(\gamma, f), \quad \gamma \in G_{\mathrm{reg}}(F),
	\]
	when we want to emphasize the dependence on \( \gamma \). If \( Q \) belongs to \( \mathcal{P}(M) \), and \( f_Q \) is the Schwartz function defined by
	\[
	f_Q(m) = \delta_Q(m)^{\frac{1}{2}} \int_K \int_{N_Q(F)} f(k^{-1} m n k) \, dn \, dk, \quad m \in M(F),
	\]
	on \( M(F) \), the function \( (f_Q)_M \) is independent of \( Q \). We denote it simply by \( f_M \). One sees easily that
	\[
	f_M(\gamma) = f_G(\gamma), \quad \gamma \in M(F) \cap G_{\mathrm{reg}}(F).
	\]
	The
	non-invariant and invariant distributions are related by Arthur's
	invariantization procedure.

	Finally, for \(L\in\mathcal L\) and \(\tau\in T(L)\), we denote by
	\[
	I_L(\tau,f)=I_L^G(\tau,f)
	\]
	the invariant spectral distribution. Its geometric Fourier kernel is
	written
	\[
	I_L(\tau,\gamma)=I_L^G(\tau,\gamma),
	\]
	where
	\[
	\gamma\in I_{\mathrm{ell}}(M(F))\cap G_{\mathrm{reg}}(F).
	\]
	We now recall Arthur’s invariantization procedure from \cite{Art81} and \cite{Art94}.

	Let \(S\) be a finite set of places of \(F\). Suppose that for every \( M \in \mathcal{L}(M_0) \), \( V(M) \) is a second complete topological vector space. Suppose that for every pair \( M \subset L \), we are given a continuous map
	\[
	\phi_M^L: U(L) \to V(M)
	\]
	such that for every \( y \in L(F_S)^1 \),
    \begin{align}\label{41}
        \phi_M^L(f^y) = \sum_{Q \in \mathcal{F}^L(M)} \phi_M^{M_Q}(f_{Q, y})
    \end{align}
	We shall sometimes write \( \phi \) for \( \phi_M^M \). In this case, (\ref{41}) says that \( \phi(f^y) = \phi(f) \) for each \( y \) and \( f \). It follows that for every \( i \) in \( V(M)' \), the distribution
	\[
	\phi'(i): f \to i(\phi(f)), \quad f \in U(M),
	\]
	is invariant.

	Furthermore, assume:
    \par\medskip
\noindent%
\refstepcounter{equation}\label{eq:phi-assumption}%
\textup{(\theequation)}\quad
For every $M \in \mathscr{L}(M_0)$, $\phi$ maps $U(M)$
onto $V(M)$; the image of the transpose, $\phi'$, is the
space of all invariant distributions on $U(M)$.

	The first statement of \eqref{eq:phi-assumption} implies that \( \phi' \) is injective; the second states that any invariant distribution on \( U(M) \) is of the form \( \phi'(\cdot) \). If \( I \) is any invariant distribution on \( U(M) \) we shall let \( \hat{I} \) be the unique element \( i \) in \( V(M)' \) such that \( \phi'(i) = I \).

	\begin{itemize}
		\item Consider the geometric side of the trace formula, \(\mathscr{C}(G(F))\), \(V(M)=\mathscr{I}(M(F_S))\). In this case, we have \[\phi^L_M(f)(\tau)=J_M^L(\tau,f)=\rm{tr}\,(\mathscr{R}^L_M(\tau,P_0)\mathscr{I}^L_P(\tau,f)),\] where \(\tau\in \Pi_{\rm{temp}}(M(F_S))\).
		\item Consider the spectral side of the trace formula, \(U(L)=C_c^\infty(L_{\rm{rel}}(F_S))\), \(V(M)=I_c^\infty(M_{\rm{reg}}(F_S))\) (\cite{Art94}, Definition 3.2). We have \[\phi_M^L(f)(\gamma)=J_M^L(\gamma,f).\]
	\end{itemize}
    To obtain the invariant trace formula, we only need to consider \(f\in\mathcal{H}(G(F_S))\subset\mathcal{H}_{\rm{ac}}(G(F_S))\). For the second case, we shall use a formula associated with the limit formula to extend the result to arbitrary elements.

	These two cases both satisfy the condition (*) and (**). Suppose that we have been given a family of distributions 
	\(\{J^L \in U(L)'\}\) and a family of nonzero complex numbers \(\{c(L)\}\), each indexed by 
	\(\mathcal{L}(M_0)\), such that for any \(L\), \(y \in L(F_S)\), and \(f \in U(L)\),
	\begin{align}\label{***}
		J^L(f^y) = \sum_{Q \in \mathcal{F}^L(M_0)} c(M_Q)c(L)^{-1} J^{M_Q}(f_{Q,y}),
	\end{align}where
	\[
	f_{Q,y}(m)
	=
	\delta_Q^L(m)^{1/2}
	\int_{K\cap L(\mathbb A)}
	\int_{N_Q(\mathbb A)}
	f(k^{-1}mnk)\,
	u'_Q(k,y)\,dn\,dk ,
	\]and \[
	u'_Q(k,y)
	=
	\int_{\mathfrak a_Q/\mathfrak a_L}
	\sum_{\substack{
			R\in\mathcal F^L(M_0)\\
			Q\subset R
	}}
	\tau_Q^R(H)\,
	\widehat{\tau}_R^L
	\bigl(H+H_Q(ky)\bigr)\,dH .
	\]
	
	\section{Fourier kernel for regular elliptic terms}
	In this section, we shall give the Fourier kernel of invariant distribution of \(\GGL(2)\) which will occur in invariant trace formula of $\GGL(3)$. For brevity, we write \(F_S=F_v=F\), for any place \(v\).
	
	\subsection{The work of Arthur}
	In \cite{Art94}, Arthur proves the existence of the Fourier kernel of different distributions in trace formula. These results are proved by Arthur in \cite{Art94}.
	\begin{theorem}\label{A1}
The invariant distribution \(I_M(\gamma)\) on
\(\mathcal{C}(G(F))\) has an expansion
\begin{equation}
I_M(\gamma,f)
=
\sum_{L\in\mathcal{L}}
|W_0^L|\,|W_0^G|^{-1}
\int_{T_{\mathrm{disc}}(L)}
I_M(\gamma,\tau)\,f_L(\tau)\,d\tau,
\end{equation}
for a smooth function
\[
I_M(\gamma,\tau)=I_M^G(\gamma,\tau),
\qquad
\gamma\in\Gamma(M(F))\cap G_{\mathrm{reg}}(F),
\qquad
\tau\in T_{\mathrm{disc}}(L),
\]
which satisfies the symmetry condition
\begin{equation}\label{eq:fourier-kernel-symmetry}
I_{wM}(w\gamma,z w^\vee\tau)
=
\chi_\tau(z)\,I_M(\gamma,\tau),
\qquad
z\in Z_\tau,
\quad
w,w^\vee\in W_0^G.
\end{equation}

If \(\gamma\) lies in a Levi subgroup \(M_1\) of \(M\),
the function also satisfies the descent relation
\begin{equation}\label{eq:fourier-kernel-descent}
I_M(\gamma,\tau)
=
\sum_{S\in\mathcal{L}(M_1)}
d_{M_1}^G(M,S)
\sum_{\substack{
w\in W_0^S\backslash W_0^G\\
wL\subseteq S
}}
I_{M_1}^S(\gamma,w\tau)).
\end{equation}

Finally, if \(F\) is Archimedean, the function satisfies
the growth estimate
\begin{equation}\label{eq:fourier-kernel-growth}
\left|
D_\tau D_\gamma I_M(\gamma,\tau)
\right|
\leq
c(\gamma)\bigl(1+\|\mu_\tau\|\bigr)^n,
\end{equation}
where \(n\) is a positive integer and \(c(\gamma)\) is a
locally bounded function on
\(\Gamma(M(F))\cap G_{\mathrm{reg}}(F)\).
Both depend on the given pair of invariant differential
operators \(D_\gamma\) and \(D_\tau\), transferred from
\(A_M(F)\) and \(i\mathfrak{a}_L^*\), respectively.
\end{theorem}
	
	\begin{corollary}\label{32}
		The weighted orbital integral \( J_M(\gamma) \) on \( \mathscr{C}(G(F)) \) has an expansion
		\[
		J_M(\gamma, f) = \sum_{L \in \mathcal{L}} \sum_{S \in \mathcal{L}(L) \cap \mathcal{L}(M)} \int_{T_{\rm{disc}(L)}} |W_0^L| |W_0^S|^{-1} I_M^S(\gamma, \tau) J_S(\tau, f) \, d\tau.
		\]
	\end{corollary}
	
	\begin{theorem}\label{A2}
The invariant distribution \(I_L(\tau)\) on
\(C_c^\infty(G_{\mathrm{reg}}(F))\) has an expansion
\begin{equation}
I_L(\tau,g)
=
\sum_{M\in\mathcal{L}}
|W_0^M|\,|W_0^G|^{-1}
\int_{\Gamma_{\mathrm{ell}}(M(F))}
I_L(\tau,\gamma)\,g_M(\gamma)\,d\gamma,
\end{equation}
for a smooth function
\[
I_L(\tau,\gamma)=I_L^G(\tau,\gamma),
\qquad
\tau\in T(L),
\qquad
\gamma\in
\Gamma_{\mathrm{ell}}(M(F))
\cap G_{\mathrm{reg}}(F),
\]
which satisfies the symmetry condition
\begin{equation}
I_{w^\vee L}(z w^\vee\tau,w\gamma)
=
\chi_\tau(z)^{-1}I_L(\tau,\gamma),
\qquad
z\in Z_\tau,
\quad
w^\vee,w\in W_0^G.
\end{equation}

If \(\tau\) lies in the image of \(T(L_1)\) in \(T(L)\)
for a Levi subgroup \(L_1\subseteq L\), the function
also satisfies the descent relation
\begin{equation}
I_L(\tau,\gamma)
=
\sum_{S\in\mathcal{L}(L_1)}
d_{L_1}^G(L,S)
\sum_{\substack{
w\in W_0^S\backslash W_0^G\\
wM\subseteq S
}}
I_{L_1}^S(\tau,w\gamma).
\end{equation}

Finally, if \(F\) is Archimedean, the function satisfies
the growth estimate
\begin{equation}
\left|
D_\tau D_\gamma I_L(\tau,\gamma)
\right|
\leq
c(\gamma)\bigl(1+\|\mu_\tau\|\bigr)^n,
\end{equation}
where \(D_\tau\), \(D_\gamma\), \(c(\gamma)\), and \(n\)
are as in \eqref{eq:fourier-kernel-growth}.
\end{theorem}
	
	\begin{corollary}
		The weighted character \( J_L(\tau) \) on \( C_c^\infty(G_{\mathrm{reg}}(F)) \) has an expansion
		\[
		J_L(\tau, g) = \sum_{M \in \mathcal{L}} \int_{I_{\mathrm{ell}}(M(F))} \sum_{S \in \mathcal{L}(M) \cap \mathcal{L}(L)} |W_0^M| |W_0^S|^{-1} I_L^S(\tau, \gamma) J_S(\gamma, g) \, d\gamma.
		\]
	\end{corollary}
	
	\begin{theorem}\label{reciprocity}
		The functions defined by Theorems \ref{A1} and \ref{A2} satisfy a reciprocity relation
		\[
		I_M(\gamma, \tau) = (-1)^{\dim(A_M \times A_L)} i^L(\tau) I_L(\tau^\vee, \gamma),
		\]
		for any pair of points \(\gamma \in \Gamma_{\mathrm{ell}}(M(F)) \cap G_{\mathrm{reg}}(F)\) and \(\tau \in T_{\mathrm{disc}}(L)\).
	\end{theorem}
    Let $M$ be a general linear group over $F$, and let $M_1$
and $L$ be Levi subgroups of $M$. Since the relevant
$R$-groups are trivial, we identify Arthur's spectral
parameters with the corresponding irreducible tempered
representations. Thus, $\tau\in T_{\mathrm{disc}}(L)$
is regarded as an irreducible tempered representation
of $L(F)$.

For
\[
\gamma\in M_1(F)\cap M_{\mathrm{reg}}(F),
\qquad
\tau\in T_{\mathrm{disc}}(L),
\]
we write
\[
I_{M_1,L}^{M}(\gamma,\tau)
:=
I_{M_1}^{M}(\gamma,\tau)
\]
for Arthur's Fourier kernel, with the additional
subscript $L$ recording the Levi subgroup on which
$\tau$ is defined.

With compatible measures and normalizations, the
relation with Hoffmann's notation is
\[
I_{M_1,L}^{M}(\gamma,\tau)
=
\Phi_{M_1,L}^{M}(\gamma,\tau),
\]
where $\tau^{\vee}$ denotes the contragredient
representation of $\tau$.

	\subsection{The Real case}
	Hoffmann obatains the fourier kernel of \(\GGL(2,\mathbb{R})\) in \cite{Hof08} by solving the differential equation 
	\[\chi(z)\Phi_M(t) = \sum_{M' \supset M} \partial_M^{M'}(t, z_{M'})\Phi_{M'}(t),\]which comes from \cite{Art88}.

	\begin{theorem}\cite{Hof08}
		\textit{(i) The Fourier transform of the orbital integral \(I_M^M(a)\) with \(a \in A \cap M_{\text{reg}}\) is given by}
		\[
		\Phi_{M,M}^M(a, \sigma) = 0, \qquad \Phi_{M,A}^M(a, \chi) = \chi(a) + w\chi(a).
		\]
		
		\noindent\textit{(ii) The Fourier transform of the orbital integral \(I_M^{M,ij}(t)\) with \(t \in T \cap M_{\text{reg}}\) vanishes on \(\Pi_{\text{disc}}(M) \setminus \Pi_2(M)\) and is given by}
		\[
		\Phi_{M,M}^{M,ij}(t, \sigma) = \Phi^{M,ij}(\sigma, t)
		\]
		\textit{for \(\sigma \in \Pi_2(M)\) and}
		\[
		\Phi_{M,A}^{M,ij}(t, \chi) = \frac{n_A^M}{2} \left( \frac{\chi(\gamma)\tilde{\chi}(t^c) - \tilde{\chi}(\gamma t^c)}{\sin \pi \lambda_{ij}} + \frac{w\chi(\gamma)w\tilde{\chi}(t^c) - w\tilde{\chi}(\gamma t^c)}{\sin \pi \lambda_{ji}} \right)
		\]
		\textit{for \(\chi \in \Pi_2(A)\), where \(t^c = c^{-1}tc\) and \(n_A^M\) is the length of \(\check{\alpha}_{ij}\) with respect to \(\eta_A^M\). Here we define \(\tilde{\chi}\) on a dense subset of \(A_{\mathbb{C}}\) by}
		\[
		\tilde{\chi}(\operatorname{diag}(a_1, a_2)) = a_1^{\lambda_1} a_2^{\lambda_2}
		\]
		\textit{using the principal branches of complex powers.}
		
		\vspace{1em}
		\noindent\textit{(iii) The Fourier transform of the invariant distribution \(I_A^M(a)\) is given for \(|a_{ij}| > 1\) by}
		\[
		\Phi_{A,M}^M(a, \sigma) = -\Phi^M(\sigma, a)
		\]
		\textit{for \(\sigma \in \Pi_2(M)\),}
		\[
		\Phi_{A,M}^M(a, \chi^M) = \frac{1}{2}\chi(a)
		\]
		\textit{for \(\chi \in \Pi_2(A'_M)\),}
		\[
		\Phi_{A,A}^M(a, \chi) = n_A^M \Big( \chi(a)b(\lambda_{ji}, a_{ji}) + \chi(a)u_{ij}(\chi) + w\chi(a)b(\lambda_{ij}, a_{ji}) + w\chi(a)u_{ij}(w\chi) \Big),
		\]
		\textit{where \(b\) is the function defined in Theorem 8.1 and \(u_{ij}(\chi)\) is the logarithmic derivative of \(r_{ij}(\chi)\) with respect to the identification \((\mathfrak{a}_{\mathbb{C}}^M)^* \to \mathbb{C}\) that sends \(\lambda\) to \(\lambda_{ij}\). If we denote the logarithmic derivative of the \(\Gamma\)-function by \(\psi\), then}
		\[
		u_{ij}(\chi) = \frac{1}{2}\psi\left(\frac{\lambda_{ij} + N}{2}\right) - \frac{1}{2}\psi\left(\frac{\lambda_{ij} + N + 1}{2}\right).
		\]
	\end{theorem}
	Note that the notation \(\Phi^M_{M_1,M_2}\), for \(M_1,M_2=\GGL(2)\) or \(A\), is the Fourier kernel \(I^G_{M_1,M_2}(\gamma,\tau)\), for \(\tau\in \Pi_2(M_2(\mathbb{R)})\).
	
	\subsection{The complex case}
	Let $M=\mathrm{GL}(2,\mathbb C)$, and let $A$ be its diagonal subgroup. We use Hoffmann's notation and compatible measure conventions from \cite{Hof08}, together with Langlands' normalization of intertwining operators. The formulas below follow from \cite{Hof97}, after passing to Arthur's invariant distributions and translating to these conventions. In this case, there is no elliptic regular term. Let \[a=\operatorname{diag}(a_1,a_2)\in A,\qquad a_1\neq a_2,\]
	and write the unitary character as
	\[
	\chi(a)=\prod_{\ell=1}^2
	|a_\ell|^{\lambda_\ell}
	\left(\frac{a_\ell}{|a_\ell|}\right)^{m_\ell},
	\qquad \lambda_\ell\in i\mathbb{R},\quad m_\ell\in\mathbb{Z}.
	\]
	Set
	\[
	z=a_{12},\quad s=\lambda_{12},\quad m=m_1-m_2,\quad
	p=\frac{s+m}{2},\quad q=\frac{s-m}{2},
	\]
	and
	\[
	n_A^M=\left|\eta_A^M\bigl(\operatorname{diag}(1,-1)\bigr)\right|.
	\]
	Then for \(|z|>1\), the Fourier kernel is
	\[
		\begin{aligned}
			\Phi^M_{A,A}(a,\chi)
			=n_A^M\Bigg\{&
			\chi(a)\left[
			\frac{b(-p,z^{-1})+b(-q,\bar z^{-1})}{2}
			+\frac{1}{|m|-s}\right]\\
			+{}&w\chi(a)\left[
			\frac{b(p,z^{-1})+b(q,\bar z^{-1})}{2}
			+\frac{1}{|m|+s}\right]\Bigg\},
	\end{aligned}
	\]
	where \(w\chi(a)=\chi(\operatorname{diag}(a_2,a_1))\) and
	\[
	b(v,u)=\sum_{r\ge1}\frac{u^r}{r+v}.
	\]
	For \(|z|<1\), interchange indices \(1\) and \(2\). For \(s=0\) or \(|z|=1, z\ne1\), take the smooth continuation of the entire right-hand side.
	
	The special spectral term in the full Fourier expansion is
	\[
	\Phi^M_{A,M}\bigl(a,(\xi\otimes\xi)^M\bigr)
	=\frac{1}{2}\,\xi(a_1a_2),
	\qquad \xi\in\widehat{\mathbb{C}^\times},
	\]
	with the spectral measure also following Hoffmann's compatible conventions.
	
	In our case, the number \(n_A^M\) equlas \(\sqrt{2}\).
	
	Now we deal with the case that \(F\) is p-adic field.
	\subsection{The non-Archimedean case}
    \subsection{The case \(\gamma\in M_0(F)\cap G_{\rm{reg}}(F)\)}
	Assume \(G=\GGL(2)\) in this subsubsection.
	By Theorem \ref{reciprocity}, to obtain \(I_{M_0,M_0}(\gamma,\chi)\), we just need to calculate \(I_{M_0,M_0}(\chi^\vee,\gamma)\).
	\begin{lemma}\label{local0}
		For any \(\gamma=\rm{diag}(a,b)\in M_{0,reg}(F)\), there exists a neighbourhood \(U_\gamma\), such that for any \(f\in C_c^\infty(U_\gamma)\), \(J_{M_0}(\gamma,f)=0\).
	\end{lemma}
	\begin{proof}
		Recall that 
		\[J_{M_0}(\gamma,f)=\int_{G(\gamma,F)\backslash G(F)}f(x^{-1}\gamma x)v_{M_o}(x)dx.\]
		By Iwasawa decomposition, \[x=nak,\quad n\in N(F),\ a\in M_0(F),\ k\in K.\]
		Then by [\cite{Art74},(7.1)], we can obtain that 
		\begin{align}
			v_{M_0}(x)=H_0(wn),\quad \text{for}\ w=\left(\begin{matrix}
				0&1\\
				-1&0
			\end{matrix}\right),
		\end{align}
		by the definition of \(v_{M_0}(x)\).
		We just need to prove that if there exists a \(U_\gamma\), such that \(x^{-1}\gamma x\in U_\gamma\), \(x\in M_0(F)K\).

		Set \(g=x^{-1}\gamma x\), we have \[(g-aI)(g-bI)=0.\]
		Define \[e_1(g)=\frac{g-bI}{a-b},\quad e_2(g)=\frac{aI-g}{a-b}.\]
		It is easy to prove that 
		\begin{itemize}
			\item \(e_1(g)+e_2(g)=I\),
			\item \(e_i^2=e_i\), \(i=1,2\),
			\item \(e_1e_2=0\),
			\item \(ge_1(g)=ae_1(g)\), and \(ge_2(g)=be_2(g)\).
		\end{itemize}
		Set \[E_1=e_1(\gamma)=\frac{\gamma-bI}{a-b}=\rm{diag}(1,0),\quad E_2=e_2(\gamma)=\frac{aI-\gamma}{a-b}=\rm{diag}(0,1),\]
		we have \[e_1(g)-E_1=\frac{g-\gamma}{a-b},\quad e_2(g)-E_2=-\frac{g-\gamma}{a-b}.\]
		Write \(\mathfrak{p}\) as the maximal ideal of \(\mathcal{O}_F\), then \((a-b)\mathfrak{p}M_2(\mathcal{O}_F)\) is open.

		Take \(U_\gamma\) is small enough such that \(g\in U_\gamma\), \(g-\gamma\in (a-b)\mathfrak{p}M_2(\mathcal{O}_F)\),
		which means
		\[e_i(g)=E_i\ (mod\ \mathfrak{p}),\quad i=1,2.\]
		Define 
		\begin{align*}
			u_1=\left(\begin{matrix}
				1\\
				0
			\end{matrix}\right),\quad u_2=\left(\begin{matrix}
				0\\
				1
			\end{matrix}\right),
		\end{align*}and \[v_i=e_i(g)u_i,\quad i=1,2.\]We have \[v_i\equiv u_i(mod\ \mathfrak{p}).\]
		Set \(k(g)=(v_1,v_2)\), we have \[k(g)\equiv I(mod\ \mathfrak{p}),\]then \[\rm{det}(k(g))\equiv 1(mod\ \mathfrak{p}),\]which means \(k(g)\in K\), if we write 
		\begin{align*}
			k(g)^{-1}=\frac{1}{\rm{det}(k(g))}\left(\begin{matrix}
				d&-b\\
				-c&a
			\end{matrix}\right),\quad \text{for}\ k(g)=\left(\begin{matrix}
				a&b\\
				c&d
			\end{matrix}\right)\in M_2(\mathcal{O}_F).
		\end{align*}
		Since \[v_i=e_i(g)u_i,\quad ge_1(g)=ae_1(g),\quad ge_2(g)=be_2(g),\]we have\[gv_1=av_1,\quad gv_2=bv_2,\]that is \[gk(g)=k(g)\gamma.\]Therefore, \[xk(g)\in G_\gamma=M_0,\]that is \[x\in M_0K.\]
	\end{proof}
	By lemma \ref{local0}, for any \(\gamma=\rm{diag}(a,b)\in M_{0,reg}(F)\) and \(f_0\in C_c^\infty(U_\gamma)\) we can write 
	\begin{align}\label{R}
		I_{M_0}(\chi^\vee,f_0)=J_{M_0}(\chi^\vee,f_0)=&\rm{tr}\,(\mathscr{R}_{M_0}(\chi^\vee)\mathscr{I}_{P_0}(\chi^\vee,f_0))\\\notag
		=&\sqrt{2}\rm{tr}\,(A'(0)\mathscr{I}_{P_0}(\chi^\vee,f_0)),
	\end{align}
	where \[A(s)=\mathscr{R}_{\overline{P_0}|P_0}(\chi^\vee)^{-1}\mathscr{R}_{\overline{P_0}|P_0}(\chi_s^\vee),\]
	and the equation (\ref{R}) is obtained by the calculation of the limit of \((G,M)\)-family.

	Since we want to calculate the fourier kernel \(I_{M_0}(\gamma,\tau)\), which has nothing to do with the selection of \(f\), we directly write \(f=f_0\in C_c^\infty(U_\gamma)\), for any \(\gamma=\rm{diag}(a,b)\in M_{0,reg}(F)\).

	Set \[\xi=\chi_1\chi_2^{-1},\quad \chi_1,\chi_2\in\Pi_2(\GGL(1,F)),\quad \chi_s = (\chi_1 |\cdot|^{s/2}) \otimes (\chi_2 |\cdot|^{-s/2}).\]By the definition of induced representation, we have
	\begin{align}\label{ind}
		\Phi' \left( \begin{pmatrix} a & b \\ 0 & d \end{pmatrix} y \right) = \chi_1(a)\chi_2(d) \left| \frac{a}{d} \right|^{\frac{s+1}{2}} \Phi'(y),
	\end{align}
	for any \(\Phi'\in \mathscr{V}(\sigma_s)\).

	We consider \(\mathscr{J}_{\overline{P_0}|P_0}(\chi_s)\Phi'(y)\), which is
	\begin{align*}
		\int_{F} \Phi'(w^{-1}n(u)y) e^{\langle s+ \rho_P, H_0(w^{-1}n(u)y) \rangle} du,\quad n(u)=\left(\begin{matrix}
			1&u\\
			&1
		\end{matrix}\right).
	\end{align*}
	Recall that Arthur obtained the weight character by the inner product formula which proved by Langlands (in the global case) :
	\begin{lemma}\cite{Art80}
		Suppose that \(P_1, P'_1 \in \mathcal{P}_\chi\), that \(\phi \in \mathcal{H}^0_{P_1}(\pi)_\chi\), \(\phi' \in \mathcal{H}^0_{P'_1}(\pi')_\chi\) and that \(\zeta\) and \(\zeta'\) are vectors in \(\mathfrak{a}^*_{P_1,\mathbb{C}}\) and \(\mathfrak{a}^*_{P'_1,\mathbb{C}}\) whose real parts are suitably regular points in \((\mathfrak{a}^*_{P_1})^+\) and \((\mathfrak{a}^*_{P'_1})^+\) respectively. Then
		\[
		\int_{G(\mathbb{Q})\backslash G(\mathbb{A})^1} \Lambda^T E(x, \phi, \zeta) \overline{\Lambda^T E(x, \phi', \zeta')} \, dx
		\]
		equals the sum over \(P_2 \in \mathcal{P}_\chi\), \(s \in \Omega(\mathfrak{a}_1, \mathfrak{a}_2)\), and \(s' \in \Omega(\mathfrak{a}'_1, \mathfrak{a}_2)\) of
		\begin{align}\label{inner}
			\mathrm{vol}(\mathfrak{a}_2^G/L_2) \cdot \frac{e^{(s\zeta + s'\bar{\zeta}')(T)}}{\prod_{\alpha \in \Delta_2} (s\zeta + s'\bar{\zeta}')(\alpha^\vee)} (M(s, \zeta)\phi, M(s', \zeta')\phi'). \qquad \square
		\end{align}
	\end{lemma}
    A local analogue was subsequently established by Arthur in \cite{Art91}.

	By the Bruhat decomposition,
\[
G(F)=P_0(F)\sqcup P_0(F)wN_0(F).
\]
For unitary inducing data, restriction to the open
Bruhat cell gives the usual unitary noncompact
realization of the induced representation, with
compatible Haar measures. We compute the above trace
in this realization, using the functions
$n\mapsto\Phi(wn)$; the complementary Bruhat cell has
measure zero in $P_0(F)\backslash G(F)$.

	Since \begin{align*}
		w^{-1}n(u)wn(x)=\left(\begin{matrix}
			-u^{-1}&1\\
			&-u
		\end{matrix}\right)w\left(\begin{matrix}
			1&x-u^{-1}\\
			&1
		\end{matrix}\right),
	\end{align*}
	by (\ref{ind}), the function \(\mathscr{J}_{\overline{P_0}|P_0}(\chi_s)\Phi(wn(x))\) equals
	\begin{align}
		\int_{F}\xi(u)|u|^{s-1}\Phi'(wn(x+u))du,
	\end{align}
	if we replace \(u\) by \(-u^{-1}\).

	Define \[h_s(x)=\xi(x)|x|^{s-1},\quad \text{and}\quad\Phi'_w(n(u))=\Phi'(wn(u))\]thus \[\mathscr{J}_{\overline{P_0}|P_0}(\chi_s)\Phi_w'=h_s*\Phi_w'.\]
	Now we write the fourier transform as \(\mathcal{F}\), such that \(\hat{f}=\mathcal{F}f\). Note that we shall use this two notations to write Fourier transform to emphasize the key point. Consider \[\mathcal{F}\mathscr{J}_{\overline{P_0}|P_0}(\chi_s)\Phi_w'(u)=c(s)\widehat{h_s}(y)\mathcal{F}\Phi_w'(u),\]by the fact that \(\widehat{h_s*\Phi_w'}(u)=\widehat{h_s}(u)\hat{\Phi}_w'(y)\), where \[\widehat{h_s}(u)=\int_F\xi(x)|x|^{s-1}\psi(xu)dx.\]
	Since the Fourier transform is unitary isomorphism, we have\[\mathcal{F}\mathscr{J}_{\overline{P_0}|P_0}(\chi_s)\mathcal{F}^{-1}=\widehat{h_s},\]which is easy to check by the definition of the trace of operator.

	Now we calculate the function \(\hat{h}_s(u)\), for any \(u\in F\).

	\noindent\textbf{Case 1:} For \(\xi\) is unramified.

	Set \(z=\xi(\varpi)\), where \(\mathfrak{p}_F=\varpi\mathcal{O}_F\). The integral \(\hat{h}_s(u)\) is the Tate integral.

	By the decomposition
	\[
	F^\times = \bigsqcup_{\ell \in \mathbb{Z}} \varpi^\ell \mathcal{O}_F^\times,
	\]for any \(x\in F^\times\), we write
	\[
	x = \varpi^l x_0,\ u=\varpi^nu_0, \qquad l = v(x),\ n=v(u) \qquad x_0,u_0 \in \mathcal{O}_F^\times.
	\]
	
	Since
	\[
	dx = q^{-\ell} d^\times x, \qquad |x|^{s-1} = q^{-\ell(s-1)}, \qquad \xi(x) = z^\ell,
	\]
	
	it follows that
	\[
	\xi(x)|x|^{s-1} dx = (z q^{-s})^\ell dx_0.
	\]
	
	Consequently,
	\begin{equation}
		\widehat{h_s}(u) = \sum_{\ell \in \mathbb{Z}} (z q^{-s})^\ell \int_{\mathcal{O}_F^\times} \psi(\varpi^{\ell+n} x_0 u_0) \, d^\times x_0. 
	\end{equation}
	Since \(\psi\) is trivial on \(\mathcal{O}_F\) and is not trivial on \(\varpi^{-1}\mathcal{O}_F\), we have
	\[
	\int_{\mathcal{O}_F^\times} \psi(\varpi^{\ell+n} x_0 u_0) \, d^\times x = 
	\begin{cases} 
		1 - q^{-1}, & \ell + n \geq 0, \\
		-q^{-1}, & \ell + n = -1, \\
		0, & \ell + n < -1.
	\end{cases}
	\]
	Therefore,
	\[
	\widehat{h_s}(u) = (1 - q^{-1}) \sum_{\ell=-n}^{\infty} (zq^{-s})^\ell - q^{-1}(zq^{-s})^{-n-1}.
	\]
	In the region of convergence
	\[
	|zq^{-s}| < 1,
	\]
	the geometric series yields
	\begin{align*}
		\widehat{h_s}(u) &= (1 - q^{-1}) \frac{(zq^{-s})^{-n}}{1 - zq^{-s}} - q^{-1}(zq^{-s})^{-n-1} \\
		&= \frac{1 - z^{-1}q^{s-1}}{1 - zq^{-s}} (zq^{-s})^{-n}.
	\end{align*}
	Since \(\xi\) is unramified,
	\[
	\xi(u)^{-1}|u|^{-s} = z^{-n}q^{ns} = (zq^{-s})^{-n},
	\]
	Thus,
	\[
	\widehat{h_s}(u) = c_s(z) \xi(u)^{-1}|u|^{-s},
	\]
	where \[c_s(z) = \frac{1 - z^{-1}q^{s-1}}{1 - zq^{-s}}.\]
	The function \(c_s(z)\) can also be written as
	\[
	\frac{L(s, \xi)}{L(1-s, \xi^{-1})},
	\]
	where the local unramified \(L\)-functions are given by
	\[
	L(s, \xi) = \frac{1}{1 - zq^{-s}}, \qquad L(1-s, \xi^{-1}) = \frac{1}{1 - z^{-1}q^{s-1}}.
	\]
	Recall that\[\mathscr{R}_{\overline{P}_0|P_0}(\chi_s)=r^{-1}_{\overline{P}_0|P_0}(\chi_s)\mathscr{J}_{\overline{P_0}|P_0}(\chi_s),\]and  Arthur gives the construction of the factor \(r_{\overline{P}_0|P_0}\) in \cite{Art13}. By calculation, we have 
	\[r_{\overline{P}_0|P_0}(\chi_s)=\frac{L(s, \xi)}{\epsilon(s, \xi, \psi) L(1+s, \xi)},\]and in the case of unramified, \[\epsilon(s, \xi, \psi) = 1.\]Therefore,\[r_{\overline{P}_0|P_0}(\chi_s)=\frac{1 - zq^{-1-s}}{1 - zq^{-s}},\]which means
    \begin{align}\label{fourier1}
    \mathcal{F}\mathscr{R}_{\overline{P_0}|P_0}(\chi_s)\Phi_w'(u)=c^R_s(z) \xi(u)^{-1}|u|^{-s} \hat{\Phi}_w'(u),
    \end{align}
	for \[c_s^R(z)=\frac{1 - z^{-1}q^{s-1}}{1 - zq^{-1-s}}.\]
	Consider \[\frac{d}{ds}\rm{log}\, (c_s^R(z)|u|^{-s})|_{s=0},\]which is 
    \begin{align}\label{unrami}
        -(\rm{log}\ q)S_q(z)+v(u)\rm{log}\ q,
    \end{align}
	where \[S_q(z)=\frac{q^{-1}z}{1 - q^{-1}z} + \frac{q^{-1}z^{-1}}{1 - q^{-1}z^{-1}}.\]
	\noindent\textbf{Case 2:} For \(\xi\) is ramified.

	Now suppose
	\[
	\xi = \xi_0 \nu_a, \qquad a = \xi(\varpi), 
	\]and let the conductor of \(\xi\) be \(c\).

	We define the Gauss integral by
	\[
	\tau(\xi_0, \psi) = \int_{\mathcal{O}_F^\times} \xi_0(t) \psi(\varpi^{-c} t) \, dt.
	\]
	Since the unit integral is supported on a single valuation shell, we have
	\[
	\int_F \xi(x) |x|^{s-1} \psi(xu) \, dx = c_s^{\mathrm{raw}}(\xi) \xi(u)^{-1} |u|^{-s},
	\]
	where
	\[
	c_s^{\mathrm{raw}}(\xi) = \tau(\xi_0, \psi) a^{-c} q^{cs}.
	\]
	Since $\xi$ is ramified, its local $L$-factor is
identically $1$. Hence the normalizing factor becomes
\[
r_{\overline{P}_0|P_0}(\chi_s)
=
\frac{L(s,\xi)}
{\epsilon(s,\xi,\psi)L(1+s,\xi)}
=
\epsilon(s,\xi,\psi)^{-1}.
\]
For our Fourier transform convention, Tate's local
functional equation gives
\[
c_s^{\mathrm{raw}}(\xi)
=
\epsilon(1-s,\xi^{-1},\psi)
=
\xi(-1)\epsilon(s,\xi,\psi)^{-1}.
\]
Therefore,
\[
c_s^{\mathrm{raw}}(\xi)
=
\xi(-1)r_{\overline{P}_0|P_0}(\chi_s).
\]
	By similar calculation in the last case, we have \[\mathcal{F}\mathscr{J}_{\overline{P_0}|P_0}(\chi_s)\Phi'_w(u)=|u|^{-s}\mathcal{F}\Phi'_w,\]that is, our desired result in this case is
    \begin{align}\label{rami}
        (\log\,q)\ n.
    \end{align}
	Next, we shall find a correspondence of the order of \(u\).

	For \(\chi=\chi_1\otimes\chi_2\), we write
	\[
	\chi_i = \theta_i \nu_{z_i},\quad i=1,2,
	\]
	where \(\theta_i\) is the unit part and \(\nu_{z_i}\) is the unramified part and \(v_i(\varpi)=z_i\). We define the relative character
	\[
	\xi = \chi_1 \chi_2^{-1} = \xi_0 \nu_z,
	\]
	where
	\[
	\xi_0 = \theta_1 \theta_2^{-1}, \qquad z = z_1 z_2^{-1}.
	\]

	Since we need to calculate \(I_{M_0}(\chi^\vee,f)=I_{M_0}(\chi^{-1},f)\), for any \[t=\rm{diag}(a,b)\in M_0(F),\]we have  
	\begin{align*}
		\chi^{-1}(t)=\chi^{-1}_1(a)\chi^{-1}_2(b)&=(\chi_1^{-1}(a)\chi_2(a))\cdot\chi_2^{-1}(ab)\\
		&=\xi^{-1}(a)\cdot(\chi_2^{-1}\circ\rm{det}).
	\end{align*}
	Thus 
	\[\chi_1^{-1} \otimes \chi_2^{-1} = (\xi^{-1} \otimes 1) \cdot (\chi_2^{-1} \circ \det).\]
	Set \(\tau=\chi_2^{-1} \circ \det,\) we have\[\mathscr{I}_{P_0}(\chi^\vee)\simeq \tau\otimes\mathscr{I}_{P_0}(\xi^{-1} \otimes 1).\]
	Now we use the compact model (like Arthur used in \cite{Art81}), set the induced representation
	\[
	V= \left\{ \Phi : K \to \mathbb{C} : \Phi(pk) = \theta_1^{-1}(a)\theta_2^{-1}(b)\Phi(k) \right\},
	\]
	for
	\[
	p = \begin{pmatrix} a & x \\ 0 & b \end{pmatrix} \in P_0 \cap K.
	\]
	Set \[K(n)=1+\varpi^nM_2(\mathcal{O}_F),\  V_n=V^{K(n)},\quad \text{for}\ n\in \mathbb{Z}_{>0}\]and \[K(0)=K,\]
	thus \[K(n+1)\subset K(n).\]
	Write \(c\) as the conductor of \(\xi\), and for the unramified case, define the conductor \(c\) equals \(0\).
	
	For the relative character \(\xi\), define another space \(V\):
	\[
	V_{\mathrm{rel}} = \{ \Phi : K \longrightarrow \mathbb{C} : \Phi(pk) = \xi_0^{-1}(a)\Phi(k) \}.
	\]
	Equivalently,
	\[
	V_{\mathrm{rel}} = \operatorname{Ind}_{P_0 \cap K}^K (\xi_0^{-1} \otimes 1).
	\]
	Define
	\[
	\mathcal{T} : V \longrightarrow V_{\mathrm{rel}}
	\]
	by
	\[
	(\mathcal{T}\Phi)(k) = \tau(k)^{-1}\Phi(k).
	\]
	
	\noindent Let
	\[
	p = \begin{pmatrix} a & x \\ 0 & b \end{pmatrix} \in P_0 \cap K.
	\]
	
	\noindent Then
	\begin{align*}
		(\mathcal{T}\Phi)(pk) &= \tau(pk)^{-1}\Phi(pk) \\
		&= \tau(p)^{-1}\tau(k)^{-1}\theta_1^{-1}(a)\theta_2^{-1}(b)\Phi(k).
	\end{align*}
	\noindent Since
	\[
	\tau(p) = \theta_2^{-1}(ab),
	\]
	it follows that
	\begin{align*}
		\tau(p)^{-1}\theta_1^{-1}(a)\theta_2^{-1}(b) &= \theta_2(ab)\theta_1^{-1}(a)\theta_2^{-1}(b) \\
		&= \theta_1^{-1}(a)\theta_2(a) \\
		&= \xi_0^{-1}(a).
	\end{align*}
	\noindent Thus,
	\[
	(\mathcal{T}\Phi)(pk) = \xi_0^{-1}(a)(\mathcal{T}\Phi)(k).
	\]
	
	\noindent Hence,
	\[
	\mathcal{T}\Phi \in V_{\mathrm{rel}}.
	\]
	\noindent Since \(\tau\) is a unitary character, multiplication by \(\tau^{-1}\) preserves the inner product. Therefore,
	\[
	\mathcal{T} : V\xrightarrow{\sim} V_{\mathrm{rel}}
	\]
	is a unitary isomorphism.

	Hence, we have
	\[
	\mathcal{T} \mathscr{I}_{P_0}(\chi^{-1},g) \mathcal{T}^{-1} = \tau(g) \mathscr{I}_{P_0}(\xi^{-1}\otimes 1,g).
	\]
	Set \(V_{\rm{rel},n}=V_{\rm{rel}}^{K(n)}\). 
	\begin{lemma}
		For any \(n\in \mathbb{Z}_{\ge0}\), \(V_{\rm{rel},n}\) is an
		invariant subspace.
	\end{lemma}
	\begin{proof}
		It is easy to check that \(K(n) \triangleleft K\).
		
		Take any \(h\in K(n)\) and \(v\in V_{\rm{rel},n}\) we have
		\begin{align*}
			&\operatorname{Ind}_{P_0 \cap K}^K (\xi_0^{-1} \otimes 1)(h)\operatorname{Ind}_{P_0 \cap K}^K (\xi_0^{-1} \otimes 1)(k)v \\=&\operatorname{Ind}_{P_0 \cap K}^K (\xi_0^{-1} \otimes 1)(k)\operatorname{Ind}_{P_0 \cap K}^K (\xi_0^{-1} \otimes 1)(k^{-1}hk)v\\
			=&\operatorname{Ind}_{P_0 \cap K}^K (\xi_0^{-1} \otimes 1)(k)v.
		\end{align*}
		Therefore, \(V_{\rm{rel},n}\) is a \(K\)-invariant subspace.
	\end{proof}
	Similarly, we can prove this lemma for \(V_n\).
	\begin{lemma}
		\(V_{\rm{rel},n}=0\) for any \(0\le n<c\).
	\end{lemma}
	\begin{proof}
        For $n=0<c$, choose $a\in\mathcal O_F^\times$ such that $\xi_0(a)\neq1$. If $\Phi\in V_{\mathrm{rel},0}=V_{\mathrm{rel}}^K$, then $\Phi$ is constant on $K$. Hence \[\Phi(1)=\Phi\bigl(\operatorname{diag}(a,1)\bigr)=\xi_0(a)^{-1}\Phi(1).\]
        Since $\xi_0(a)\neq1$, we obtain $\Phi(1)=0$, and therefore $\Phi=0$. Thus $V_{\mathrm{rel},0}=0$.

		Now assume that $1\leq n<c$. By the definition of the conductor, there exists
		\[
		a \in 1 + \varpi^n \mathcal{O}_F
		\]
		such that
		\[
		\xi_0(a) \neq 1.
		\]
		
		\noindent Take
		\[
		p = \begin{pmatrix} a & 0 \\ 0 & 1 \end{pmatrix} \in P_0 \cap K(n).
		\]
		
		\noindent If \(\Phi \in V_{\rm{rel},n}\), then 
		\[
		\Phi(pk) = \xi_0^{-1}(a)\Phi(k).
		\]
		
		\noindent On the other hand, since \(K(n) \triangleleft K\),
		\[
		pk = k(k^{-1}pk), \qquad k^{-1}pk \in K(n).
		\]
		
		\noindent Since \(\Phi\) is right \(K(n)\)-invariant,
		\[
		\Phi(pk) = \Phi(k).
		\]
		
		\noindent Hence,
		\[
		(1 - \xi_0^{-1}(a))\Phi(k) = 0.
		\]
		
		\noindent Since \(\xi_0(a) \neq 1\), we obtain
		\[
		\Phi(k) = 0
		\]
		for all \(k \in K\). Therefore,
		\[
		V_{\rm{rel},n} = 0, \qquad n < c.
		\]
	\end{proof}
	If we similarly define the normalized intertwining operator \(\mathscr{R}_{\rm{rel},\overline{P}_0|P_0}\) on \(V_{\rm{rel}}\) and function \(A_{\rm{rel}}(s)\), we have 
	\[\mathcal{T}\mathscr{R}_{\overline{P}_0|P_0}\mathcal{T}^{-1}=\mathscr{R}_{\rm{rel},\overline{P}_0|P_0},\quad \mathcal{T}A(s)\mathcal{T}^{-1}=A_{\rm{rel}}(s),\]and \[\mathcal{T}A'(0)\mathcal{T}^{-1}=A'_{\rm{rel}}(0).\]
	Recall that we have the filtration \[0\subset V_{\rm{rel},c}\subset V_{\rm{rel},c+1}\subset \cdots\subset V_{\rm{rel},n}\subset\cdots,\]define \(V_{\rm{rel},c-1}=0\), and there exists a subspace \(E_{\rm{rel},n}\) of \(V_{\rm{rel},n}\) such that \[E_{\text{rel},n} = V_{\text{rel},n} \ominus V_{\text{rel},n-1},,\quad n\ge c.\]Then we have 
    \[V_{\text{rel},n} = \bigoplus_{j=c}^{n} E_{\text{rel},j},\]
    and a filtration \[0\subset E_{\rm{rel},c}\subset E_{\rm{rel},c+1}\subset \cdots\subset E_{\rm{rel},n}\subset\cdots.\]
	\begin{lemma}\label{n>c}
		For \(n\ge c\) and \(n\ge 1\), \[V_{\rm{rel},n} \simeq \operatorname{Ind}_{P_0(\mathcal{O}_F/\varpi^n\mathcal{O}_F)}^{\operatorname{GL}(2,\mathcal{O}_F/\varpi^n\mathcal{O}_F)} (\xi_{0,n}^{-1} \otimes 1).\]
	\end{lemma}
	\begin{proof}
		It is obvious that the kernel of
		\[
		\rho_n : \mathcal{O}_F^\times \longrightarrow (\mathcal{O}_F/\varpi^n\mathcal{O}_F)^\times,
		\]is
		\(1 + \varpi^n\mathcal{O}.\)

		Therefore, by the first isomorphism theorem,
		\[
		\mathcal{O}_F^\times / (1 + \varpi^n\mathcal{O}_F) \simeq (\mathcal{O}_F/\varpi^n\mathcal{O}_F)^\times.
		\]
		Since
		\[
		\xi_0|_{1+\varpi^n\mathcal{O}_F} = 1,
		\]
		we have that if
		\[
		u \equiv v \pmod{\varpi^n}, \qquad u, v \in \mathcal{O}_F^\times,
		\]
		then
		\[
		uv^{-1} \in 1 + \varpi^n\mathcal{O}_F.
		\]
		Consequently,
		\[
		\xi_0(u)\xi_0(v)^{-1} = \xi_0(uv^{-1}) = 1.
		\]Hence,
		\[
		\xi_0(u) = \xi_0(v).
		\]
		In other words, \(\xi_0(u)\) depends only on the residue class of \(u\) modulo \(\varpi^n\). Thus, we define
		\[
		\xi_{0,n} : (\mathcal{O}_F/\varpi^n\mathcal{O}_F)^\times \longrightarrow \mathbb{C}^\times
		\]
		by
		\[
		\xi_{0,n}(\bar{u}) = \xi_0(u).
		\]In other words, we can descend \(\xi_0\) to \(\bar{\xi}_0\).

		Consider the map
		\[r_n:\GGL(2,\mathcal{O}_F)\rightarrow \GGL(2,\mathcal{O}_F/\varpi^n\mathcal{O}_F),\]whose kernel is
		\[\{k:k \equiv I \pmod{\varpi^n}\}=K(n).\]
		Thus we have
		\[
		\operatorname{GL}(2,\mathcal{O}_F) / K(n) \simeq \operatorname{GL}(2,\mathcal{O}_F / \varpi^n \mathcal{O}_F).
		\]
		Write
		\[
		G_n = \operatorname{GL}(2,\mathcal{O}_F / \varpi^n \mathcal{O}_F),
		\]
		Then,
		\[
		K / K(n) \simeq G_n,
		\]which is a finite group.

		For \(\Phi\in V_{\rm{rel},n}\),
		\[
		\Phi(kh) = \Phi(k), \qquad h \in K(n).
		\]
		In other words, \(\Phi\) is constant on each right coset
		\[
		kK(n).
		\]
		Therefore, if
		\[
		k_1K(n) = k_2K(n),
		\]we have
		\[
		\Phi(k_2) = \Phi(k_1h) = \Phi(k_1),
		\]for some \(h \in K(n)\).

		Thus, \(\Phi(k)\) depends only on the coset
		\[
		kK(n) \in K/K(n).
		\]
		Consequently, we may define
		\[
		\bar{\Phi} : K/K(n) \to \mathbf{C}
		\]
		by
		\[
		\bar{\Phi}(kK(n)) = \Phi(k).
		\]
		Since
		\[
		K/K(n) \simeq G_n,
		\]
		we may also regard it as
		\[
		\bar{\Phi} : G_n \to \mathbf{C}.
		\]
		In the compact model, \(\Phi\) satisfies
		\[
		\Phi(pk) = (\xi_0^{-1} \otimes 1)(p)\Phi(k), \quad p \in P_0(F) \cap K.
		\]
		For \(p = \begin{pmatrix} a & x \\ 0 & b \end{pmatrix} \in P_0(F) \cap K\), where \(a, b \in \mathcal{O}_F^\times\), we have \(\delta_{P_0}(p) = 1\). Thus the condition becomes
		\[
		\Phi\left(\begin{pmatrix} a & x \\ 0 & b \end{pmatrix} k\right) = \xi_0^{-1}(a)\Phi(k).
		\]
		Write
		\[
		P_n = P_0(\mathcal{O}_F/\varpi^n\mathcal{O}_F) \subset G_n.
		\]
		Since \(\xi_0\) descends to \(\xi_{0,n} : (\mathcal{O}_F/\varpi^n\mathcal{O}_F)^\times \to \mathbf{C}^\times\), we can descend \(\Phi\) to
		\[
		\bar{\Phi}(\bar{p}\bar{k}) = (\xi_{0,n}^{-1} \otimes 1)(\bar{p})\bar{\Phi}(\bar{k}), \qquad \bar{p} \in P_n.
		\]
		This is precisely the defining condition for a finite-group induced representation:
		\[
		\bar{\Phi} \in \operatorname{Ind}_{P_n}^{G_n}(\xi_{0,n}^{-1} \otimes 1).
		\]
	\end{proof}
	Define\[R_n=\mathcal{O}_F/\varpi^n\mathcal{O}_F.\]
	\begin{lemma}
		For any \(n\ge c\), the induced representation \(V_{\rm{rel},n}\) is multiplicity-free.
	\end{lemma}
	\begin{proof}
		\textbf{Case 1: \(c\ge1\)}.

		Since \(\rm{Ind}_{P_n}^{G_n}(\xi_{0,n}^{-1}\otimes 1)\) is a representation induced from a one-dimensional character, we have 
		\[\rm{dim}\,V_{\rm{rel,n}}=[G_n:P_n].\]
		Set \[\mathbf{P}^1(R_n) = \{ (x, y) \in R_n^2 : R_nx + R_ny = R_n \} / R_n^\times,\]where \(R_nx\) is the principal ideal generated by x in \(R_n\).

		We first prove that \[G_n/P_n\simeq \mathbf{P}^1(R_n).\]
		The left action of \(G_n\) on the projective line \(\mathbf{P}^1(R_n)\) is defined by
		\[
		g\cdot[x : y]= [g(x, y)^T].
		\]
		Take the base point \([1 : 0]\). For \(g = \begin{pmatrix} a & b \\ c & d \end{pmatrix}\), we have
		\[
		g[1 : 0]= [a : c].
		\]
		Now consider the stabilizer of the base point:
		\[
		g[1 : 0]= [1 : 0].
		\]
		This holds if and only if
		\[
		[a : c] = [1 : 0],
		\]
		which is equivalent to \(c = 0\). Thus
		\[
		\operatorname{Stab}_{G_n}([1: 0]) = P_n.
		\]
		By the orbit-stabilizer theorem,
		\[
		\mathbf{P}^1(R_n) \simeq G_n/P_n.
		\]It is easy to check that 
		\[|R_n^2|-|\varpi R_n^2|=q^{2n}-q^{2n-2},\]and
		\[|R_n^\times|=q^{n}-q^{n-1},\]
		thus \[\rm{dim}\,V_{\rm{rel},n}=|\mathbf{P}^1(R_n)|=\frac{q^{2n}-q^{2n-2}}{q^{n}-q^{n-1}}=q^{n}+q^{n-1},\]and 
		\begin{align*}
			\rm{dim}\,E_{\rm{rel},n}=\rm{dim}\,V_{\rm{rel},n}-\rm{dim}\,V_{\rm{rel},n-1}=\begin{cases}
				q^c+q^{c-1},&if\ n=c\\
				q^n-q^{n-2},&if\ n>c.
			\end{cases}
		\end{align*}
		Next, we compute
		\(
		\operatorname{End}_{G_n}(V_{\rm{rel},n}).
		\) Define \[g_j = \begin{pmatrix} 1 & 0 \\ \varpi^j & 1 \end{pmatrix}, \qquad 0 \leq j \leq n.\]
		For
		\[
		g = \begin{pmatrix} a & b \\ c & d \end{pmatrix} \in G_n,
		\]
		define
		\[
		j(g) = v_n(c) = \min\{v(\tilde{c}), n\},
		\]
		where \(\tilde{c} \in \mathcal{O}\) is any lift of \(c\), and we set \(v_n(0) = n\).
		
		If
		\[
		b_1 = \begin{pmatrix} \alpha & \beta \\ 0 & \delta \end{pmatrix}, \qquad b_2 = \begin{pmatrix} r & s \\ 0 & t \end{pmatrix} \in P_n,
		\]
		then
		\[
		(b_1 g b_2)_{21} = \delta c r.
		\]
		Since \(\delta, r \in R_n^\times\), it follows that
		\[
		j(b_1 g b_2) = j(g).
		\]
		Hence \(j(g)\) is a double coset invariant.
		
		On the other hand, by
		\[
		G_n / P_n \simeq \mathbf{P}^1(R_n),
		\]
		the double cosets \(P_n \backslash G_n / P_n\) correspond to the orbits of \(P_n\) on \(\mathbf{P}^1(R_n)\). Let \([x : y] \in \mathbf{P}^1(R_n)\) and set \(j = v_n(y)\). If \(j > 0\), then \(y\) is not a unit, so \(x \in R_n^\times\). Thus
		\[
		[x : y] = [1 : x^{-1}y].
		\]
		Since \(v_n(x^{-1}y) = j\), there exists \(u \in R_n^\times\) such that
		\[
		x^{-1}y = u\varpi^j.
		\]
		Therefore,
		\[
		\begin{pmatrix} 1 & 0 \\ 0 & u^{-1} \end{pmatrix} [1 : u\varpi^j] = [1 : \varpi^j].
		\]
		If \(j = 0\), then \(y \in R_n^\times\), so
		\[
		[x : y] = [xy^{-1} : 1].
		\]
		Furthermore,
		\[
		\begin{pmatrix} 1 & 1 - xy^{-1} \\ 0 & 1 \end{pmatrix} [xy^{-1} : 1] = [1 : 1],
		\]
		which shows that this point lies in the \(P_n\)-orbit of \([1 : 1]\). Consequently, each orbit contains a unique element of the form
		\[
		[1 : \varpi^j], \qquad 0 \leq j \leq n.
		\]
		This point is precisely the projective point determined by the first column of \(g_j\). Combining this with the double coset invariance of \(j(g)\), we obtain
		\[
		G_n= \bigsqcup_{j=0}^n P_n g_j P_n.
		\]
		By Mackey theory,
		\[
		\operatorname{End}_{G_n} \left( \operatorname{Ind}_{P_n}^{G_n}(\xi_{0,n}^{-1} \otimes 1)\right) \simeq \bigoplus_{j=0}^n \operatorname{Hom}_{H_j} ((\xi_{0,n}^{-1} \otimes 1)|_{H_j}, (\xi_{0,n}^{-1} \otimes 1)^{g_j}|_{H_j}),
		\]
		where
		\[
		H_j = P_n \cap g_j P_n g_j^{-1},\quad 
		(\xi_{0,n}^{-1} \otimes 1)^{g_j}(h) = (\xi_{0,n}^{-1} \otimes 1)(g_j^{-1} h g_j).
		\]
		Since \(\xi_{0,n}^{-1} \otimes 1\) and \((\xi_{0,n}^{-1} \otimes 1)^{g_j}\) are both one-dimensional characters,
		\[
		\dim \operatorname{Hom}_{H_j} ((\xi_{0,n}^{-1} \otimes 1)|_{H_j}, (\xi_{0,n}^{-1} \otimes 1)^{g_j}|_{H_j}) = 
		\begin{cases} 
			1, & (\xi_{0,n}^{-1} \otimes 1)|_{H_j} = (\xi_{0,n}^{-1} \otimes 1)^{g_j}|_{H_j}, \\
			0, & \text{otherwise}.
		\end{cases}
		\]
		Therefore, it suffices to determine when the two characters agree on \(H_j\). For
		\[
		b = \begin{pmatrix} a & x \\ 0 & d \end{pmatrix} \in P_n,
		\]
		by a direct computation,
		\[
		g_j^{-1} b g_j = \begin{pmatrix} a + \varpi^j x & x \\ \varpi^j(d - a - \varpi^j x) & d - \varpi^j x \end{pmatrix}.
		\]
		The element \(b\) lies in \(H_j\) if and only if
		\[
		g_j^{-1} b g_j \in P_n.
		\]
		This is equivalent to its lower-left entry being zero:
		\[
		\varpi^j(d - a - \varpi^j x) = 0 \quad \text{in } R_n.
		\]
		For such \(b\), we have
		\[
		(\xi_{0,n}^{-1}\otimes1)(b) = \xi_{0,n}^{-1}(a)
		\]
		and
		\[
		(\xi_{0,n}^{-1}\otimes1)^{g_j}(b) = \xi_{0,n}^{-1}(a + \varpi^j x).
		\]
		Hence the ratio of the two characters is
		\begin{align*}
			\frac{(\xi_{0,n}^{-1}\otimes1)(b)}{(\xi_{0,n}^{-1}\otimes1)^{g_j}(b)} &= \frac{\xi_{0,n}^{-1}(a)}{\xi_{0,n}^{-1}(a + \varpi^j x)} \\
			&= \xi_{0,n} \left( \frac{a + \varpi^j x}{a} \right) \\
			&= \xi_{0,n} \left( 1 + \frac{\varpi^j x}{a} \right).
		\end{align*}
		If \(j \geq c\), then
		\[
		\frac{\varpi^j x}{a} \in \varpi^c R_n.
		\]
		Since \(\xi_{0,n}\) is trivial on \(1 + \varpi^c \mathcal{O}_F\), we have
		\[
		\xi_{0,n} \left( 1 + \frac{\varpi^j x}{a} \right) = 1.
		\]
		Therefore,
		\[
		(\xi_{0,n}^{-1} \otimes 1)|_{H_j} = (\xi_{0,n}^{-1} \otimes 1)|_{H_j}.
		\]
		Thus,
		\[
		\dim \operatorname{Hom}_{H_j}((\xi_{0,n}^{-1} \otimes 1)|_{H_j}, (\xi_{0,n}^{-1} \otimes 1)^{g_j}|_{H_j}) = 1.
		\]
		Conversely, suppose \(1\le j < c\). Since \(c\) is the conductor of \(\xi_{0,n}\), \(\xi_{0,n}\) is not identically \(1\) on
		\[
		1 + \varpi^j \mathcal{O}_F.
		\]
		Hence there exists \(x \in R_n\) such that
		\[
		\xi_{0,n}(1 + \varpi^j x) \neq 1.
		\]
		Take
		\[
		a = 1, \qquad d = 1 + \varpi^j x.
		\]
		Then
		\[
		d - a - \varpi^j x = 0.
		\]
		Therefore,
		\[
		b = \begin{pmatrix} 1 & x \\ 0 & 1 + \varpi^j x \end{pmatrix}
		\]
		lies in \(H_j\).

		However, on this element,
		\[
		\frac{(\xi_{0,n}^{-1} \otimes 1)(b)}{(\xi_{0,n}^{-1} \otimes 1)^{g_j}(b)} = \xi_{0,n}(1 + \varpi^j x) \neq 1.
		\]
		Therefore,
		\[
		(\xi_{0,n}^{-1} \otimes 1)|_{H_j} \neq (\xi_{0,n}^{-1} \otimes 1)^{g_j}|_{H_j}
		\]
		and
		\[
		\dim \operatorname{Hom}_{H_j}((\xi_{0,n}^{-1} \otimes 1)|_{H_j}, (\xi_{0,n}^{-1} \otimes 1)^{g_j}|_{H_j}) = 0.
		\]
        For $j=0$, choose $u\in R_n^\times$ with $\xi_{0,n}(u)\ne1$. Set\[b=\begin{pmatrix}1&u-1\\0&u\end{pmatrix}.\] Since \[g_0^{-1}bg_0=\begin{pmatrix}u&u-1\\0&1\end{pmatrix}\in P_n,\]we have $b\in H_0$. Moreover,\[(\xi_{0,n}^{-1}\otimes1)(b)=1,\qquad(\xi_{0,n}^{-1}\otimes1)^{g_0}(b)=\xi_{0,n}(u)^{-1}\ne1.\]Thus the two characters differ on $H_0$.

		In summary,
		\[
		\dim \operatorname{Hom}_{H_j}((\xi_{0,n}^{-1} \otimes 1)|_{H_j}, (\xi_{0,n}^{-1} \otimes 1)^{g_j}_{H_j}) = 
		\begin{cases} 
			0, & j < c, \\
			1, & c \leq j \leq n.
		\end{cases}
		\]
		Consequently,
		\begin{align*}
			\dim \operatorname{End}_{G_n}(V_{\rm{rel},n}) &= \sum_{j=0}^n \dim \operatorname{Hom}_{H_j}((\xi_{0,n}^{-1} \otimes 1)|_{H_j}, (\xi_{0,n}^{-1} \otimes 1)^{g_j}|_{H_j}) \\
			&= \sum_{j=c}^n 1 \\
			&= n - c + 1.
		\end{align*}
		Let \(Q_j : V_{\rm{rel},n} \to V_{\rm{rel},j}\) be the orthogonal projection. Because \(V_{\rm{rel},j}\) is \(K\)-stable, \(Q_j\) is \(K\)-equivariant, so \(Q_j \in \operatorname{End}_K(V_{\rm{rel},n}) = \operatorname{End}_{G_n}(V_{\rm{rel},n})\). Since the subspaces are nested,
		\[
		V_{\rm{rel},c} \subseteq V_{\rm{rel},c+1} \subseteq \dots \subseteq V_{\rm{rel},n},
		\]
		it follows that
		\[
		Q_i Q_j = Q_{\min(i,j)}.
		\]
		In particular, these projections commute:
		\[
		Q_i Q_j = Q_j Q_i.
		\]
		
		We now show that \(Q_c, Q_{c+1}, \dots, Q_n\) are linearly independent. Suppose
		\[
		\sum_{j=c}^n a_j Q_j = 0.
		\]
		Take a nonzero vector \(v_n \in E_{\rm{rel},n}\). Since \(v_n \perp V_{\rm{rel},n-1}\), we have \(Q_j v_n = 0\) for \(j < n\), while \(Q_n v_n = v_n\). Applying the equation to \(v_n\) yields \(a_n v_n = 0\), so \(a_n = 0\). Repeating this argument with a nonzero vector \(v_{n-1} \in E_{\rm{rel},n-1}\) gives \(a_{n-1} = 0\). By downward induction, we obtain
		\[
		a_c = a_{c+1} = \dots = a_n = 0.
		\]
		Hence these projections are linearly independent.

		Since we have already shown that
		\[
		\dim \operatorname{End}_{G_n}(V_{\rm{rel},n}) = n - c + 1,
		\]
		it follows that
		\[
		\operatorname{End}_{G_n}(V_{\rm{rel},n}) = \operatorname{span}_{\mathbb{C}}\{Q_c, Q_{c+1}, \dots, Q_n\}.
		\]
		Consequently,
		\(
		\operatorname{End}_{G_n}(V_{\rm{rel},n})
		\)
		is a commutative algebra.

		Since \(G_n\) is a finite group, \(V_{\rm{rel},n}\) is completely reducible over \(\mathbb{C}\). We write
		\[
		V_{\rm{rel},n} \simeq \bigoplus_{\rho \in \widehat{G_n}} \rho \otimes \mathbb{C}^{m_\rho},
		\]
		where \(m_\rho\) denotes the multiplicity of the irreducible representation \(\rho\) in \(V_{\rm{rel},n}\).
		
		Therefore,
		\[
		\operatorname{End}_{G_n}(V_{\rm{rel},n}) \simeq \bigoplus_{\rho \in \widehat{G_n}} M_{m_\rho}(\mathbb{C}).
		\]
		If \(m_\rho \geq 2\) for some \(\rho\), then
		\[
		M_{m_\rho}(\mathbb{C})
		\]
		is not a commutative algebra, it follows that
		\[
		m_\rho \leq 1
		\]
		for all \(\rho\). Hence,
		\[
		V_{\rm{rel},n} \text{ is multiplicity-free.}
		\]
		\textbf{Case 2: c=0.}

		If \(n=0\), it is easy to check that \(\rm{dim}\,E_{\rm{rel},0}=1.\)

		If \(n=1\), similar as the proof of case 1, \[\rm{dim}\, V_{\rm{rel},1}=q+1,\quad \rm{dim}\,E_{\rm{rel},1}=q.\]

		If \(n\ge 2\), the proof is the same as case 1.
		
	\end{proof}
	\begin{lemma}\label{schur}
		For any \(n\ge c\), \(E_{\rm{rel},n}\) is irreducible and multiplicity-free.
	\end{lemma}
	\begin{proof}
		For \(c\ge1\), since 
		\[\dim \operatorname{End}_{G_n}(V_{\rm{rel},n}) = n - c + 1,\]and
		\[V_{\rm{rel},n}=E_{\rm{rel},c}\oplus E_{\rm{rel},c+1}\oplus\cdots\oplus E_{\rm{rel},n},\]it is obvious that \[E_{\rm{rel},n} \text{ is multiplicity-free and irreducible.}\]
		For \(c=0\), the proof is similar.
	\end{proof}
	\begin{lemma}
		For any \(n\ge c\), 
		\begin{align}\label{En}
			A(s)|_{E_{\rm{rel},n}}=\begin{cases}
				\rm{Id}_{E_{\rm{rel},0}},&c=0,n=0,\\
				C_z^{-1}(s)q^{ns}\rm{Id}_{E_{\rm{rel},n}},&c=0,n\ge1,\\
				q^{ns}\rm{Id}_{E_{\rm{rel},n}},&c\ge1,n\ge1.
			\end{cases}
		\end{align}
	\end{lemma}
	\begin{proof}
		\textbf{Case 1: \(c\ge1.\)}

		We only need to choose some special function to get the constants in (\ref{En}) by Lemma \ref{schur}.

		Fix a nontrivial additive character \(\psi\) of conductor \(\mathcal{O}_F\), 
		\[
		\psi : F \longrightarrow \mathbb{C}^\times,
		\]
		such that
		\[
		\psi|_{\mathcal{O}_F} = 1, \qquad \psi|_{\varpi^{-1}\mathcal{O}_F} \neq 1.
		\]
		Take the Haar measure \(dx\) normalized so that
		\[
		\operatorname{vol}(\mathcal{O}_F) = 1.
		\]
		Define
		\[
		\Phi_n(x) = \mathbf{1}_{\mathcal{O}_F}(x)\psi(\varpi^{-n}x).
		\]
		It is a nonzero, compactly supported, locally constant function, hence it gives a vector in the induced representation.
		
		We now show that
		\[
		\Phi_n \in E_{\rm{rel},n}.
		\]
		Let
		\[
		\bar{n}(x) = \begin{pmatrix} 1 & 0 \\ x & 1 \end{pmatrix}.
		\]
		Recall that we can take any \(\Phi_w(n(x))=\Phi(wn(x))\), since \(wn(x)w^{-1}w=\bar{n}(-x)w\), we can take
		\[
		f(x) = \Phi(\bar{n}(x)).
		\]
		Take any
		\[
		k = \begin{pmatrix} a & b \\ c' & d \end{pmatrix} \in K(n).
		\]
		Since
		\[
		K(n) = 1 + \varpi^n M_2(\mathcal{O}_F),
		\]
		we have
		\[
		a, d \in 1 + \varpi^n \mathcal{O}_F, \qquad b, c' \in \varpi^n \mathcal{O}_F.
		\]
		For \(x \in \mathcal{O}\), we have the decomposition
		\[
		\bar{n}(x)k = \begin{pmatrix} \frac{\rm{det}\,k}{bx+d} & b \\ 0 & bx+d \end{pmatrix} \bar{n}\left(\frac{ax+c'}{bx+d}\right).
		\]
		Write
		\[
		x_k = \frac{ax+c'}{bx+d}.
		\]
		Since \(x \in \mathcal{O}_F\), \(b \in \varpi^n \mathcal{O}_F\), and \(d \in 1 + \varpi^n \mathcal{O}_F\), we have
		\[
		bx+d \in 1 + \varpi^n \mathcal{O}_F.
		\]
		Furthermore,
		\begin{align*}
			x_k - x &= \frac{ax+c' - x(bx+d)}{bx+d} \\
			&= \frac{x(a-d) + c' - bx^2}{bx+d}.
		\end{align*}
		The numerator lies in \(\varpi^n \mathcal{O}_F\) and the denominator is a unit, so
		\[
		x_k - x \in \varpi^n \mathcal{O}_F.
		\]
		On the other hand,
		\[
		\rm{det}\, k \in 1 + \varpi^n \mathcal{O}_F,
		\]
		whence
		\[
		\frac{\rm{det}\, k}{bx+d} \in 1 + \varpi^n \mathcal{O}_F.
		\]
		Since \(n \geq c\), we have
		\[
		\xi_{0,n}|_{1+\varpi^n\mathcal{O}_F} = 1.
		\]
		Moreover, \(bx+d\) and \(\rm{det}\,k/(bx+d)\) are units, so the values of both the modular character and the unramified character equal \(1\). Hence, for \(x \in \mathcal{O}_F\),
		\[
		(\rm{Ind}_{P_0}^{G_n}(\xi_{0,n}\otimes 1)(k)\Phi_n)(x) = \Phi_n(x_k).
		\]
		Therefore, \(x_k\) can be written as
		\[
		x_k = x + \varpi^n u, \qquad u \in \mathcal{O}_F.
		\]
		Thus,
		\begin{align*}
			\Phi_n(x_k) &= \psi(\varpi^{-n}x_k) \\
			&= \psi(\varpi^{-n}x + u) \\
			&= \psi(\varpi^{-n}x),
		\end{align*}
		since \(\psi|_{\mathcal{O}_F} = 1\). Therefore,
		\[
		\Phi_n(x_k) = \Phi_n(x),
		\]which means \[\Phi_n\in V_{\rm{rel},n}.\]
		Let \(n > c\) and take any \(h \in V_{\rm{rel},n-1}\). For any
		\[
		t \in \varpi^{n-1}\mathcal{O}_F,
		\]
		we have
		\[
		\bar{n}(t) \in K(n-1).
		\]
		Since \(h\) is \(K(n-1)\)-invariant and
		\[
		\bar{n}(x)\bar{n}(t) = \bar{n}(x+t),
		\]
		it follows that
		\[
		h(x+t) = h(x), \qquad t \in \varpi^{n-1}\mathcal{O}_F.
		\]
		Hence \(h|_{\mathcal{O}_F}\) is constant on each coset
		\[
		a + \varpi^{n-1}\mathcal{O}_F.
		\]
		Since \(\Phi_n\) is supported on \(\mathcal{O}_F\), it suffices to compute the inner product
		\[
		\int_{\mathcal{O}_F} \Phi_n(x)\overline{h(x)}\,dx.
		\]
		Decompose \(\mathcal{O}_F\) into cosets modulo \(\varpi^{n-1}\mathcal{O}_F\). For each coset \(a + \varpi^{n-1}\mathcal{O}_F\), since \(h\) is constant on it, we have
		\begin{align*}
			\int_{a+\varpi^{n-1}\mathcal{O}_F} \psi(\varpi^{-n}x)\overline{h(x)}\,dx
			&= \overline{h(a)}\psi(\varpi^{-n}a) \int_{\varpi^{n-1}\mathcal{O}_F} \psi(\varpi^{-n}t)\,dt.
		\end{align*}
		Let
		\[
		t = \varpi^{n-1}u.
		\]
		Then
		\[
		dt = q^{-(n-1)}du,
		\]
		so
		\begin{align*}
			\int_{\varpi^{n-1}\mathcal{O}_F} \psi(\varpi^{-n}t)\,dt &= q^{-(n-1)} \int_{\mathcal{O}_F} \psi(\varpi^{-1}u)\,du \\
			&= 0.
		\end{align*}
		The last equality follows from the fact that \(\psi(\varpi^{-1}u)\) is a nontrivial additive character on \(\mathcal{O}_F\). Therefore,
		\[
		\langle \Phi_n, h \rangle = 0
		\]
		for any \(h \in V_{\rm{rel},n-1}\), i.e.,
		\[
		\Phi_n \perp V_{\rm{rel},n-1}.
		\]
		Combining this with \(\Phi_n \in V_{\rm{rel},n}\), we obtain
		\[
		\Phi_n \in E_{\rm{rel},n}.
		\]
		When \(n = c\), by the convention \(V_{\rm{rel},c-1} = 0\), so from \(\Phi_c \in V_{\rm{rel},c}\) we directly obtain
		\[
		\Phi_c \in E_{\rm{rel},c}.
		\]Consider
		\begin{align*}
			\hat{\Phi}_n(y) &= \int_{\mathcal{O}_F} \psi(\varpi^{-n}x)\psi(xy)\,dx \\
			&= \int_{\mathcal{O}_F} \psi\big((y + \varpi^{-n})x\big)\,dx.
		\end{align*}
		Hence
		\[
		\hat{\Phi}_n(y) = \mathbf{1}_{-\varpi^{-n}+\mathcal{O}_F}(y).
		\]
		Therefore,
		\[
		\operatorname{supp} \hat{\Phi}_n = -\varpi^{-n} + \mathcal{O}_F.
		\]
		If
		\[
		y = -\varpi^{-n} + u, \qquad u \in \mathcal{O}_F,
		\]
		then
		\[
		y = -\varpi^{-n}(1 - \varpi^n u).
		\]
		Since \(n \geq c \geq 1\), we have
		\[
		1 - \varpi^n u \in \mathcal{O}_F^\times.
		\]
		Thus, on the entire support of the Fourier transform,
		\[
		v(y) = -n.
		\]
		Consequently,
		\[
		|y|^{-s} = q^{-ns}
		\]
		holds identically on the entire support of \(\hat{\Phi}_n\).

		\textbf{Case 2: \(c=0\).}
		
		For \(n=0\), since \[\rm{dim}\, E_{\rm{rel},0}=1,\] it is easy to check that \[A(s)|_{E_{\rm{rel},0}}=Id_{E_{\rm{rel},0}}.\]
		For \(n\ge1\), the proof is the same as Case 1.

		Now by the Tate integral we discussed before, after calculating \(\mathcal{F}A(s)\), we can obtain the desired result.
	\end{proof}
	By the last lemma, we have
	\begin{lemma}\label{315}
		For \(n\ge c\),
		\begin{align*}
			\rm{tr}(\mathscr{R}_{\rm{rel},M_0}(\xi_0^{-1}\otimes 1)|_{E_{\rm{rel},n}})=\begin{cases}
				\sqrt{2}c(\rm{log}\,q)(q^c+q^{c-1}),&c\ge1, n=c,\\
				\sqrt{2}n(\rm{log}\,q)(q^n-q^{n-2}),&c\ge1,n>c\\
				0,&c=0,n=0\\
				\sqrt{2}q(\rm{log}\,q)(S_q(z)+n),&c=0,n=1\\
				\sqrt{2}(q^n-q^{n-2})(\rm{log}\,q)(S_q(z)+n),&c=0,n\ge 2.
			\end{cases}
		\end{align*}
	\end{lemma}
	\begin{lemma}
		For any \(\gamma\in M_0(F)_{\rm{reg}}\) and \(f\in C_c^\infty(U_\gamma)\), \(\rm{tr} ( \mathscr{R}_{\rm{rel},M_0}(\xi_{0,n}^{-1}\otimes 1) \mathscr{I}_{P_0}(\xi_{0,n}^{-1}\otimes 1,f) )\) equals
		\begin{align*}
			\begin{cases}
				\sqrt{2}(\rm{log}\,q) \sum_{n \ge c} n \operatorname{tr} \left( \mathscr{I}_{P_0}(\xi_{0,n}^{-1}\otimes 1,f)|_{E_{\mathrm{rel},n}} \right), & c \ge 1, \\
				\sqrt{2}(\rm{log}\,q) \sum_{n \ge 1} (S_q(z) + n) \operatorname{tr} \left( \mathscr{I}_{P_0}(\xi_{0,n}^{-1}\otimes 1,f)|_{E_{\mathrm{rel},n}} \right), & c = 0.
			\end{cases}
		\end{align*}
	\end{lemma}
	For any \(x \in G(F)\), take an orthonormal basis \(\{v_{n,j}\}\) of \(V_n\), and define
	\[
	F_n(x) = \sum_j \langle \mathscr{I}_{P_0}(\xi_{0,n}^{-1}\otimes 1,x)v_{n,j}, v_{n,j} \rangle.
	\]
	When \(\mathscr{I}_{P_0}(\xi_{0,n}^{-1}\otimes 1,x)\) preserves \(V_n\), this is the usual restricted trace:
	\[
	F_n(x) = \rm{tr} \left( \mathscr{I}_{P_0}(\xi_{0,n}^{-1}\otimes 1,x)|_{V_n} \right).
	\]
	Similarly, take an orthonormal basis \(\{e_{n,j}\}\) of \(E_n\) and define 
	\[
	T_n(x) = \sum_{j} \langle \mathscr{I}_{P_0}(\xi_{0,n}^{-1}\otimes 1,x) e_{n,j}, e_{n,j} \rangle.
	\]
	Since
	\[
	E_{\rm{rel},n} = V_{\rm{rel},n} \ominus V_{\rm{rel},n-1},
	\]
	taking the union of orthonormal bases of \(V_{n-1}\) and \(E_n\) as a basis of \(V_n\), we can obtain that
	\[
	T_n(x) = F_n(x) - F_{n-1}(x).\]
	Consider \[\mathscr{I}_{P_0}(\xi_{0,n}^{-1}\otimes 1,f)=\int_{G(F)}\mathscr{I}_{P_0}(\xi_{0,n}^{-1}\otimes 1,x)f(x)dx,\]
	which equals the sum of 
	\begin{align}\label{weyl}
		&&\frac{1}{2} \int_{M_0(F)\backslash G(F)} \mathscr{I}_{P_0}(\xi_{0,n}^{-1}\otimes 1,g)^{-1}\left[ \int_{M_0(F)_{\mathrm{reg}}} f(g^{-1}mg) |D_G(m)| \mathscr{I}_{P_0}(\xi_{0,n}^{-1}\otimes 1,m) \, dm \right] \\\notag
		&&\mathscr{I}_{P_0}(\xi_{0,n}^{-1}\otimes 1,g) \, dg
	\end{align}
	and 
	\begin{align*}
		\int_{G(F)_{\rm{ell}}}|D_G(g)|^{\frac{1}{2}}\mathscr{I}_{P_0}(\xi_{0,n}^{-1}\otimes 1,g)f_G(g)dg.
	\end{align*}
	\begin{lemma}
		There exists a neighborhood \(U'_{\gamma}\), such that \(U_\gamma\subset U_\gamma'\), for \(f\in C_c^\infty(U_\gamma)\), the integral (\ref{weyl}) equals 
		\begin{align}\label{weyl2}
			\frac{1}{2} \int_{K} \left[ \int_{M_0(F)_{\mathrm{reg}}} f(k^{-1}mk) |D_G(m)| \mathscr{I}_{P_0}(\xi_{0,n}^{-1}\otimes 1,m) \, dm \right] \, dk.
		\end{align}
	\end{lemma}
	\begin{proof}
		We first prove that if \(g^{-1}mg\in U_\gamma\), for \(g\in G(F)\), then \(g\in M_0(F)K\).
		
		The proof is similar as Lemma \ref{local0}. For
		\[
		\gamma = \begin{pmatrix} \gamma_1 & 0 \\ 0 & \gamma_2 \end{pmatrix} \in M_0(F)_{\mathrm{reg}}, \quad \gamma_1 \neq \gamma_2,
		\]we shall prove that there exist locally analytic maps
		\[
		\lambda_i : U_\gamma \longrightarrow F^\times, \quad i = 1, 2,
		\]
		and
		\[
		k : U_\gamma \longrightarrow K,
		\]
		such that for every \(y \in U_\gamma\),
		\[
		k(y)^{-1} y k(y) = \begin{pmatrix} \lambda_1(y) & 0 \\ 0 & \lambda_2(y) \end{pmatrix}.
		\]
		In other words,
		\[
		y = k(y) t(y) k(y)^{-1}, \quad k(y) \in K, \quad t(y) \in M_0(F)_{\mathrm{reg}}.
		\]
		The proof is as follows.
		For \(y \in G(F)\), denote its characteristic polynomial by
		\[
		P_y(X) = \det(XI - y).
		\]
		Consider the function
		\[
		\Phi(y, \lambda) = \det(\lambda I - y).
		\]
		At the point \((\gamma, \gamma_1)\),
		\[
		\Phi(\gamma, \gamma_1) = 0,\quad \frac{\partial \Phi}{\partial \lambda}(\gamma, \gamma_1) = \gamma_1 - \gamma_2 \neq 0.
		\]
		Similarly,
		\[
		\Phi(\gamma, \gamma_2) = 0, \quad \frac{\partial \Phi}{\partial \lambda}(\gamma, \gamma_2) = \gamma_2 - \gamma_1 \neq 0.
		\]
		By the Non-Archimedean Implicit Function Theorem, there exists a neighborhood \(U_1\) of \(\gamma\) such that for every \(y \in U_1\), the equation \(P_y(\lambda) = 0\) has a unique root \(\lambda_1(y)\) in some small neighborhood of \(\gamma_1\). Similarly, there exists a neighborhood \(U_2\) such that there is a unique root \(\lambda_2(y)\) in a small neighborhood of \(\gamma_2\).

		Take two disjoint open neighborhoods \(V_1\) and \(V_2\) such that
		\[
		\gamma_1\in V_1 , \quad \gamma_2\in V_2.
		\]
		After shrinking \(U_1, U_2\), we can ensure
		\[
		\lambda_1(y) \in V_1, \quad \lambda_2(y) \in V_2.
		\]
		Let
		\[
		U' = U_1 \cap U_2\neq \emptyset.
		\]
		Then for any \(y \in U'\),
		\[
		\lambda_1(y) \neq \lambda_2(y),
		\]
		so \(y\) is a regular semisimple element.
		Define the spectral projections
		\[
		p_1(y) = \frac{y - \lambda_2(y) I}{\lambda_1(y) - \lambda_2(y)},\quad p_2(y) = \frac{y - \lambda_1(y) I}{\lambda_2(y) - \lambda_1(y)}.
		\]
		By the Cayley-Hamilton theorem, we have
		\begin{itemize}
			\item \(p_1(y) + p_2(y) = I\),
			\item \(p_i(y)^2 = p_i(y)\), \(i=1,2\),
			\item \(p_1(y) p_2(y) = 0\),
			\item \(y p_i(y) = \lambda_i(y) p_i(y)\), \(i=1,2\).
		\end{itemize}
		Therefore, \(p_i(y)\) is the projection onto the \(\lambda_i(y)\)-eigenspace.
		When \(y = \gamma\),
		\[
		p_1(\gamma) = \begin{pmatrix} 1 & 0 \\ 0 & 0 \end{pmatrix}, \quad p_2(\gamma) = \begin{pmatrix} 0 & 0 \\ 0 & 1 \end{pmatrix}.
		\]
		Let \(e_1, e_2\) be the standard basis vectors, and define
		\[
		v_1(y) = p_1(y) e_1, \quad v_2(y) = p_2(y) e_2.
		\]
		Therefore,
		\[
		v_1(\gamma) = e_1, \quad v_2(\gamma) = e_2.
		\]
		Define the matrix
		\[
		k(y) = (v_1(y), v_2(y)).
		\]
		Then
		\[
		k(\gamma) = I.
		\]
		The map \(y \mapsto k(y)\) is continuous and locally analytic. Since \(K = \mathrm{GL}(2,\mathcal{O}_F)\) is an open subgroup of \(G(F)\) containing the identity, we can further shrink \(U'\) to obtain an open neighborhood
		\[
		U_\gamma \subset U'
		\]
		such that
		\[
		k(y) \in K \quad \text{for all } y \in U_\gamma.
		\]
		In fact, for any preassigned integer \(N\), we can take \(U_\gamma\) small enough such that
		\[
		k(y) \in K(N) = I + \varpi^N M_2(\mathcal{O}).
		\]
		It is east to check that
		\[
		y v_i(y) = \lambda_i(y) v_i(y).
		\]
		Therefore,
		\[
		y k(y) = k(y) \begin{pmatrix} \lambda_1(y) & 0 \\ 0 & \lambda_2(y) \end{pmatrix}.
		\]
		Let
		\[
		t(y) = \begin{pmatrix} \lambda_1(y) & 0 \\ 0 & \lambda_2(y) \end{pmatrix}.
		\]
		Since \(\lambda_1(y) \neq \lambda_2(y)\), we have
		\[
		t(y) \in M_0(F)_{\mathrm{reg}}.
		\]
		Thus,
		\[
		y = k(y) t(y) k(y)^{-1}.
		\]
		Now we take \(y=g^{-1}mg\), \(m=\rm{diag}(m_1,m_2)\), thus for \(w_0\in \{1,w\}\), \(t(y)=w_0^{-1}\rm{diag}(m_1,m_2)w_0\). Therefore 
		\begin{align*}
			g^{-1}mg=k(y)w_0^{-1}mw_0k(y)^{-1},
		\end{align*}
		which means \[w_0k(y)^{-1}g^{-1}\in M_0(F),\]that is \[g\in M_0(F)K.\]
		Since \(E_{\rm{rel},n}\) is a invariant subspace under \(K\), we can write (\ref{weyl}) as (\ref{weyl2}).
	\end{proof}
	Now we write (\ref{weyl2}) as 
	\[\frac{1}{2}\int_{M_0(F)_{\rm{reg}}}|D_G(m)|^{\frac{1}{2}}\left[|D_G(m)|^{\frac{1}{2}}\int_{M_0(F)\backslash G(F)}f(g^{-1}mg)dg \right]\mathscr{I}_{P_0}(\xi_{0,n}^{-1}\otimes 1,m) dm,\]which is 
	\begin{align}\label{inverse}
		\frac{1}{2}\int_{M_0(F)_{\rm{reg}}}|D_G(m)|^{\frac{1}{2}}\mathscr{I}_{P_0}(\xi_{0,n}^{-1}\otimes 1,m) f_{M_0}(m)dm.
	\end{align}
	Set \(\mathscr{R}_{\rm{rel},M_0}(\xi_{0,n}^{-1}\otimes 1)|_{E_{\rm{rel},n}}=a_n \rm{Id}_{E_{\rm{rel},n}}\), for any \(n\in\mathbb{Z}_{\ge0}\), we have:
	\begin{lemma}\label{318}
		For any \(\gamma\in M_0(F)_{\rm{reg}}\), \(f\in C_c^\infty(U_\gamma)\) and \(N\in\mathbb{Z}_{\ge0}\),
		\begin{align*}
			\rm{tr} ( \mathscr{R}_{\rm{rel},M_0}(\xi_{0,N}^{-1}\otimes 1) \mathscr{I}_{P_0}(\xi_{0,N}^{-1}\otimes 1,f)|_{V_{\rm{rel},N}} )=\sum_{n\le N}a_n(F_n(\gamma)-F_{n-1}(\gamma)).
		\end{align*}
	\end{lemma}
	By Lemma \ref{315}, \[a_n = \begin{cases}
		\sqrt{2}(\rm{log}\,q) n, & c \ge 1, n \ge c, \\
		0, & c = 0, n = 0, \\
		\sqrt{2}(\rm{log}\,q) \left( S_q(z) + n \right), & c = 0, n \ge 1.
	\end{cases}\]
	It is easy to check that \[F_n(m)=F_n(w^{-1}mw),\quad \text{for}\ m=\rm{diag}(m_1,m_2),\]thus we can assume \[m = z_0h_u, \qquad h_u = \begin{pmatrix} u & 0 \\ 0 & 1 \end{pmatrix}.\]
	Set \[r=v(u)\ge 0,\quad d=v(u-1).\]
	\begin{lemma}
		For every \(n\) and \(c\),
		\begin{align*}
			F_n(h_u) = \begin{cases} 
				\xi_{0,n}(u)^{-1} q^{\min\{d,n\}} + q^{\min\{d,n-1\}}, & n \ge \max\{1, c\},r = 0, \\ 
				q^{-r/2} \left( 1 + \xi(u)^{-1} \right), &n \ge \max\{1, c\}, r > 0,\\
				1, &c=n=0, r = 0, \\ 
				\frac{q^{-r/2}}{1+q^{-1}} \left[ 1 + z^{-r} + (1-q^{-1}) \sum_{j=1}^{r-1} z^{-j} \right], &c=n=0, r > 0. 
			\end{cases}
		\end{align*}
	\end{lemma}
	\begin{proof}
		\textbf{Case 1: $r=0$.}
		In this case, \(u\in \mathcal{O}_F^\times\), \(h_u\in K\).
		
		Recall that 
		\[G_n/P_n\simeq\mathbf{P}^1(R_n) = \{ [1 ,t] : t \in R_n \} \sqcup \{ [\varpi s,1] : s \in \mathcal{O}_F/\varpi^{n-1}\mathcal{O}_F \}.\]
		Consider \(h_u[1,t] = [1,u^{-1}t],\) which is \[(u-1)t \equiv 0 \pmod{\varpi^n}.\]
		If \(d < n\), we have
		\[
		t \in \varpi^{n-d} R_n,
		\]
		which has \(q^d\) solutions. 
		
		If \(d \ge n\), all \(t \in R_n\) are solutions, which has \(q^n\) solutions. Therefore, the number of fixed points is
		\[
		q^{\rm{min}\{d,n\}}.
		\]
		For each coset $gP_n\in G_n/P_n$, choose a basis function $e_g$ supported on $P_ng^{-1}$ and normalized by $e_g(g^{-1})=1$.

		Take the representative element
		\[
		g_t = \begin{pmatrix} 1 &  \\ t & 1 \end{pmatrix}, \qquad g_t^{-1} = \begin{pmatrix} 1 &  \\ -t & 1 \end{pmatrix},
		\]we have 
		\[g_t^{-1}h_ug_t=\left(\begin{matrix}
			u&\\
			(1-u)t&1
		\end{matrix}\right)\equiv\left(\begin{matrix}
			u&\\
			&1
		\end{matrix}\right)\pmod{\varpi^n},\]which means
		\begin{align*}
			(\mathscr{I}_{P_0}(\xi^{-1}_{0,n}\otimes1,h_u)e_{t})(g_t^{-1})=e_t(g^{-1}_th_ug_t\cdot g_t^{-1})=\xi(u)^{-1}.
		\end{align*}
		Consider \(h_u[\varpi s, 1] =[u\varpi s : 1]\), which is 
		\[(u-1)s\equiv0 \pmod{\varpi^{n-1}}.\]
		The number of solutions is
		\[
		q^{\min\{d,n-1\}}.
		\]
		Taking the representative element
		\[
		g_s = \begin{pmatrix} \varpi s & 1 \\ 1 &  \end{pmatrix}.
		\]
		Then
		\[
		g_s^{-1} h_u g_s = \begin{pmatrix} 1 & 0 \\ (u-1)\varpi s & u \end{pmatrix}\equiv\left(\begin{matrix}
			1&\\
			&u
		\end{matrix}\right)\pmod{\varpi^{n}},
		\]that is 
		\begin{align*}
			(\mathscr{I}_{P_0}(\xi^{-1}_{0,n}\otimes1,h_u)e_{s})(g^{-1}_s)=e_s(g^{-1}_sh_ug_s\cdot g_s^{-1})=1.
		\end{align*}
		Therefore, \[F_n(h_u) =\xi(u)^{-1} q^{\min\{d,n\}} + q^{\min\{d,n-1\}}, \quad r = 0, \quad n \ge \max(1,c)).\]
		\textbf{Case 2: $r>0$.}

		Let
		\[
		\phi_{u,n} = \frac{\mathbf{1}_{h_u K(n)}}{\rm{vol}(K(n))}.
		\]
		The integral formula for the distribution of the ordinary character gives
		\[
		\rm{tr}(\mathscr{I}_{P_0}(\xi^{-1}_{0,n}\otimes 1,h_u)|_{V_{\rm{rel},n}}) = \frac{1}{\rm{vol}(K(n))} \int_{K(n)} \Theta_{\mathscr{I}_{P_0}(\xi^{-1}_{0,n}\otimes 1)}(h_u k) \, dk.
		\]
		Take arbitrarily
		\[
		k = \begin{pmatrix} 1 + \varpi^n A & \varpi^n B \\ \varpi^n E & 1 + \varpi^n D \end{pmatrix} \in K(n).
		\]
		The characteristic polynomial of the matrix \(h_u k\) is
		\[
		p_k(X) = X^2 - (1 + \varpi^n D + u(1 + \varpi^n A))X + u \det k.
		\]
		Since \(r > 0\), we have \(u \in \varpi \mathcal{O}_F\). Furthermore,
		\[
		p_k(1) = \varpi^n D(u-1) + u \varpi^{2n} (AD - BE) \in \varpi^n \mathcal{O}_F,
		\]
		and
		\[
		p_k'(1) = 1 - u - \varpi^n (D + uA) \in \mathcal{O}_F^\times.
		\]
		By Hensel's lemma, there exists an eigenvalue
		\[
		\lambda_1 \in 1 + \varpi^n \mathcal{O}_F.
		\]
		The other eigenvalue is
		\[
		\lambda_2 = \frac{u \det k}{\lambda_1}.
		\]
		Since \(\det k\) and \(\lambda_1\) both belong to \(1 + \varpi^n \mathcal{O}_F\), we obtain
		\[
		\lambda_2 \in u(1 + \varpi^n \mathcal{O}_F),
		\]
		so
		\[
		v(\lambda_1) = 0, \quad v(\lambda_2) = r, \quad |\lambda_1 - \lambda_2| = 1.
		\]
		
		In particular, \(h_u k\) is a split regular semisimple element, and its discriminant satisfies
		\[
		|D_G(h_u k)|^{1/2} = \frac{|\lambda_1 - \lambda_2|}{|\lambda_1 \lambda_2|^{1/2}} = q^{r/2}.
		\]
		
		Now using the character formula of the ordinary normalized principal series. The contributions of the two fixed flags are respectively
		\[
		\frac{|\lambda_1 / \lambda_2|^{1/2} \xi(\lambda_1)^{-1}}{|1 - \lambda_1 / \lambda_2|}\quad \text{and}\quad \frac{|\lambda_2 / \lambda_1|^{1/2} \xi(\lambda_2)^{-1}}{|1 - \lambda_2 / \lambda_1|}.
		\]
		The absolute value factors in both fractions are equal to
		\[
		\frac{|\lambda_1 \lambda_2|^{1/2}}{|\lambda_1 - \lambda_2|} = |D_G(h_u k)|^{-1/2}.
		\]
		Therefore
		\[
		\Theta_{\mathscr{I}_{P_0}(\xi^{-1}_{0,n}\otimes 1)}(h_u k) = \frac{\xi(\lambda_1)^{-1} + \xi(\lambda_2)^{-1}}{|D_G(h_u k)|^{1/2}}.
		\]
		When \(n \ge c\), \(\xi\) is trivial on \(1 + \varpi^n \mathcal{O}_F\). Therefore, we have
		\[
		\xi(\lambda_1) = 1, \quad \xi(\lambda_2) = \xi(u).
		\]
		We can obtain
		\[
		\Theta_{\mathscr{I}_{P_0}(\xi^{-1}_{0,n}\otimes 1)}(h_u k) = q^{-r/2} (1 + \xi(u)^{-1}) \quad \text{for all } k \in K(n).
		\]
		Therefore,
		\[
		\rm{tr}(\mathscr{I}_{P_0}(\xi_{0,n}\otimes 1)(h_u)|_{V_{\rm{rel},n}}) = q^{-r/2} (1 + \xi(u)^{-1}), \quad r > 0, \quad n \ge \max\{1,c\}.
		\]
		\textbf{Case 3: $c=0$, \(n=0\).}
		In this case, \(\xi\) is trivial on \(\mathcal O_F^\times\) and \(\rm{dim}\,V_{\rm{rel},0}=1\).

		Take the spherical vector \(e_0\), satisfying
		\[
		e_0(k) = 1 \quad \text{for}\ k \in K,
		\]
		and normalize \(\operatorname{vol}(K) = 1\). Under the inner product of the compact model, \(\|e_0\| = 1\), so
		\[
		\rm{tr}(\mathscr{I}_{P_0}(\xi^{-1}\otimes 1,h_u)|_{V_{\rm{rel},0}}) = \langle \mathscr{I}_{P_0}(\xi^{-1}\otimes 1,h_u)e_0, e_0 \rangle = \int_K e_0(kh_u) \, dk. 
		\]
		If \(r = 0\), then \(h_u \in K\), hence
		\[
		e_0(kh_u) = 1.
		\]
		Therefore
		\[
		F_0(h_u) = 1, \quad r = 0.
		\]
		Assume \(r > 0\), and write
		\[
		u = \varpi^r u_0, \quad u_0 \in \mathcal{O}_F^\times,
		\]
		and let
		\[
		a_r = \operatorname{diag}(\varpi^r,1).
		\]
		Since
		\[
		h_u = a_r \operatorname{diag}(u_0,1)
		\]
		and \(e_0\) is right \(K\)-invariant, we have
		\[
		\operatorname{tr}(\mathscr{I}_{P_0}(\xi^{-1}\otimes 1,h_u)|_{E_{\rm{rel},0}}) = \int_K e_0(ka_r) \, dk.
		\]
		Let the bottom row of \(k \in K\) be \((x,y)\). It is a primitive row vector, therefore
		\[
		\max\{|x|, |y|\} = 1.
		\]
		Let
		\[
		j = \min\{v(y), r\},
		\]
		where we define \(v(0) = +\infty\).
		The bottom row of the matrix \(ka_r\) is
		\[
		(\varpi^rx,  y).
		\]
		Its maximum norm is \(q^{-j}\): if \(y\) is a unit, the maximum norm is 1, if \(y\) is not a unit, then \(x\) is a unit, the maximum norm is
		\[
		\max\{q^{-r}, q^{-v(y)}\} = q^{-\min\{r, -v(y)\}}.
		\]
		Perform the Iwasawa decomposition
		\[
		ka_r = \begin{pmatrix} a & * \\ 0 & d_0 \end{pmatrix} k', \quad k' \in K.
		\]
		Since the bottom row of \(k'\) is a primitive row vector, the maximum norm of the bottom row on the right-hand side is \(|d_0|\), hence
		\[
		v(d_0) = j.
		\]
		Also since
		\[
		v(\det(ka_r)) = r,
		\]
		and \(\det k'\) is a unit, we have
		\[
		v(a) = r - j.
		\]
		Therefore
		\[
		e_0(ka_r) = \left| \frac{a}{d_0} \right|^{\frac{1}{2}} \xi(a)^{-1}.
		\]
		Since \(\xi\) is unramified,
		\[
		\xi(a) = z^{v(a)} = z^{r-j},
		\]
		that is
		\[
		e_0(ka_r) = q^{-r/2 + j} z^{-(r-j)}.
		\]
		It remains only to compute the distribution of \(j\).

		Take an integer \(L > r\). The Haar measure pushed forward to \(\mathrm{GL}(2,\mathcal{O}_F/\varpi^L\mathcal{O}_F)\) is the uniform measure, and its bottom row is uniformly distributed among all primitive row vectors. The total number of primitive row vectors is
		\[
		q^{2L} - q^{2L-2}.
		\]
		When \(j = 0\), \(y\in\mathcal{O}_F^\times\), \(x\) is arbitrary, so
		\[
		\operatorname{Pr}(j=0) = \frac{(q^L - q^{L-1})q^L}{q^{2L} - q^{2L-2}} = \frac{1}{1 + q^{-1}}.
		\]
		For \(1 \le j < r\), we require \(v(y) = j\), and \(x\in\mathcal{O}_F^\times\). Therefore
		\[
		\operatorname{Pr}(j) = \frac{(q^{L-j} - q^{L-j-1})(q^L - q^{L-1})}{q^{2L} - q^{2L-2}} = \frac{(1 - q^{-1})q^{-j}}{1 + q^{-1}}. 
		\]
		Finally, \(j = r\) means \(v(y) \ge r\), and \(x\in\mathcal{O}_F^\times\). Therefore
		\[
		\operatorname{Pr}(j=r) = \frac{q^{L-r}(q^L - q^{L-1})}{q^{2L} - q^{2L-2}} = \frac{q^{-r}}{1 + q^{-1}}.
		\]
		Thus, we have
		\begin{align*}
			\int_K e_0(ka_r) \, dk &= \frac{q^{-r/2}z^{-r}}{1+q^{-1}} + \sum_{j=1}^{r-1} \frac{(1-q^{-1})q^{-j}}{1+q^{-1}} q^{-r/2+j} z^{-(r-j)} + \frac{q^{-r}}{1+q^{-1}} q^{r/2}\\
			&=\frac{q^{-r/2}}{1+q^{-1}} \left[ z^{-r} + (1-q^{-1}) \sum_{j=1}^{r-1} z^{-(r-j)} + 1 \right]\\
			&=\frac{q^{-r/2}}{1+q^{-1}} \left[ 1 + z^{-r} + (1-q^{-1}) \sum_{j=1}^{r-1} z^{-j} \right].
		\end{align*}
		Finally, we have
		\[
		F_0(h_u) = \frac{q^{-r/2}}{1+q^{-1}} \left[ 1 + z^{-r} + (1-q^{-1}) \sum_{j=1}^{r-1} z^{-j} \right], \quad c=0, \ r>0.
		\]
	\end{proof}
	Then, the next lemma is clear.
	\begin{lemma}\label{320}
		If \(r=0\),
		\[
		F_n(h_u)-F_{n-1}(h_u) = \begin{cases} 
			q^n \xi(u)^{-1} + q^{n-1}, & n = c \ge 1, d \ge n, \\
			q^d(\xi(u)^{-1} + 1), & n = c \ge 1, d \le n - 1, \\
			1, & n = 1, c = 0, d = 0, \\
			q, & n = 1, c = 0, d \ge 1, \\
			q^{n-2}(q^2 - 1), & n \ge \max\{2, c + 1\}, d \ge n, \\
			q^{n-2}(q - 1), & n \ge \max\{2, c + 1\}, d = n - 1, \\
			0, & n \ge \max\{2, c + 1\}, d \le n - 2.
		\end{cases}
		\]
		If \(r > 0\),
		\[
		F_n(h_u)-F_{n-1}(h_u) = \begin{cases} 
			q^{-r/2}(1 + \xi(u)^{-1}), & n = c \ge 1, \\
			\frac{q^{-r/2}}{1+q^{-1}}
			\left[
			q^{-1}(1+z^{-r})
			-(1-q^{-1})\sum_{j=1}^{r-1}z^{-j}
			\right]., & n = 1, c = 0, \\
			0, & n \ge \max\{2, c + 1\}.
		\end{cases}
		\]
	\end{lemma}
	By (\ref{inverse}) and Lemma \ref{318}, Lemma \ref{320}, we have
	\begin{lemma}\label{fourier}
		For \(\chi=\chi_1\otimes \chi_2\), \(m=z_0h_u\in M_0(F)_{\rm{reg}}\), the Fourier kernel \(I_{M_0,M_0}(m,\chi)=I_{M_0,M_0}(\chi^\vee,m)\) equals 
		\begin{equation*}
			\sqrt{2} (\log q) \cdot \chi_2(z_0)^{-1} \chi_1(z_0u)^{-1} \times \begin{cases} 
				(1 - q^{-1}) \frac{q^{-1} \sum_{j=0}^r z^j - \sum_{j=1}^{r-1} z^j}{(1 - q^{-1} z)(1 - q^{-1} z^{-1})}, & c = 0, r > 0, \\[1.5ex]
				c(1 + \xi_0(u_0)z^r), & c \ge 1, r > 0, \ u = \varpi^r u_0, \\[1.5ex]
				c(1 + \xi_0(u)), & c \ge 1, r = 0, d < c, \\[1.5ex]
				2(d + 1) - \frac{q + 1}{q - 1} \left( 1 - q^{-(d - c + 1)} \right), & c \ge 1, r = 0, d \ge c, \\[1.5ex]
				(2 - q^{-d}) S_q(z) + 2d + 1 - \frac{2(1 - q^{-d})}{q - 1}, & c = 0, r = 0.
			\end{cases}
		\end{equation*}
	\end{lemma}
	Now, for \(\gamma\in M_{0}(F)_{\rm{reg}}\) we can write \(I_{M_0}(\gamma,f)\) as
	\begin{align*}
		\frac{1}{2}\int_{T_{\rm{disc}}(M_0)}I_{M_0,M_0}(\gamma,\tau)f_M(\tau)d\tau+\int_{T_{\rm{disc}}(G)}I_{M_0,G}(\gamma,\tau)f_G(\tau)d\tau,\quad \text{for}\ \gamma\in M_0(F)_{\rm{reg}}.
	\end{align*}
    Also,
	\begin{align}\label{f12}
		I_{G,M_0}(\gamma,\chi)=I_{G,M_0}(\chi^{-1}, \gamma) = |D_G(\gamma)|^{\frac{1}{2}} \Theta_{\chi}( \gamma^{-1}),
	\end{align}
	by Theorem \ref{reciprocity}.

	By the descent formula in \cite{Art94},
	\begin{align}\label{GG}
		I_{G,G}(\gamma, \tau) = \sum_{S \in \mathcal{L}(M_0)} d_{M_0}^G(G, S) \sum_{\substack{w \in W_0^S \backslash W_0^G \\ wG \subseteq S}} I_{M_0, wG}^S (\gamma, w\tau),
	\end{align}
	we have \[I_{G,G}(\gamma,\tau_\pi)=0.\]
    
    \subsubsection{The case \(\gamma\in G_{\rm{ell,reg}}\)}
    In this subsubsection, we calculate the Fourier kernel \(I_G(\gamma,\tau)\), for the elliptic regular terms \(\gamma \) in \(\GGL(2,F)\).
    \begin{lemma}
        For \(h\in C_c^\infty(G_{\rm{ell,reg}}(F))\),\[\rm{tr}(\mathscr{I}_{P_0}(\chi^\vee,h))=0,\]where \(\chi=\chi_1\otimes\chi_2\), \(\chi_1,\chi_2\in\widehat{F^\times}\)
    \end{lemma}
    \begin{proof}
By the definition of trace \[\operatorname{tr}\mathscr I_{P_0}(\chi^\vee,h)
=
\int_{M_0(F)}
\delta_{P_0}(m)^{1/2}\chi^\vee(m)
\int_K\int_{N_0(F)}
h(k^{-1}mnk)\,dn\,dk\,dm.
\]
For \(m=\operatorname{diag}(a,b)\) and \(n\in N_0(F)\), the characteristic polynomial of \(mn\) is
\[
\det(XI_2-mn)=(X-a)(X-b).
\]
Hence \(k^{-1}mnk\) is not regular elliptic. The support condition on
\(h\) therefore implies
\[
h(k^{-1}mnk)=0
\qquad
(m\in M_0(F),\ n\in N_0(F),\ k\in K),
\]
which means
\[
\operatorname{tr}\mathscr I_{P_0}(\chi^\vee,h)=0.
\]
\end{proof}
Therefore, like Lemma \ref{local0}, we have
\begin{lemma}\label{simil}
    For any \(\gamma\in G_{\rm{ell,reg}}(F)\), and \(h\in C_c^\infty(G_{\rm{ell,reg}}(F))\), \(J_{G}(\chi^\vee,h)=0\).
\end{lemma}
\begin{lemma}
    For \(\gamma\in G_{\rm{ell,reg}}(F)\), the Fourier kernel \(I_{G,M_0}^G(\gamma,\chi)=-I_{M_0,G}^G(\chi^\vee,\gamma)
\) equals 
    \begin{align*}
\frac{\sqrt2\log q}{1-q^{-1}}
|\Delta|^{1/2}\chi_1(d)^{-1}
\int_F\frac{\xi(t)}{|t^2-Tt+d|}\,dt.
    \end{align*}
\end{lemma}
\begin{proof}
    By Lemma \ref{simil}, \[I_{M_0}^G(\chi^\vee,h)=J_{M_0}^G(\chi^\vee,h)=\sqrt{2}\rm{tr}(A'(0)\mathscr{I}_{P_0}(\chi^\vee,h)).\]
    According to (\ref{unrami}) and (\ref{rami}), we have 
    \[\log|u|+\begin{cases} 
    (\log q)S_q(z),&\xi\text{ unramified},\\
    0,&\xi\text{ ramified}.\end{cases}\]
    But \[\rm{tr}((\log q)S_q(z)\mathscr{I}_{P_0}(\chi^\vee,h))=(\log q)S_q(z)\rm{tr}(\mathscr{I}_{P_0}(\chi^\vee,h))=0.\]
    Set \(g_s(x) = |x|^{s-1}.\) Note that we can let \(\xi\) equals trivial character in the equalty (\ref{fourier1}) to get
\begin{equation*}
\mathcal{F}g_s(u) = c_s(1)|u|^{-s}, \qquad c_s(1) = \frac{1 - q^{s-1}}{1 - q^{-s}}.
\end{equation*}

Therefore
\begin{equation*}
\begin{aligned}
\mathcal{F}^{-1}(|\cdot|^{-s})(x) &= c_s(1)^{-1} g_s(x) \\
&= \frac{1 - q^{-s}}{1 - q^{s-1}}|x|^{s-1}.
\end{aligned}
\end{equation*}Then, 
\[
\begin{aligned}
\mathcal F^{-1}(\log|\cdot|)(x)
&=-\left.\frac{d}{ds}
\left(\frac{1-q^{-s}}{1-q^{s-1}}|x|^{s-1}\right)\right|_{s=0}\\
&=-\frac{\log q}{1-q^{-1}}\,|x|^{-1},
\qquad x\neq0.
\end{aligned}
\]
Choose the conjugate representative
\[\gamma=\left(\begin{matrix}0&-d\\1&T\end{matrix}\right),\qquad T=\operatorname{tr}\gamma,\quad d=\det\gamma,\quad
\Delta=T^2-4d.\]
Since 
\[
\begin{aligned}
wn(x)\gamma=
\begin{pmatrix}d/x&-1\\0&x\end{pmatrix}
wn\!\left(T-\frac d x\right),
\end{aligned}
\]by the definition of induced representation, 
\[
\bigl(\mathscr I_{P_0}(\chi^\vee,\gamma)\Phi\bigr)(x)
=
\chi_1(d)^{-1}|d|^{1/2}
\frac{\xi(x)}{|x|}
\Phi\!\left(T-\frac d x\right).
\]Take the trace, we have 
\[
\begin{aligned}
I_{M_0,G}^G(\chi^\vee,\gamma)
&=
-\frac{\sqrt2\log q}{1-q^{-1}}
|D_G(\gamma)|^{1/2}|d|^{1/2}\chi_1(d)^{-1}
\int_F\frac{\xi(x)}{|x|\,|T-d/x-x|}dx\\
&=
-\frac{\sqrt2\log q}{1-q^{-1}}
|\Delta|^{1/2}\chi_1(d)^{-1}
\int_F\frac{\xi(x)}{|x^2-Tx+d|}\,dx,
\end{aligned}
\]since the support of \(h\) is the subset of \(G_{\rm{ell,reg}}(F)\).
\end{proof}
Also,
	\begin{align}
		I_{G,G}(\gamma,\tau)=I_{G,G}(\tau^{-1}, \gamma) = |D_G(\gamma)|^{\frac{1}{2}} \Theta_{\tau}( \gamma^{-1}),
	\end{align}
	by Theorem \ref{reciprocity}.

	\section{Fourier kernel for all elements}
	\subsection{The non-Archimedean case}
	In this subsection, we write the Fourier kernel \(I_{M_0}(m,\chi)\) and \(I_{M_0}(m,\tau)\) for all \(m\in M_0(F)\), mainly by the method that Arthur introduced in \cite{Art88}.

	Recall that
	\[
	m = z_0 \begin{pmatrix} u & 0 \\ 0 & 1 \end{pmatrix}, \quad u \neq 1, \quad d = v_F(u-1),\quad \text{for}\ m\in M_0(F)_{\rm{reg}},
	\]
	and denote
	\[
	L_q = \sqrt{2} \log q, \quad e = v(2).
	\]
	Let
	\[
	S_q(z) = \frac{q^{-1}z}{1 - q^{-1}z} + \frac{q^{-1}z^{-1}}{1 - q^{-1}z^{-1}}.
	\]
	Arthur introduces the definition of weighted integrals of \(\gamma\) which is not regular and semisimple in \cite{Art88} by
	\[J_M(\gamma, f) = \lim_{a \to 1} \sum_{L \in \mathcal{L}(M)} r_M^L(\gamma, a) J_L(a\gamma, f),\]where he defines the function \(r_M^L(\gamma,a)\) to cancel the poles. By Cor \ref{32}, we can write 
	\[
	I_{M_0,L}(z_0I_2, \tau) = \lim_{u \to 1} \left[ I_{M_0,L}(z_0u, \tau) + \rho(u) I_{G,L}(z_0u, \tau) \right], \quad L = M_0, G,
	\]where
	\[
	\rho(u)=r_{M_0}^G(z_0h_u)= \sqrt{2} \log |u - u^{-1}|.
	\]
	Since we shall let \(u\rightarrow 1\), then under the non-Archimedean topology,
	\[
	u \to 1 \iff d = v(u-1) \to +\infty,
	\]and in this case, \(r=v(u)\rightarrow0\), which means we do not have to consider the case \(r>0\). Also, \[\chi_1(u)=\chi_2(u)=\xi(u)\rightarrow1,\quad \text{if}\ u\rightarrow 1.\]
	Therefore, we only need to consider the cases \(c=0\) and \(c\ge1\), \(d\ge c\) in Lemma \ref{fourier}.

	\textbf{Case 1: $L=M_0$.}

	Consider the orbital integral \(I_G(z_0h_u,f)\), its Fourier kernel is 
	\[B_\chi(u) := I_{G,M_0}(z_0h_u, \chi) = \chi_1(z_0u)^{-1}\chi_2(z_0)^{-1} + \chi_1(z_0)^{-1}\chi_2(z_0u)^{-1}.\]
	When \(u\) is sufficiently close to 1,
	\[
	B_\chi(u) = 2\omega_\chi(z_0)^{-1}.
	\]
	On the other hand, since
	\[
	u - u^{-1} = \frac{(u-1)(u+1)}{u},
	\]
	when \(d > e\),
	\[
	v(u+1) = v(2 + (u-1)) = e.
	\]
	Thus
	\[
	\rho(u) = -\sqrt{2}(d+e)\log q = -L_q(d+e).
	\]
	Therefore, 
	\[
	\rho(u)B_\chi(u) = -2L_q\omega_\chi(z_0)^{-1}(d+e).
	\]

	First compute \(c=0\). By Lemma \ref{fourier}, when \(u\) is sufficiently close to 1,
	\[
	I_{M_0,M_0}(z_0u, \chi) = L_q \omega_\chi(z_0)^{-1} \left[ (2 - q^{-d})S_q(z) + 2d + 1 - \frac{2(1 - q^{-d})}{q-1} \right].
	\]
	Write
	\[
	(2 - q^{-d})S_q(z) + 2d + 1 - \frac{2(1-q^{-d})}{q-1}
	\]
	\[
	= 2d + 1 + 2S_q(z) - \frac{2}{q-1} + q^{-d} \left[ \frac{2}{q-1} - S_q(z) \right].
	\]
	Therefore,
	\[
	I_{M_0,M_0}(\gamma_u, \chi) + \rho(u)B_\chi(u)
	\]
	\[
	= L_q \omega_\chi(z_0)^{-1} \left[ 1 + 2S_q(z) - \frac{2}{q-1} - 2e + q^{-d} \left( \frac{2}{q-1} - S_q(z) \right) \right].
	\]
	Let \(d \to \infty\), the last term tends to zero, so
	\[
	I_{M_0,M_0}(z_0, \chi) = L_q \omega_\chi(z_0)^{-1} \left[ \frac{q-3}{q-1} + 2S_q(z) - 2e \right], \quad c=0.
	\]
	Next compute \(c \ge 1\). For \(d\ge c\),
	\[
	I_{M_0,M_0}(z_0h_u, \chi) = L_q \omega_\chi(z_0)^{-1} \left[ 2(d+1) - \frac{q+1}{q-1} \left( 1 - q^{-(d-c+1)} \right) \right].
	\]
	Write
	\begin{align*}
		2(d+1) - \frac{q+1}{q-1} \left( 1 - q^{-(d-c+1)} \right)=2d + 2 - \frac{q+1}{q-1} + \frac{q+1}{q-1} q^{-(d-c+1)}.
	\end{align*}
	Since
	\[
	2 - \frac{q+1}{q-1} = \frac{q-3}{q-1},
	\]
	we have,
	\begin{align*}
		&I_{M_0,M_0}(z_0h_u, \chi) + \rho(u)B_\chi(u)\\
		=&L_q \omega_\chi(z_0)^{-1} \left[ \frac{q-3}{q-1} - 2e + \frac{q+1}{q-1} q^{-(d-c+1)} \right].
	\end{align*}
	Fixing \(c\), letting \(d \to \infty\), we obtain
	\[
	I_{M_0,M_0}(z_0I_2, \chi) = L_q \omega_\chi(z_0)^{-1} \left[ \frac{q-3}{q-1} - 2e \right], \quad c \ge 1.
	\]
	The conductor \(c\) only appears in the remainder term that tends to zero. Therefore, the final value in the ramified case no longer depends on the relative conductor \(c\), but still depends on \(\chi\) through the central character.
	
	Combining the two cases, this is the formula for all unitary principal series parameters and all scalar matrices:
	\begin{align*}
		I_{M_0,M_0}(z_0 I_2, \chi) = \sqrt{2}(\log q) \omega_\chi(z_0)^{-1} \times\begin{cases} \frac{q-3}{q-1} + 2S_q(z) - 2v(2), & c=0, \\ \frac{q-3}{q-1} - 2v(2), & c \ge 1. \end{cases}
	\end{align*}
	\textbf{Case 2: \(L=G\).}

	Recall that in this case, \[
	I_{M_0,G}(z_0I_2, \tau) = \lim_{u \to 1} \left( I_{M_0,G}(z_0h_u, \tau) + \rho(u) I_{G,M_0}(z_0h_u, \tau) \right), \quad \tau \in T_{\mathrm{disc}}(G).
	\]By (\ref{GG}),
	\[I_{G,G}(z_0h_u,\tau_\pi)=0,\]
	which means \[I_{M_0,G}(z_0I_2, \tau) =\lim_{u \to 1} I_{M_0,G}(z_0h_u, \tau).\]
	By Theorem \ref{reciprocity}, 
	\begin{align*}
		I_{M_0,G}(z_0h_u, \tau) =& (-1)^{\dim(\mathfrak{a}_{M_0} / \mathfrak{a}_G)} I_{G,G}(\tau^\vee, z_0h_u)\\
		=&-|D_G(z_0h_u)|^{1/2} \Theta_{\tau^\vee}(z_0h_u).
	\end{align*}
	Thus,
	\[I_{M_0,G}(z_0I_2, \tau) =-\rm{lim}_{u\to1}q^{-d}\Theta_{\tau^\vee}(z_0h_u) ,\]where \[\Theta_{\tau^\vee}(z_0h_u) = \omega_{\tau^\vee}(z_0) \Theta_{\tau^\vee} \begin{pmatrix} u & 0 \\ 0 & 1 \end{pmatrix} = \omega_{\tau}(z_0)^{-1} \Theta_{\tau^\vee} \begin{pmatrix} u & 0 \\ 0 & 1 \end{pmatrix}.\]
	Let \(\pi \in \Pi_2(G)\).
	When \(u\) is sufficiently close to 1, there exists
	\[
	X_u = \begin{pmatrix} \log u & 0 \\ 0 & 0 \end{pmatrix}
	\]
	such that
	\[
	h_u = \exp X_u.
	\]
	The Lie algebra \(\mathfrak{gl}(2,F)\) has two nilpotent orbits:
	\[
	\mathcal{O}_0 = \{0\}, \quad \mathcal{O}_{\mathrm{reg}}.
	\]
	By \cite{Howe}, for irreducible admissible representation \(\rho\), in a sufficiently small neighborhood of the identity element, we have
	\[\Theta_\rho(\exp X) = \sum_{\mathcal{O}} c_{\mathcal{O}}(\rho) \widehat{\mu}_{\mathcal{O}}(X),\]then in our case,
	\[
	\Theta_\pi(\exp X) = c_0(\pi^\vee) + c_{\mathrm{reg}}(\pi^\vee) \widehat{\mu}_{\mathrm{reg}}(X).
	\]
	Here, the Fourier transform of the zero orbit has been normalized to
	\[
	\widehat{\mu}_0(X) = 1.
	\]
	Next, we shall consider the Whittaker model. The representation \(\pi\) of \(\GGL(2,F)\) is generic, then \(\pi^\vee\) is. The existence and uniqueness of the Whittaker model of \(\GGL(n)\) have been proved in \cite{shal}, that is
	\[
	\dim \mathrm{Wh}_\pi(\pi^\vee) = 1.
	\]
	If residual characteristic of \(F\) is \(2\), by \cite{MW87}, we have
	\[c_{\mathrm{reg}}(\pi^\vee) = \dim \mathrm{Wh}_\pi(\pi^\vee) = 1.\]And if residual characteristic of \(F\) is \(2\), by \cite{Var14}, we have the same identity.

	Thus,
	\[
	\Theta_{\pi^\vee}(h_u) = c_0(\pi^\vee) + \widehat{\mu}_{\mathrm{reg}}(X_u).
	\]
	Now, take the normalized induced representation
	\[
	\pi_0 = \rm{Ind}_{P_0}^G(1 \otimes 1),
	\]where \(1\) means the trivial representation of \(F^\times\).

	Since
	\[
	1 \neq |\cdot|^{-1},
	\]
	this principal series is irreducible and generic, hence
	\[
	c_{\mathrm{reg}}(\pi_0) = 1.
	\]
	For \(h_u = \operatorname{diag}(u, 1)\), the normalized principal series character formula gives
	\[
	\Theta_{\pi_0}(h_u) = |D_G(h_u)|^{-1/2}(1+1).
	\]
	Since
	\[
	|D_G(h_u)|^{-1/2} = q^{d},
	\]
	we have
	\[
	\Theta_{\pi_0}(h_u) = 2q^d.
	\]
	On the other hand, the local expansion gives
	\[
	\Theta_{\pi_0}(h_u) = c_0(\pi_0) + \widehat{\mu}_{\mathrm{reg}}(X_u).
	\]
	Comparing the two formulas, we obtain
	\[
	\widehat{\mu}_{\mathrm{reg}}(X_u) = 2q^d - c_0(\pi_0).
	\]
	In particular,
	\[
	\widehat{\mu}_{\mathrm{reg}}(X_u) = 2q^d + O(1).
	\]
	Since
	\[
	c_{\mathrm{reg}}(\pi^\vee) = 1,
	\]
	we have
	\[
	\begin{aligned}
		\Theta_\pi(h_u) &= c_0(\pi^\vee) + \widehat{\mu}_{\mathrm{reg}}(X_u) \\
		&= c_0(\pi^\vee) + 2q^d - c_0(\pi_0) \\
		&= 2q^d + (c_0(\pi^\vee) - c_0(\pi_0)).
	\end{aligned}
	\]
	Therefore,
	\[
	\Theta_\pi(h_u) = 2q^d + O(1).
	\]
	Now we have
	\[
	\begin{aligned}
		\Theta_{\pi^\vee}(z_0h_u) &= \Theta_{\pi^\vee}(z_0 h_u) \\
		&= \omega_{\pi^\vee}(z_0) \Theta_{\pi^\vee}(h_u) \\
		&= \omega_\pi(z_0)^{-1} (2q^d + O(1)).
	\end{aligned}
	\]
	Thus,
	\[
	\begin{aligned}
		\lim_{u \to 1} q^{-d} \Theta_{\pi^\vee}(z_0h_u) &= \omega_\pi(z_0)^{-1} \lim_{u \to 1} q^{-d} (2q^d + O(1)) \\
		&= 2\omega_\pi(z_0)^{-1}.
	\end{aligned}
	\]
	Therefore,
	\[
	\begin{aligned}
		I_{M_0,G}(z_0I_2, \tau_\pi) &= \lim_{u \to 1} (I_{M_0,G}(z_0h_u, \tau_\pi) + \rho(u) I_{G,G}(z_0h_u, \tau_\pi)) \\
		&= -2\omega_\pi(z_0)^{-1}.
	\end{aligned}
	\]
	Now, let
	\[
	\tau_\eta = (M_0, \eta \otimes \eta,1),
	\]
	where \(\eta\) is a character of \(F^\times\), and \(\eta\otimes \eta\) is due to the definition of \(T_{\rm{disc}}\). The corresponding normalized induced representation is
	\[
	\pi_\eta = \rm{Ind}_{P_0}^G(\eta \otimes \eta).
	\]
	Its dual representation is
	\[
	\pi_\eta^\vee = \rm{Ind}_{P_0}^G(\eta^{-1} \otimes \eta^{-1}).
	\]
	The normalized principal series character formula gives
	\[
	\Theta_{\pi^\vee_\eta}(z_0h_u) = |D^G(h_u)|^{-1/2} (\eta^{-1}(z_0 u)\eta^{-1}(z_0) + \eta^{-1}(z_0)\eta^{-1}(z_0 u)).
	\]
	Therefore
	\[
	\Theta_{\pi^\vee_\eta}(z_0h_u) = 2q^d \eta(z_0)^{-2}\eta(u)^{-1}.
	\]
	There exists an integer \(N \ge 1\) such that
	\[
	\eta(1+\mathfrak{p}^N\mathcal{O}_F) = 1.
	\]
	When \(u\) is sufficiently close to 1,
	\[
	\eta(u) = 1.
	\]
	Thus,
	\[
	\Theta_{\pi^\vee_\eta}(z_0h_u) = 2q^d \eta(z_0)^{-2},
	\]
	so
	\[
	\lim_{u \to 1} q^{-d} \Theta_{\pi^\vee_\eta}(z_0h_u) = 2\eta(z_0)^{-2}.
	\]
	Let
	\[
	\sigma_\eta = \eta \otimes \eta.
	\]
	Its Weyl group is
	\[
	W_{\sigma_\eta} = \{1, w\}.
	\]
	Arthur's decomposition yields
	\[
	W_{\sigma_\eta}^0 = \{1, w\}, \quad R_{\sigma_\eta} = \{1\}.
	\]
	On the one-dimensional vector space
	\[
	\mathfrak{a}_{M_0}^G = \mathfrak{a}_{M_0} / \mathfrak{a}_G,
	\]
	the action of the reflection \(w\) is multiplication by \(-1\). Therefore,
	\[
	\epsilon(w) = -1,\qquad |\det(1-w)|_{\mathfrak{a}_{M_0}^G}| = 2.
	\]
	The identity fixes all of $\mathfrak a_{M_0}^G$ and is therefore not regular. Only \(w\) contributes to \(i^G(\tau_\eta)\). Thus,
	\[
	\begin{aligned}
		i^G(\tau_\eta) &= \frac{1}{|W_{\sigma_\eta}^0|} \epsilon(w) |\det(1-w)|_{\mathfrak{a}_{M_0}^G}|^{-1} \\
		&= \frac{1}{2}\times (-1) \times\frac{1}{2} \\
		&= -\frac{1}{4}.
	\end{aligned}
	\]
	Therefore,
	\begin{align*}
		I_{M_0,G}(z_0I_2, \tau_\eta) &= -i^G(\tau_\eta) \lim_{u \to 1} q^{-d} \Theta_{\pi_\eta}(z_0h_u) \\
		&= -\left(-\frac{1}{4}\right) 2\eta(z_0)^{-2} \\
		&= \frac{1}{2}\eta(z_0)^{-2}.
	\end{align*}
	Now we have
	\begin{lemma}
		For non-Archimedean local field \(F\), \(z_0I_2\in A_G(F)\), the Fourier kernel associated to \(M_0\) is 
		\begin{align*}
			I_{M_0,M_0}(z_0I_2,\chi)=\sqrt{2}(\log q) \omega_\chi(z_0)^{-1} \times\begin{cases} \frac{q-3}{q-1} + 2S_q(z) - 2v(2), & c=0, \\ \frac{q-3}{q-1} - 2v(2), & c \ge 1, \end{cases}
		\end{align*}for \(\chi\in \Pi_2(M_0(F))\) and 
		\begin{align*}
			I_{M_0,G}(z_0I_2, \tau_\pi) = -2\omega_\pi(z_0)^{-1},
		\end{align*}for \(\pi\in \Pi_2(G(F))\),
		\begin{align*}
			I_{M_0,G}(z_0I_2, \tau_\eta) =\frac{1}{2}\eta(z_0)^{-2},
		\end{align*}for \(\eta\in \Pi_2(\GGL(1,F))\). 
	\end{lemma}
	
	\subsection{The Archimedean Case}
	In this subsection, we use the notation of Hoffmann. Let \(F = \mathbb{R}\) or \(\mathbb{C}\), \(M = \mathrm{GL}(2,F)\), and \(A\) be the diagonal subgroup. Set
	\[
	h_u = \operatorname{diag}(u, u^{-1}), \quad z_0 \in F^\times.
	\]
	Let \(a_u = z_0 h_u\), so that \(z = (a_u)_{12} =\frac{z_0u}{z_0u^{-1}}= u^2\). For \(L = A, M\), the Fourier kernel at scalar elements is defined by
	\begin{equation*}
		\begin{aligned}
			\Phi^M_{A,L}(z_0 I_2, \tau) 
			= \lim_{u \to 1^+} \Big[ 
			\Phi^M_{A,L}(a_u, \tau)+\sqrt{2}\log(u - u^{-1}) \, \Phi^M_{M,L}(a_u, \tau) 
			\Big].
		\end{aligned}
	\end{equation*}
	
	Let \(\chi = \chi_1 \otimes \chi_2\), \(s = \lambda_{12} = \lambda_1 - \lambda_2\), and \(\psi = \Gamma'/\Gamma\). 
	
	For \(L=A\), we have \[\Phi^M_{M,A}(a_u, \chi) = \chi(a_u) + w\chi(a_u) = \chi(z_0 I_2)(u^s + u^{-s}),\]
	then the Foruier kernel is
	\begin{equation*}
		\begin{aligned}
			\Phi^M_{A,A}(z_0I_2,\chi)
			=\lim_{u\to1^+}\Big[&\Phi^M_{A,A}(a_u,\chi)+\sqrt2\,\chi(z_0 I_2)(u^s+u^{-s})\log(u-u^{-1})\Big].
		\end{aligned}
	\end{equation*}
	For \(F = \mathbb{R}\), write
	\[
	\chi_i(x) = \operatorname{sgn}(x)^{\varepsilon_i} |x|^{\lambda_i},
	\qquad \varepsilon_i \in \{0, 1\}, \quad \lambda_i \in i\mathbb{R},
	\]
	and set \(N \equiv \varepsilon_1 + \varepsilon_2 \pmod 2\), with \(N \in \{0, 1\}\). 
	
	We have
	\begin{equation*}
		\begin{aligned}
			\Phi^M_{A,A}(z_0I_2,\chi)
			=&\sqrt2\, \chi(z_0 I_2)\lim_{u\to1^+}\Bigg\{u^s\left[
			b(-s,u^{-2})
			+\frac12\psi\!\left(\frac{N-s}{2}\right)
			-\frac12\psi\!\left(\frac{N+1-s}{2}\right)\right]\\
			+&u^{-s}\left[
			b(s,u^{-2})
			+\frac12\psi\!\left(\frac{N+s}{2}\right)
			-\frac12\psi\!\left(\frac{N+1+s}{2}\right)\right]+(u^s+u^{-s})\log(u-u^{-1})\Bigg\}\\
			=&\sqrt{2}\, \chi(z_0 I_2) \Bigg[ 
			2\psi(1) - 2\log 2 -\psi\!\left(\frac{1 + N + s}{2}\right)-\psi\!\left(\frac{1 + N - s}{2}\right) 
			\Bigg].
		\end{aligned}
	\end{equation*}
	
	For \(F = \mathbb{C}\), using the standard complex modulus, write
	\[
	\chi_i(z) = |z|^{\lambda_i} (z/|z|)^{m_i},
	\qquad \lambda_i \in i\mathbb{R}, \quad m_i \in \mathbb{Z}.
	\]
	Set \(m = m_1 - m_2\), \(p = \frac{s + m}{2}\), and \(q = \frac{s - m}{2}\). Then
	\begin{equation*}
		\begin{aligned}
			\Phi^M_{A,A}(z_0I_2,\chi)
			=&\sqrt2\, \chi(z_0 I_2)\lim_{u\to1^+}\Bigg\{u^s\left[
			\frac{b(-p,u^{-2})+b(-q,u^{-2})}{2}
			+\frac1{|m|-s}\right]\\
			+&u^{-s}\left[
			\frac{b(p,u^{-2})+b(q,u^{-2})}{2}
			+\frac1{|m|+s}\right]+(u^s+u^{-s})\log(u-u^{-1})\Bigg\}\\
			=& \sqrt{2}\, \chi(z_0 I_2) \Bigg[ 
			2\psi(1) -\psi\!\left(\frac{|m| + s}{2}\right)-\psi\!\left(\frac{|m| - s}{2}\right) 
			\Bigg].
		\end{aligned}
	\end{equation*}
	For \(m = s = 0\), the limit is taken at the removable singularity, yielding \(0\).
	
	For \(L=M\), \(F=\mathbb{R}\), since \[\Phi^M_{M,M}(a_u,\sigma)=0,\] we have
	\begin{equation*}
		\Phi^M_{A,M}(z_0 I_2, \sigma) = -2 \, \omega_\sigma(z_0),
		\qquad \sigma \in \Pi_2(M),
	\end{equation*}
	where \(\sigma(z_0 I_2) = \omega_\sigma(z_0) \operatorname{Id}\). and \[
	\Phi^M_{A,M}\bigl(z_0I_2,(\xi\otimes\xi)^M\bigr)
	=\frac12\,\xi(z_0^2),
	\qquad \xi\in\widehat{F^\times}.
	\]
	
	For \(F=\mathbb{C}\), we have
	\begin{equation*}
		\Phi^M_{A,M}\bigl(z_0 I_2, (\xi \otimes \xi)^M\bigr) 
		= \frac{1}{2} \, \xi(z_0^2),
		\qquad \xi \in \widehat{F^\times}.
	\end{equation*}
	There is no discrete series term for the complex field.
	
	Therefore, we have 
	\begin{lemma}
		For \(F=\mathbb R\) or \(\mathbb C\) and \(z_0\in F^\times\), the Fourier kernel is
		\begin{align*}
			&I_{M_0,M_0}(z_0I_2,\chi)
			=\sqrt2\omega_\chi(z_0)^{-1}\times
			\begin{cases}
				\displaystyle
				2\psi(1)-2\log2
				-\psi\left(\frac{2-N+s}{2}\right)
				-\psi\left(\frac{2-N-s}{2}\right),\\
				&F=\mathbb R,\\
				\displaystyle
				2\psi(1)
				-\psi\left(\frac{|m|+s}{2}\right)
				-\psi\left(\frac{|m|-s}{2}\right),\\
				&F=\mathbb C,
			\end{cases}
		\end{align*}
		for \(\chi\in\Pi_2(M_0(F))\). When \(F=\mathbb R\),
		\[
		I_{M_0,G}(z_0I_2,\tau_\pi)
		=-2\omega_\pi(z_0)^{-1},
		\]
		for \(\pi\in\Pi_2(G(F))\), when \(F=\mathbb{R}\ \text{or}\ \mathbb{C}\),
		\[
		I_{M_0,G}(z_0I_2,\tau_\eta)
		=\frac{1}{2}\eta(z_0)^{-2},
		\]
		for \(\eta\in\Pi_2(\mathrm{GL}(1,F))\).
	\end{lemma}

	\section{Invariant distribution}
	For $P=G$, $\gamma_G\in G_{ell}$, we have \[I_G(\gamma_G,f)=|D_G(\gamma_G)|^{\frac{1}{2}}\int_{G_{\gamma_G}(F)\backslash G(F)}f(x^{-1}\gamma_G x)dx.\]
	\subsection{The maximal Levi}
    In this subsection, we discuss the case \(p\)-adic field \(F\).
    
	Arthur proved the descent formula:
	\begin{proposition}\label{descent}\cite{Art81}
		Suppose that \( \{ \phi_M^L \} \) satisfies (*) and (**). Then for every family \( \{ J^L \} \) of distributions satisfying (\ref{***}) there is a unique family \( \{ I^L \in U(L)' \} \) of invariant distributions such that for every \( f \in U(L) \),
		\[
		J^L(f) = \sum_{M \in \mathcal{L}^L(M_0)} c(M) c(L)^{-1} \hat{I}^M(\phi_M^L(f)).
		\]
	\end{proposition}
	For \(P=P_{21}\), \(\gamma_{21}\in M_{21,\rm{ell}}\) by Prop \ref{descent}, we have
	\[I_{M_{21}}(\gamma_{21},f)=J_{M_{21}}(\gamma_{21},f)-\hat{I}_{M_{21}}^{M_{21}}(\gamma_{21},\phi_{M_{21}}^G(f)).\]
	Since \(\phi^M_M\) is surjective, choose
\[
h_{21}\in\mathcal {H}_{\rm{ac}}(M_{21}(F))
\]
such that
\[
\phi_{M_{21}}^{M_{21}}(h_{21})
=
\phi_{M_{21}}^G(f).
\]
In particular,
\[
\operatorname{tr}\tau(h_{21})
=
\phi_{M_{21}}^{M_{21}}(h_{21})(\tau)
=
J_{M_{21}}(\tau,f),
\qquad
\tau\in\Pi_{\mathrm{temp}}(M_{21}(F)).
\]
Consequently,
\[
\begin{aligned}
\widehat I_{M_{21}}^{M_{21}}
\bigl(\gamma_{21},\phi_{M_{21}}^G(f)\bigr)
&=
\widehat I_{M_{21}}^{M_{21}}
\bigl(\gamma_{21},
      \phi_{M_{21}}^{M_{21}}(h_{21})\bigr)\\
&=
I_{M_{21}}^{M_{21}}(\gamma_{21},h_{21}).
\end{aligned}
\]

Write
\[
\gamma_{21}=\gamma^1\times\gamma^2,
\qquad
\gamma^1\in\mathrm{GL}(2,F)_{\mathrm{ell,reg}},
\quad
\gamma^2\in\mathrm{GL}(1,F).
\]
Then
\[
\begin{aligned}
I_{M_{21}}^{M_{21}}(\gamma_{21},h_{21})
=&|D_{M_{21}}(\gamma_{21})|^{1/2}
\int_{\mathrm{GL}(2,F)_{\gamma^1}\backslash\mathrm{GL}(2,F)}
h_{21}(x^{-1}\gamma^1x,\gamma^2)\,dx,
\end{aligned}
\]
where
\[
D_{M_{21}}(\gamma_{21})
=
D_{\mathrm{GL}(2)}(\gamma^1).
\]

For
\[
\chi=\chi_1\otimes\chi_2\otimes\chi_3
\in\Pi_{\mathrm{unit}}(M_0(F)),
\]
write
\[
\chi^{M_{21}}
=
\mathscr I_{P_0\cap M_{21}}^{M_{21}}(\chi).
\] 
Fourier inversion on \(M_{21}\)
gives
\[
\begin{aligned}
I_{M_{21}}^{M_{21}}(\gamma_{21},h_{21})
=&
|D_{M_{21}}(\gamma_{21})|^{\frac{1}{2}}
\int_{\Pi_2(M_{21}(F))}
\Theta_\tau(\gamma_{21}^{-1})
\operatorname{tr}\tau(h_{21})\,d\tau\\
&+
\frac12
\int_{\Pi_{\mathrm{unit}}(M_0(F))}
I_{M_{21},M_0}^{M_{21}}(\gamma_{21},\chi)
\operatorname{tr}\chi^{M_{21}}(h_{21})\,d\chi.
\end{aligned}
\]
The kernels in this expansion are compatible with the product
decomposition. More precisely, if \(M_0'\) denotes the diagonal
Levi subgroup of \(\mathrm{GL}(2)\), then
\[
I_{M_{21},M_0}^{M_{21}}(\gamma_{21},\chi)
=
I_{\mathrm{GL}_2,M_0'}^{\mathrm{GL}_2}
(\gamma^1,\chi_1\otimes\chi_2)\,
\chi_3(\gamma^2)^{-1}.
\]
Likewise, for
\[
\tau=\tau_1\otimes\eta,
\qquad
\tau_1\in\Pi_2(\mathrm{GL}(2,F)),
\quad
\eta\in\Pi_{\mathrm{unit}}(\mathrm{GL}(1,F)),
\]
we have
\[
\Theta_\tau(\gamma_{21}^{-1})
=
\Theta_{\tau_1}\bigl((\gamma^1)^{-1}\bigr)
\eta(\gamma^2)^{-1}.
\]
Thus the preceding expansion uses the previously computed
\(\mathrm{GL}(2)\) Fourier kernels together with Fourier inversion
on \(\mathrm{GL}(1,F)\).

Substituting
\[
\operatorname{tr}\tau(h_{21})=J_{M_{21}}(\tau,f),
\qquad
\operatorname{tr}\chi^{M_{21}}(h_{21})
=
J_{M_{21}}(\chi^{M_{21}},f),
\]
we obtain
\begin{equation}\label{eq:maximal-levi-invariant}
\begin{aligned}
I_{M_{21}}(\gamma_{21},f)
=&
J_{M_{21}}(\gamma_{21},f)\\
-&
|D^{M_{21}}(\gamma_{21})|^{1/2}
\int_{\Pi_2(M_{21}(F))}
\Theta_\tau(\gamma_{21}^{-1})
J_{M_{21}}(\tau,f)\,d\tau\\
-&
\frac12
\int_{\Pi_{\mathrm{unit}}(M_0(F))}
I_{M_{21},M_0}^{M_{21}}(\gamma_{21},\chi)
J_{M_{21}}(\chi^{M_{21}},f)\,d\chi.
\end{aligned}
\end{equation}
The expression is independent of the choice of \(h_{21}\),
since \(I_{M_{21}}^{M_{21}}(\gamma_{21},\cdot)\) factors through
the invariant Fourier transform.

\subsection{The minimal Levi}\label{subsec:minimal-levi}
In this subsection, we discuss the case \(p\)-adic field \(F\).

Fix
\[
\gamma\in M_0(F)\cap G_{\mathrm{reg}}(F),
\qquad
f\in\mathcal H(G(F)).
\]
Let \(\mathscr W\) be a set of representatives for
\(W_0^G/W_0^{M_{21}}\), and put
\[
\gamma_w=w^{-1}\gamma w,
\qquad w\in\mathscr W,
\]
where Weyl group elements are represented by permutation matrices.
The three intermediate Levi subgroups containing \(M_0\) are
\(wM_{21}w^{-1}\), with \(w\in\mathscr W\).
Then
\begin{equation}\label{eq:minimal-levi-recursion}
\begin{aligned}
I_{M_0}(\gamma,f)
=&
J_{M_0}(\gamma,f)
-
\widehat I_{M_0}^{M_0}
\bigl(\gamma,\phi_{M_0}^G(f)\bigr)-
\sum_{w\in\mathscr W}
\widehat I_{M_0}^{M_{21}}
\bigl(\gamma_w,\phi_{M_{21}}^G(f)\bigr).
\end{aligned}
\end{equation}
Choose
\[
h_0\in\mathcal{H}_{\rm{ac}}(M_0(F)),
\qquad
h\in\mathcal{H}_{\rm{ac}}(M_{21}(F)),
\]
such that
\[
\phi_{M_0}^{M_0}(h_0)=\phi_{M_0}^G(f),
\qquad
\phi_{M_{21}}^{M_{21}}(h)=\phi_{M_{21}}^G(f).
\]
Since \(M_0(F)\) is abelian, Fourier inversion gives
\begin{equation}\label{eq:minimal-levi-torus-term}
\begin{aligned}
\widehat I_{M_0}^{M_0}
\bigl(\gamma,\phi_{M_0}^G(f)\bigr)
&=
I_{M_0}^{M_0}(\gamma,h_0)
=
h_0(\gamma)\\
&=
\int_{\Pi_{\mathrm{unit}}(M_0(F))}
\tau(\gamma^{-1})\operatorname{tr}\tau(h_0)\,d\tau\\
&=
\int_{\Pi_{\mathrm{unit}}(M_0(F))}
\tau(\gamma^{-1})J_{M_0}(\tau,f)\,d\tau.
\end{aligned}
\end{equation}
We next compute the contribution associated with \(M_{21}\).
By the choice of \(h\),
\[
\widehat I_{M_0}^{M_{21}}
\bigl(\gamma,\phi_{M_{21}}^G(f)\bigr)
=
I_{M_0}^{M_{21}}(\gamma,h).
\]
Write
\[
M_0=M_0'\times\mathrm{GL}_1,
\qquad
\gamma=\gamma^1\times\gamma^2,
\]
where
\[
\gamma^1\in M_0'(F)\cap\mathrm{GL}(2,F)_{\mathrm{reg}},
\qquad
\gamma^2\in\mathrm{GL}(1,F).
\]
Therefore,
\[
J_{M_0}^{M_{21}}(\gamma,h^1\otimes h^2)
=
J_{M_0'}^{\mathrm{GL}_2}(\gamma^1,h^1)\,h^2(\gamma^2)
\]
for
\[
h^1\in\mathcal {H}(\mathrm{GL}(2,F)),
\qquad
h^2\in\mathcal {H}(\mathrm{GL}(1,F)).
\]
Then
\[
J_{M_0}^{M_{21}}
(\tau'\otimes\eta,h^1\otimes h^2)
=
J_{M_0'}^{\mathrm{GL}_2}(\tau',h^1)\,
\operatorname{tr}\eta(h^2),
\]
and
\[
I_{M_0}^{M_{21}}(\gamma,h^1\otimes h^2)
=
I_{M_0'}^{\mathrm{GL}_2}(\gamma^1,h^1)\,h^2(\gamma^2).
\]
Consequently, the Fourier kernels satisfy
\begin{equation}\label{eq:minimal-levi-product-kernels}
\begin{aligned}
I_{M_0,M_0}^{M_{21}}
(\gamma,\tau'\otimes\eta)
&=
I_{M_0',M_0'}^{\mathrm{GL}_2}
(\gamma^1,\tau')\,\eta(\gamma^2)^{-1},\\
I_{M_0,M_{21}}^{M_{21}}
(\gamma,\sigma\otimes\eta)
&=
I_{M_0',\mathrm{GL}_2}^{\mathrm{GL}_2}
(\gamma^1,\sigma)\,\eta(\gamma^2)^{-1},
\end{aligned}
\end{equation}
where
\[
\tau'\in T_{\mathrm{disc}}(M_0'),
\qquad
\sigma\in T_{\mathrm{disc}}(\mathrm{GL}(2)),
\qquad
\eta\in\Pi_{\mathrm{unit}}(\mathrm{GL}(1,F)).
\]
The invariant Fourier inversion formula on \(M_{21}\) now gives
\begin{equation}\label{eq:minimal-levi-fourier-inversion}
\begin{aligned}
I_{M_0}^{M_{21}}(\gamma,h)
=&
\frac12
\int_{T_{\mathrm{disc}}(M_0)}
I_{M_0,M_0}^{M_{21}}(\gamma,\tau)\,
h_{M_0}(\tau)\,d\tau\\
&+
\int_{T_{\mathrm{disc}}(M_{21})}
I_{M_0,M_{21}}^{M_{21}}(\gamma,\tau)\,
h_{M_{21}}(\tau)\,d\tau.
\end{aligned}
\end{equation}
By the choice of \(h\),
\[
\begin{aligned}
h_{M_0}(\tau)
&=
\operatorname{tr}\tau^{M_{21}}(h)
=
J_{M_{21}}(\tau^{M_{21}},f),\\
h_{M_{21}}(\tau)
&=
J_{M_{21}}(\tau,f).
\end{aligned}
\]
Substituting into
\eqref{eq:minimal-levi-fourier-inversion}, we obtain
\begin{equation}\label{eq:minimal-levi-maximal-correction}
\begin{aligned}
\widehat I_{M_0}^{M_{21}}
\bigl(\gamma,\phi_{M_{21}}^G(f)\bigr)
=&
\frac12
\int_{T_{\mathrm{disc}}(M_0)}
I_{M_0,M_0}^{M_{21}}(\gamma,\tau)\,
J_{M_{21}}(\tau^{M_{21}},f)\,d\tau\\
+&
\int_{T_{\mathrm{disc}}(M_{21})}
I_{M_0,M_{21}}^{M_{21}}(\gamma,\tau)\,
J_{M_{21}}(\tau,f)\,d\tau.
\end{aligned}
\end{equation}
Since
\[
T_{\mathrm{disc}}(M_0)
=
\Pi_{\mathrm{unit}}(M_0(F)),
\]
combining
\eqref{eq:minimal-levi-recursion},
\eqref{eq:minimal-levi-torus-term}, and
\eqref{eq:minimal-levi-maximal-correction}
gives
\begin{equation}\label{eq:minimal-levi-invariant}
\begin{aligned}
I_{M_0}(\gamma,f)
=&
J_{M_0}(\gamma,f)-
\int_{\Pi_{\mathrm{unit}}(M_0(F))}
\tau(\gamma^{-1})J_{M_0}(\tau,f)\,d\tau\\
-&
\frac12\sum_{w\in\mathscr W}
\int_{\Pi_{\mathrm{unit}}(M_0(F))}
I_{M_0,M_0}^{M_{21}}(\gamma_w,\tau)\,
J_{M_{21}}(\tau^{M_{21}},f)\,d\tau\\
-&
\sum_{w\in\mathscr W}
\int_{T_{\mathrm{disc}}(M_{21})}
I_{M_0,M_{21}}^{M_{21}}(\gamma_w,\tau)\,
J_{M_{21}}(\tau,f)\,d\tau.
\end{aligned}
\end{equation}
Thus,
	\begin{theorem}\label{main1}
		Suppose \(f\in\mathcal{H}(G(F_v))\). For \(\gamma_G\in G_{\rm{ell}}\), \[I_G(\gamma_G,f)=|D_G(\gamma_G)|^{\frac{1}{2}}\int_{G_{\gamma_G}(F)\backslash G(F)}f(x^{-1}\gamma_G x)dx.\] 
		For \(\gamma_{21}\in M_{21}(F)_{\rm{ell}}\), the invariant distribution \(I_{M_{21}}(\gamma_{21},f)\) equals 
        \begin{align*}
            &J_{M_{21}}(\gamma_{21},f)\\
-&
|D^{M_{21}}(\gamma_{21})|^{1/2}
\int_{\Pi_2(M_{21}(F))}
\Theta_\tau(\gamma_{21}^{-1})
J_{M_{21}}(\tau,f)\,d\tau\\
-&
\frac12
\int_{\Pi_{\mathrm{unit}}(M_0(F))}
I_{M_{21},M_0}^{M_{21}}(\gamma_{21},\chi)
J_{M_{21}}(\chi^{M_{21}},f)\,d\chi.
        \end{align*}
		For \(\gamma\in M_0(F)\),
		\(I_{M_0}(\gamma,f)\) equals 
		\begin{align*}
			&J_{M_0}(\gamma,f)-
\int_{\Pi_{\mathrm{unit}}(M_0(F))}
\tau(\gamma^{-1})J_{M_0}(\tau,f)\,d\tau\\
-&
\frac12\sum_{w\in\mathscr W}
\int_{\Pi_{\mathrm{unit}}(M_0(F))}
I_{M_0,M_0}^{M_{21}}(\gamma_w,\tau)\,
J_{M_{21}}(\tau^{M_{21}},f)\,d\tau\\
-&
\sum_{w\in\mathscr W}
\int_{T_{\mathrm{disc}}(M_{21})}
I_{M_0,M_{21}}^{M_{21}}(\gamma_w,\tau)\,
J_{M_{21}}(\tau,f)\,d\tau.
		\end{align*}
	\end{theorem}
    \subsection{The global case}
	For any finite set \(S\), we can take any subset \(S_1\), \(S_2\subset S\) such that \(S=S_1\sqcup S_2\), we have the splitting formula in \cite{Art88(1)}:
	\begin{proposition}
		Suppose that
		\[
		\gamma = \gamma_1 \gamma_2, \qquad \gamma_i \in M(F_{S_i}),
		\]is any element in \(M(F_S)\). Then for any function \(f \in \mathscr{H}_{\mathrm{ac}}(G(F_S))\) of the form
		\[
		f = f_1 f_2, \qquad f_i \in \mathscr{H}_{\mathrm{ac}}(G(F_{S_i})),
		\]
		we have
		\[
		I_M(\gamma, f) = \sum_{L_1, L_2 \in \mathcal{L}(M)} d_M^G(L_1, L_2) \widehat{I}_M^{L_1}(\gamma_1, f_{1, L_1}) \widehat{I}_M^{L_2}(\gamma_2, f_{2, L_2}).
		\]
	\end{proposition}
	We can now write the invariant distribution in the case of \(\GGL(3,F_S)\). Since Hoffmann and Wakatsuki has calculated the coefficient \(a^M(\gamma,S)\) in \cite{HW}, we can get the geometric side of invariant trace formula of \(\GGL(3)\):
	\begin{theorem}\label{main2}
		For \(f\in\mathcal{H}(G(\mathbb{A})\), the geometric side of invariant trace formula of \(\GGL(3)\) is 
		\[
		I_{\rm{geom}}(f) = \lim_{S} \sum_{M \in \mathcal{L}} |W_0^M| |W_0^G|^{-1} \sum_{\gamma \in \Gamma(M)_S} a^M(\gamma,S) I_M(\gamma, f)
		\]
	\end{theorem}

	\section{The spectral side}
	Recall that Arthur (in \cite{Art88(1)}) writes the invariant trace formula as
	\begin{equation*}
		I(f) = \sum_{t \ge 0} I_t(f),
	\end{equation*}
	where
	\begin{equation*}
		I_t(f) = I_t^G(f) = \sum_{M \in \mathcal{L}} |W_0^M| |W_0^G|^{-1} \int_{\Pi(M,t)} a^M(\pi) I_M(\pi, f) \, \mathrm{d}\pi.
	\end{equation*}\(\Pi(M, t)\) can be written according to the discrete data in \(M_1 \subseteq M\) as
	\begin{equation*}
		\Pi(M, t) = \bigsqcup_{M_1 \in \mathcal{L}^M} \{ \pi_{1, \lambda}^M : \pi_1 \in \Pi_{\mathrm{disc}}(M_1, t), \quad \lambda \in i\mathfrak{a}_{M_1}^* / i\mathfrak{a}_M^* \},
	\end{equation*}
	where \(\pi_{1, \lambda}^M\) is the normalized unitary induction from \(M_1\) to \(M\). The corresponding measure is
	\begin{equation*}
		\begin{split}
			\int_{\Pi(M,t)} \Phi(\pi) \, d\pi = & \sum_{M_1 \in \mathcal{L}^M} \frac{|W_0^{M_1}|}{|W_0^M|} \sum_{\pi_1 \in \Pi_{\mathrm{disc}}(M_1, t)}\int_{i\mathfrak{a}_{M_1}^* / i\mathfrak{a}_M^*} \Phi(\pi_{1,\lambda}^M) \, d\lambda.
		\end{split}
	\end{equation*}
	The spectral coefficients decompose as
	\begin{equation*}
		a^M(\pi_{1,\lambda}^M) = a_{\mathrm{disc}}^{M_1}(\pi_1)r_{M_1}^M(\pi_1, \lambda).
	\end{equation*}
	
	Then we have 
	\begin{align}\label{dist}
		I_t(f) = & \sum_{M \in \mathcal{L}(T)} \sum_{M_1 \subset M}\frac{|W_0^{M_1}|}{|W_0^G|} \sum_{\sigma \in \Pi_{\mathrm{disc}}(M_1, t)} a_{\mathrm{disc}}^{M_1}(\sigma)\int_{i\mathfrak{a}_{M_1}^* / i\mathfrak{a}_M^*} r_{M_1}^M(\sigma, \lambda) I_M \left( \mathrm{Ind}_{P_1\cap M}^M(\sigma_\lambda), f \right) d\lambda.
	\end{align}
	Write \(P_0'\) as the borel subgroup of \(\GGL(2)\).
	
	\subsection{The case \(M=G\)}
	In this case, let \(F=F_v=F_S\) be local field. We choose \(f\in\mathcal{H}(G(F))\).
    
    For any Levi subgroup \(M_1\), \(\lambda\in i\mathfrak{a}_{M_1}^* / i\mathfrak{a}_M^*\) and \(\sigma\in \Pi_{\rm{disc}}(M_1,t)\), we have \[I_G(\mathscr{I}_{M_1}(\sigma_\lambda),f)=\rm{tr}(\mathscr{I}_{M_1}(\sigma_\lambda)(f)).\]
	Let
	\[
	m_{M_{21}}(z) = r_{\bar{P}_{21}|P_{21}}(\sigma_{z})=r_{\bar{P}_{21}|P_{21}}\left(\sigma\otimes \left( | \det |^{z/3} \otimes | \cdot |^{-2z/3} \right)\right) .
	\]
	Then
	\begin{align*}
		r_{M_{21}}^G(\sigma, z) &= \sqrt{\frac{3}{2}}\lim_{u \to 0} \left[ \frac{1}{u} - \frac{m_{M_{21}}(z+u)}{u m_{M_{21}}(z)} \right] \\
		&=  -\sqrt{\frac{3}{2}} \frac{m_{M_{21}}'(z)}{m_{M_{21}}(z)}  .
	\end{align*}
	For
	\[
	\chi = \chi_1 \otimes \chi_2 \otimes \chi_3, \qquad z_{ij} = \lambda_i - \lambda_j,
	\]
	let \(m_{ij}(z)\) be the global rank-one normalization factor for the root \(e_i - e_j\), and denote
	\[
	\ell_{ij} = \frac{m'_{ij}(z_{ij})}{m_{ij}(z_{ij})}.
	\]
	Let
	\[
	u = \Lambda_1 - \Lambda_2, \qquad v = \Lambda_2 - \Lambda_3,
	\]
	and
	\[
	a = \frac{m_{12}(z_{12} + u)}{m_{12}(z_{12})}, \quad b = \frac{m_{23}(z_{23} + v)}{m_{23}(z_{23})}, \quad c = \frac{m_{13}(z_{13} + u + v)}{m_{13}(z_{13})}.
	\]
	Therefore
	\[
	r_{M_0}^G(\chi, \lambda) = \sqrt{3} \lim_{u, v \to 0} \left[ \frac{1 + abc}{uv} - \frac{b + ac}{v(u+v)} - \frac{a + bc}{u(u+v)} \right].
	\]
	By performing a Taylor expansion with respect to 
	\(u\) and \(v\), we obtain
	\[
	a = 1 + \ell_{12}u + O(u^2), \qquad b = 1 + \ell_{23}v + O(v^2),
	\]
	\[
	c = 1 + \ell_{13}(u + v) + O((u + v)^2),
	\]
	then
	\[
	r_{M_0}^G(\chi, \lambda) =\sqrt{3} (\ell_{12}\ell_{23} + \ell_{12}\ell_{13} + \ell_{23}\ell_{13}).
	\]

	\subsection{The case  \(M=M_{21}\)}
	By the classification of the representation of \(\GGL(2)\), we have 
	\begin{proposition}\label{gl2}
		Let \(F\) be a local field of characteristic zero and set
		\[
		G = \operatorname{GL}(2,F), \qquad \nu = |\cdot|_F.
		\]
		Then every irreducible unitary representation of \(G\) belongs to exactly one of the following classes:
		\[
		\left\{
		\begin{aligned}
			& \pi \in \Pi_{\mathrm{temp}}(G(F)), \\
			& \pi \simeq \mathscr{I}_{P'_0}^G (\xi \nu^s \otimes \xi \nu^{-s}), \quad \xi \in \Pi_{\rm{unit}}(\GGL(1,F)), \quad 0 < s < \frac{1}{2}, \\
			& \pi \simeq \xi \circ \det, \quad \xi \in \Pi_{\rm{unit}}(\GGL(1,F)).
		\end{aligned}
		\right.
		\]
	\end{proposition}
	This classification is the result in (\cite{JL70}).
	\subsubsection{Tempered representation}
	If \(\pi\) is a tempered representation of \(M\), then
	\begin{align}\label{tem}
			I_M(\pi, f) = 
		\begin{cases}
			\operatorname{tr} \pi(f), & M = G, \\
			0, & M \neq G,
		\end{cases}
	\end{align}
by \cite{Art81}.
	\subsubsection{Complementary Series}
	By the definition of complementary series, set \[
	\rho_s = \mathscr{I}_{P'_0}^{GL_2} \left( \xi v^s \otimes \xi v^{-s} \right), \quad 0 < s < \frac{1}{2},
	\]and define 
	\[\sigma_s = \mathscr{I}_{P'_0}^{\operatorname{GL}(2)} (\xi v^s \otimes \xi v^{-s}) \otimes \eta,\qquad \eta\in\Pi_{\rm{unit}}(\GGL(1,F)),\]then \(\sigma_s\) is the unitary but not tempered representation of \(M_{21}(F)\).

	Since \(\rm{dim}(\mathfrak{a}_{M_{21}}^G)=1\), for \(z\in i\mathfrak{a}_{M_{21}}^G\), write \[\lambda_z=
\left(\frac z3,\frac z3,-\frac{2z}3\right),\qquad
	\sigma_{s,z} = \sigma_s \otimes \left( | \det |^{z/3} \otimes | \cdot |^{-2z/3} \right).
	\]Thus we can write the induced representation as 
	\[
    \mathscr{I}_{P_0}^{\operatorname{GL}(3)} \left( \xi v^{s+z/3} \otimes \xi v^{-s+z/3} \otimes \eta v^{-2z/3} \right).
	\]
    Write \(\mathfrak a_{M,F}\) as the image of \(M(F)\) under \(H_M(x)\).
    \begin{lemma}\label{lem:standard-vanishing-M21}
For \(0\le s\le\frac12\), we have
\[
I_{M_{21}}^G(\sigma_s,X,f)=0,
\qquad
X\in\mathfrak a_{M_{21},F},
\quad f\in\mathcal H(G(F)).
\]
\end{lemma}
\begin{proof}
Fix \(X\) and \(f\), and set
\[
\Omega=\{s\in\mathbb C:|\operatorname{Re}s|<3/4\}.
\]
The intertwining operator \(\mathscr R_{\overline P_{21}\mid P_{21}}(\sigma_{s,z})\)
factors into rank-one operators associated with
\(\alpha_{13}\) and \(\alpha_{23}\), with respective character ratios
\[
\xi\eta^{-1}\nu^{z+s},
\qquad
\xi\eta^{-1}\nu^{z-s}.
\]
For unitary character, the corresponding normalized
rank-one operators and their inverses are holomorphic in
\(|\operatorname{Re}u|<1\) (see
\cite[\S I.10]{MW89}), applied in both
directions, together with the functional equation.

For \(s\in\Omega\) and \(z\in i\mathbb R\), both \(z+s\)
and \(z-s\) lie in this strip. Hence the block-intertwining
operators and their logarithmic derivatives are holomorphic
in \(s\) along the unitary integration contour. The internal
root \(e_1-e_2\), with parameter \(2s\), does not occur in
these operators.

It follows that
\[
J_M^G(\sigma_s,X,f)
=
\int_{i\mathfrak a_{M,F}^*}
J_M^G(\sigma_{s,z},f)
e^{-\langle z,X\rangle}\,dz
\]
is holomorphic on \(\Omega\). For non-Archimedean \(F\), the integration domain is the compact torus
\[
i\mathfrak a_{M,F}^*
=
i\mathfrak a_M^*/\Lambda_{M,F},
\]
where
\[
\Lambda_{M,F}
=
\left\{
\lambda\in i\mathfrak a_M^*:
\langle\lambda,H\rangle\in2\pi i\mathbb Z
\ \text{for all }H\in\mathfrak a_{M,F}
\right\},
\qquad
\mathfrak a_{M,F}=H_M(M(F)).
\] 
In the Archimedean case,
the estimates for induced representations and weighted
characters justify integration locally uniformly in \(s\) (see \cite[\S6 and \S12, (12.7)]{Art89}).

The correction term
\[
\widehat I_M^M
\bigl(\sigma_s,X,\phi_M^G(f)\bigr)
\]
is entire in \(s\). Indeed, multiply \(\phi_M^G(f)\) by
a compactly supported function of its central variable
that equals \(1\) near \(X\). The resulting invariant
transform has a compactly supported Hecke lift on \(M\),
and its evaluation at \(X\) is unchanged.
The assertion follows from the Paley-Wiener property
of ordinary induced characters (see
\cite[\S1]{Art88}).

Consequently,
\[
I_{M_{21}}^G(\sigma_s,X,f)
=
J_M^G(\sigma_s,X,f)-\widehat I_M^M
\bigl(\sigma_s,X,\phi_M^G(f)\bigr)
\]
is holomorphic on \(\Omega\). For the real parameters
in the statement, this is \(I_M^G(\sigma_s,X,f)\), with
the standard-representation interpretation at \(s=\frac{1}{2}\)(see \cite[\S\S5--7]{Art89}).

For \(s\in i\mathbb R\), the representation \(\sigma_s\)
is tempered. Since \(M\ne G\), Arthur's tempered vanishing
lemma gives
\[
I_{M_{21}}^G(\sigma_s,X,f)=0,\qquad s\in i\mathbb R;
\]
see \cite[Lemma~3.1]{Art88}.
The identity theorem therefore implies
\(I_{M_{21}}^G(\sigma_s,X,f)=0\) on \(\Omega\).
Since \([0,\frac12]\subset\Omega\), the result follows.
\end{proof}

	\subsubsection{One-dimensional Representation}\label{quotent}
	Set \[q_\xi = \mathscr{I}_{P'_0}^{\GGL(2)} \left( \xi \nu^{1/2} \otimes \xi \nu^{-1/2} \right).\]
	It is easy to check that
	\[q_\xi \simeq \left[ \mathscr{I}_{P_0'}^{\GGL(2)} \left( \nu^{1/2} \otimes \nu^{-1/2} \right) \right] \otimes (\xi \circ \det).\]
	Define \[q_1=\mathscr{I}_{P_0'}^{\GGL(2)} \left( \nu^{1/2} \otimes \nu^{-1/2} \right) .\]
	We have a short exact sequence
	\[
	0 \longrightarrow \kappa_\xi \longrightarrow q_\xi \longrightarrow a_\xi \longrightarrow 0,
	\]where \(a_\xi(g)=\xi(\rm{det}\,g)\), and \(
	\kappa_\xi = \ker(q_\xi \longrightarrow a_\xi)
	\).

	For p-adic case, we have the short exact sequence
	\[0 \longrightarrow \operatorname{St}_{\GGL(2)} \longrightarrow \mathscr{I}_{P_0'}^{\GGL(2)} \left( \nu^{1/2} \otimes \nu^{-1/2} \right) \longrightarrow \mathbf{1} \longrightarrow 0,\]where \( \operatorname{St}_{\GGL(2)}\) is the Steinberg representation, and 
	\[0 \longrightarrow \operatorname{St}_{\GGL(2)} \otimes (\xi \circ \det) \longrightarrow q_\xi \longrightarrow \xi \circ \det \longrightarrow 0.\]According to Theorem 3.3 in \cite{JL70}, since \[\frac{v^{\frac{1}{2}}}{v^{-\frac{1}{2}}}=v,\]\(q_\xi\) is reducible.
    Also, the Steinberg representation is tempered, so is \(\operatorname{St}_{\GGL(2)} \otimes (\xi \circ \det) \).

	For \(F=\mathbb{R}\), by Theorem 5.11, we have  \[\kappa_\xi=\sigma(\xi v^{\frac{1}{2}},\xi v^{-\frac{1}{2}}),\]over the subspace \(B_s(\xi v^{\frac{1}{2}},\xi v^{-\frac{1}{2}})\) 
	and \(\kappa_\xi\) is the unique irreducible subrepresentation of \(q_\xi\).
	
	Then, by Corollary 5.14 in \cite{JL70}, \[\kappa_\xi\simeq \pi(w),\]which is discrete series, thus \(\kappa_\xi\) is tempered.
	
	For \(F=\mathbb{C}\), by Theorem 6.2 in \cite{JL70}, we can similarly get that \(\kappa_\xi\) is tempered.

	Define
	\[\sigma_a = a_\xi \otimes \eta, \qquad \sigma_\kappa = \kappa_\xi \otimes \eta, \qquad \sigma_q = q_\xi \otimes \eta,\]for \(\eta\in\Pi_{\rm{unit}}({\GGL(1,F)})\), we then have the short exact sequence
	\[0 \longrightarrow \sigma_\kappa \longrightarrow \sigma_q \longrightarrow \sigma_a \longrightarrow 0.\]
	Define 
	\[\sigma_a=(a_\xi |\rm{det}|^{\frac{z}{3}})\otimes (\eta |\cdot|^{-\frac{2z}{3}}),\quad \sigma_{\kappa,z}= (\kappa_\xi  |\rm{det}|^{\frac{z}{3}}) \otimes (\eta|\cdot|^{-\frac{2z}{3}}), \quad \sigma_{q,z} =(q_\xi  |\rm{det}|^{\frac{z}{3}})  \otimes (\eta|\cdot|^{-\frac{2z}{3}}),\]for \(z\in i\mathfrak{a}_{M_{21}}^G\).

	Since \(\kappa_\xi\) is tempered, the representation \(\sigma_\kappa\) is tempered. we have \[I_{M_{21}}^G(\sigma_\kappa, f) = 0,\qquad \text{since}\ M_{21}\neq G.\]
    Write three normalization factors
    \[
    \begin{aligned}
    	r_q(z) = r_{\overline{P}_{21}|P_{21}}(\sigma_{q,z}), \quad
    	r_\kappa(z) = r_{\overline{P}_{21}|P_{21}}(\sigma_{\kappa,z}), \quad
    	r_a(z) = r_{\overline{P}_{21}|P_{21}}(\sigma_{a,z}),
    \end{aligned}
    \]and 
    \begin{align*}
    	R_q(z) = \mathscr{R}_{\overline{P}_{21}|P_{21}}(\sigma_{q,z}), \quad R_\kappa(z) = \mathscr{R}_{\overline{P}_{21}|P_{21}}(\sigma_{\kappa,z}), \quad R_a(z) = \mathscr{R}_{\overline{P}_{21}|P_{21}}(\sigma_{a,z}).
    \end{align*}
    Define
    \[
    b_{\xi,\eta}(z) = \frac{r_a(z)}{r_\kappa(z)},\qquad h_{\xi,\eta}(z) = \frac{d}{dz} \log b_{\xi,\eta}(z) = \frac{b'_{\xi,\eta}(z)}{b_{\xi,\eta}(z)}.
    \]We have the short exact sequence for the induced representation
    \[0 \longrightarrow \mathscr{I}_{P_{21}}^G(\sigma_{\kappa,z}) \longrightarrow \mathscr{I}_{P_{21}}^G(\sigma_{q,z}) \longrightarrow \mathscr{I}_{P_{21}}^G(\sigma_{a,z}) \longrightarrow 0.\]
   For the restriction representation and quotient representation, we have 
   \[\mathscr{J}_{\overline{P}_{21}|P_{21}}(\sigma_{q,z})|_{\sigma_{\kappa,z}}=\mathscr{J}_{\overline{P}_{21}|P_{21}}(\sigma_{\kappa,z}),\qquad \overline{\mathscr{J}_{\overline{P}_{21}|P_{21}}(\sigma_{q,z})}=\mathscr{J}_{\overline{P}_{21}|P_{21}}(\sigma_{a,z}).\]
	According to the definition of normalized factor, we have 
	\[r_q(z)=r_a(z),\qquad r_a(z)=b_{\xi,\eta}(z)r_{\kappa}(z),\]then 
	\begin{align*}
		R_q(z)|_{\pi_\kappa} &= \frac{r_\kappa(z)}{r_q(z)} R_\kappa(z) \\
		&= \frac{r_\kappa(z)}{r_a(z)} R_\kappa(z),
	\end{align*}and 
	\begin{align*}
		\overline{R_q(z)} &= r_q(z)^{-1} \mathscr{J}_{\overline{P}_{21}|P_{21}}(\sigma_{a,z}) \\
		&= \frac{r_a(z)}{r_q(z)} R_a(z).
	\end{align*}
	By calculation, 
	\begin{align*}
		\left(R_q(z)^{-1} R'_q(z)\right) \big|_{\sigma_\kappa}=& R_\kappa(z)^{-1} b_{\xi,\eta}(z) \left( -b_{\xi,\eta}(z)^{-2} b'_{\xi,\eta}(z) R_\kappa(z) + b_{\xi,\eta}(z)^{-1} R'_\kappa(z) \right) \\
		=& -\frac{b'_{\xi,\eta}(z)}{b_{\xi,\eta}(z)} \mathrm{Id} + R_\kappa(z)^{-1} R'_\kappa(z)\\
		=&R_\kappa(z)^{-1} R'_\kappa(z)-h_{\xi,\eta}(z)\mathrm{Id}.
	\end{align*}
	On the other hand,
	\begin{align*}
		\overline{R_q(z)^{-1} R'_q(z)}=R_a(z)^{-1} R'_a(z).
	\end{align*}Therefore, 
	\begin{align*}
		\mathscr{R}_{M_{21}}(\sigma_{q,z})|_{\sigma_\kappa}= \mathscr{R}_{M_{21}}(\sigma_{\kappa,z}) + \sqrt{\frac{3}{2}}  h_{\xi,\eta}(z) \mathrm{Id}.
	\end{align*}and \[\overline{\mathscr{R}_{M_{21}}(\sigma_{q,z})} = \mathscr{R}_{M_{21}}(\sigma_{a,z},).\]
	Then 
	\begin{align*}
		J_{M_{21}}(\sigma_{q,z}, f) &= \mathrm{tr}\left( \mathscr{R}_{M_{21}}(\sigma_{q,z}) \mathscr{I}_{P_{21}}(\sigma_{q,z}, f) \right) \\
		&= \mathrm{tr}\left( \mathscr{R}_{M_{21}}(\sigma_{q,z}) |_{\sigma_\kappa} \mathscr{I}_{P_{21}}(\sigma_{\kappa,z}, f) \right) + \mathrm{tr}\left( \overline{\mathscr{R}_{M_{21}}(\sigma_{q,z}) } \mathscr{I}_{P_{21}}(\sigma_{a,z}, f) \right) \\
		&= \mathrm{tr}\left( \left[ \mathscr{R}_{M_{21}}(\sigma_{\kappa,z}) + \sqrt{\frac{3}{2}} h_{\xi,\eta}(z) \mathrm{Id} \right] \mathscr{I}_{P_{21}}(\sigma_{\kappa,z}, f) \right) + J_{M_{21}}(\sigma_{a,z}, f) \\
		&= J_{M_{21}}(\sigma_{\kappa,z}, f) + J_{M_{21}}(\sigma_{a,z}, f) + \sqrt{\frac{3}{2}} h_{\xi,\eta}(z) \mathrm{tr}\left( \mathscr{I}_{P_{21}}(\sigma_{\kappa,z}, f) \right),
	\end{align*}that is,
	\begin{align}\label{result}
		J_{M_{21}}(\sigma_{a,z}, f) =J_{M_{21}}(\sigma_{q,z}, f) -J_{M_{21}}(\sigma_{\kappa,z}, f) - \sqrt{\frac{3}{2}} h_{\xi,\eta}(z) \mathrm{tr}\left( \mathscr{I}_{P_{21}}(\sigma_{\kappa,z}, f) \right).
	\end{align}
	By \cite{MS04}, for \(\GGL(n_1)\times \GGL(n_2)\),
	$$
	r_{\overline{P}|P}(\pi_1 \otimes \pi_2, z) = \frac{L_F(z, \pi_1 \times \widetilde{\pi}_2)}{L_F(1 + z, \pi_1 \times \widetilde{\pi}_2)\epsilon_F(z, \pi_1 \times \widetilde{\pi}_2, \psi_F)}.
	$$
	In our case, $$
	\begin{gathered}
		r_a(z) = \frac{L_F(z, a_\xi \times \eta^{-1})}{L_F(1 + z, a_\xi \times \eta^{-1})\epsilon_F(z, a_\xi \times \eta^{-1}, \psi)},\qquad
		r_\kappa(z) = \frac{L_F(z, \kappa_\xi \times \eta^{-1})}{L_F(1 + z, \kappa_\xi \times \eta^{-1})\epsilon_F(z, \kappa_\xi \times \eta^{-1}, \psi)},
	\end{gathered}
	$$and 
	$$
	b_{\xi,\eta}(z) = \frac{L_F(z, a_\xi \times \eta^{-1})L_F(1 + z, \kappa_\xi \times \eta^{-1})}{L_F(1 + z, a_\xi \times \eta^{-1})L_F(z, \kappa_\xi \times \eta^{-1})} \cdot \frac{\epsilon_F(z, \kappa_\xi \times \eta^{-1}, \psi)}{\epsilon_F(z, a_\xi \times \eta^{-1}, \psi)}.
	$$
	Since \(\sigma_q=q_\xi\otimes\eta=\sigma_{1/2}\),
Lemma~\ref{lem:standard-vanishing-M21} gives
\[
I_{M_{21}}^G(\sigma_q,X,f)=0,
\qquad X\in\mathfrak a_{M_{21},F}.
\]
	
	The sequence \[0 \longrightarrow \sigma_\kappa \longrightarrow \sigma_q \longrightarrow \sigma_a \longrightarrow 0\] is exact, we have
	\begin{align}\label{exact}
		\mathrm{tr}(\sigma_q)=\mathrm{tr}(\sigma_\kappa)+\mathrm{tr}(\sigma_a). 
	\end{align} 
	Recall that
	\begin{align}\label{desc}
		I_{M_{21}}(\sigma_q,f)=J_{M_{21}}(\sigma_q,f)-\hat{I}_{M_{21}}^{M_{21}}(\sigma_q,\phi_{M_{21}}(f)),
	\end{align} according to \cite{Art88(1)}, we can still find a function \(u_{M_{21}}\in\mathscr{H}_{\rm{ac}}(M_{21}(F))\), such that
	\begin{align}\label{surj}
		\hat{I}^{M_{21}}_{M_{21}}(\sigma_q,\phi^G_{M_{21}}(f))=I^{M_{21}}_{M_{21}}(\sigma_q,\phi^{M_{21}}_{M_{21}}(u_{M_{21}}))=\rm{tr}\,\sigma_q(u_{M_{21}}).
	\end{align}
	Therefore, by (\ref{desc}),  we can write (\ref{result}) as 
	\begin{align}
		I_{M_{21}}(\sigma_a,f)=-\sqrt{\frac{3}{2}}\int_{i\mathfrak{a}_{M_{21}}^*/i\mathfrak{a}_G^*}h_{\xi,\eta}(z)\rm{tr}(\mathscr{I}_{P_{21}}(\sigma_{\kappa,z},f))dz.
	\end{align}
	Now we give the explicit formula of \(h_{\xi,\eta}(z)\).

	\textbf{The Non-Archimedean Case}
	
	Let \(F\) be a non-Archimedean field with residue field of cardinality \(q\). Let
	\[
	\mu = \xi \eta^{-1}.
	\]
	
	\textbf{1. \(\mu\) unramified}
	
	Set
	\[
	\alpha = \mu(\varpi), \quad |\alpha| = 1.
	\]
	Then
	\[
	L_F(s, \mu) = \frac{1}{1 - \alpha q^{-s}}, \quad \epsilon_F(s, \mu, \psi_F) = 1,
	\]
	where we choose the standard additive character of conductor zero. For the one-dimensional representation \(a_\xi\), we have
	\[
	L_F(z, a_\xi \times \eta^{-1}) = L_F\left(z + \frac{1}{2}, \mu\right) L_F\left(z - \frac{1}{2}, \mu\right).
	\]
	Thus,
	\[
	r_a(z) = \frac{L_F\left(z + \frac{1}{2}, \mu\right) L_F\left(z - \frac{1}{2}, \mu\right)}{L_F\left(z + \frac{1}{2}, \mu\right) L_F\left(z + \frac{1}{2}, \mu\right)} = \frac{L_F\left(z - \frac{1}{2}, \mu\right)}{L_F\left(z + \frac{1}{2}, \mu\right)}.
	\]
	On the other hand,
	\[
	\kappa_\xi = \mathrm{St}_2 \otimes (\xi \circ \det).
	\]
	Therefore,
	\[
	L_F(z, r_\kappa \times \eta^{-1}) = L_F\left(z + \frac{1}{2}, \mu\right).
	\]
	Since the epsilon factor of the Steinberg representation is
	\[
	\epsilon_F(z, r_\kappa \times \eta^{-1}, \psi_F) = -\alpha q^{1/2 - z},
	\]
	we obtain
	\[
	r_\kappa(z) = \frac{L_F\left(z + \frac{1}{2}, \mu\right)}{-\alpha q^{1/2 - z} L_F\left(z + \frac{3}{2}, \mu\right)}.
	\]
	Then
	\begin{equation}
		\begin{aligned}
			b_{\xi, \eta}(z) &= -\alpha q^{1/2 - z} \frac{L_F\left(z - \frac{1}{2}, \mu\right)}{L_F\left(z + \frac{1}{2}, \mu\right)} \\
			&= \frac{1 - \alpha q^{-z-1/2}}{1 - \alpha^{-1} q^{z-1/2}},
		\end{aligned}
	\end{equation}and
	\begin{equation}
		h_{\xi, \eta}(z) = (\log q) \left[ \frac{\alpha q^{-z-1/2}}{1 - \alpha q^{-z-1/2}} + \frac{\alpha^{-1} q^{z-1/2}}{1 - \alpha^{-1} q^{z-1/2}} \right].
	\end{equation}
	In particular, if \(\xi = \eta\), then \(\alpha = 1\), and thus
	\[
	h_{\xi, \xi}(0) = \frac{2 \log q}{\sqrt{q} - 1}.
	\]
	
	\textbf{2. \(\mu\) ramified}
	When
	\(\mu = \xi \eta^{-1}\)
	is ramified, by [ \cite{JL70} \S3, Proposition 3.6] and the local factor formula following the proof of Proposition 3.5, we have
	\[
	L_F(z, \alpha_\xi \times \eta^{-1}) = L_F(z, r_\kappa \times \eta^{-1}) = 1,
	\]and
	\[
	\epsilon_F(z, \alpha_\xi \times \eta^{-1}, \psi_F) = \epsilon_F(z, r_\kappa \times \eta^{-1}, \psi_F).
	\]
	Therefore,
	\[
	r_u(z) = r_\kappa(z),
	\]
	from which
	\begin{equation}
		b_{\xi, \eta}(z) = 1, \quad h_{\xi, \eta}(z) = 0.
	\end{equation}
	
	\textbf{Real Field \(F = \mathbb{R}\)}
	
	Let
	\[
	\mu = \xi \eta^{-1},
	\]
	and write
	\[
	\mu(x) = \mathrm{sgn}^\epsilon(x)|x|^{it}, \quad \epsilon \in \{0, 1\}.
	\]
	Let
	\[
	w = z + it.
	\]
	Using
	\[
	\Gamma_\mathbb{R}(s) = \pi^{-s/2}\Gamma(s/2), \quad \Gamma_\mathbb{C}(s) = 2(2\pi)^{-s}\Gamma(s),
	\]
	we have:
	\[
	L_F(z, \alpha_\xi \times \eta^{-1}) = \Gamma_\mathbb{R}\left(w + \frac{1}{2} + \epsilon\right) \Gamma_\mathbb{R}\left(w - \frac{1}{2} + \epsilon\right).
	\]
	Using
	\[
	\Gamma_\mathbb{R}(s+2) = \frac{s}{2\pi} \Gamma_\mathbb{R}(s)
	\]
	and the standard additive character, by [\cite{JL70}, \S5, p.93] 
	\[
	\epsilon_F(z, \alpha_\xi \times \eta^{-1}, \psi_F) = (-1)^\epsilon,
	\]
	we obtain
	\begin{equation}
		r_a(z) = (-1)^\epsilon \frac{2\pi}{w - \frac{1}{2}+\epsilon}.
	\end{equation}
	Now we calculate \(r_\kappa(z)\).
	
	By [\cite{Hum25}, \S 2.4.3, (2.13)], the local \(L\)-factor satisfies
	\begin{align*}
		L_F(z, r_\kappa \times \eta^{-1}) = \Gamma_{\mathbb{C}}\left(w + \frac{1}{2}\right),
	\end{align*}
	and by [\cite{Hum25}, \S 2.4.4, (2.18)]
	\[
	\epsilon_F(z, r_\kappa \times \eta^{-1}, \psi_F) = -1.
	\]
	Therefore,
	\begin{equation*}
		r_\kappa(z) = -\frac{2\pi}{w + \frac{1}{2}}.
	\end{equation*}
	Recall that
	\[
	b_{\xi, \eta}(z) = \frac{r_u(z)}{r_\kappa(z)},
	\]
	we obtain
	\begin{equation*}
		b_{\xi, \eta}(z) = \begin{cases} 
			-\frac{w + \frac{1}{2}}{w - \frac{1}{2}}, & \epsilon = 0, \\ 
			1, & \epsilon = 1.
		\end{cases}
	\end{equation*}
	Therefore,
	\begin{equation*}
		h_{\xi, \eta}(z) = \begin{cases} 
			\frac{1}{w + \frac{1}{2}} - \frac{1}{w - \frac{1}{2}}, & \epsilon = 0, \\ 
			0, & \epsilon = 1.
		\end{cases}
	\end{equation*}
	In particular, when \(z = iy\),
	\begin{equation*}
		h_{\xi, \eta}(iy) = \frac{1}{(y + t)^2 + \frac{1}{4}}, \quad \epsilon = 0.
	\end{equation*}
	
	\textbf{Complex Field \(F = \mathbb{C}\)}
	
	We set
	\[
	|u|_\mathbb{C} = u\bar{u}, \quad \omega(u) = \frac{u}{(u\bar{u})^{1/2}}.
	\]
	Write \(\mu = \xi \eta^{-1}\) as
	\[
	\mu(u) = \omega(u)^m |u|_\mathbb{C}^{it}, \quad m \in \mathbb{Z}.
	\]
	And let
	\[
	w = z + it.
	\]
	The two Gamma factors corresponding to \(\alpha_\xi\) are
	\begin{equation*}
		L_F(z, \alpha_\xi \times \eta^{-1}) = \Gamma_\mathbb{C}\left(w + \frac{|m|}{2} + \frac{1}{2}\right) \Gamma_\mathbb{C}\left(w + \frac{|m|}{2} - \frac{1}{2}\right).
	\end{equation*}
	At this time,
	\[
	\kappa_\xi = (\xi\omega) \times (\xi\omega^{-1}).
	\]
	Therefore,
	\begin{equation*}
		L_F(z, r_\kappa \times \eta^{-1}) = \Gamma_\mathbb{C}\left(w + \frac{|m+1|}{2}\right) \Gamma_\mathbb{C}\left(w + \frac{|m-1|}{2}\right).
	\end{equation*}
	When \(m \neq 0\),
	\[
	\left\{ \frac{|m+1|}{2}, \frac{|m-1|}{2} \right\} = \left\{ \frac{|m|}{2} + \frac{1}{2}, \frac{|m|}{2} - \frac{1}{2} \right\}.
	\]
	The two \(L\)-factors and the corresponding epsilon factors are the same, therefore
	\[
	r_a(z) = r_\kappa(z).
	\]
	When \(m = 0\), using
	\begin{align*}
		\Gamma_{\mathbb{C}}(s+1) = \frac{s}{2\pi} \Gamma_{\mathbb{C}}(s),
	\end{align*}
	we have
	\[
	\frac{r_a(z)}{r_\kappa(z)} =- \frac{w + \frac{1}{2}}{w - \frac{1}{2}}.
	\]
	Therefore,
	\begin{equation*}
		b_{\xi, \eta}(z) = \begin{cases} 
			-\frac{w + \frac{1}{2}}{w - \frac{1}{2}}, & m = 0, \\ 
			1, & m \neq 0,\qquad 
		\end{cases}
		h_{\xi, \eta}(z) = \begin{cases} 
			\frac{1}{w + \frac{1}{2}} - \frac{1}{w - \frac{1}{2}}, & m = 0, \\ 
			0, & m \neq 0.
		\end{cases}
	\end{equation*}
	On the imaginary axis \(z = iy\), if \(m = 0\), then
	\begin{equation*}
		h_{\xi, \eta}(iy) = \frac{1}{(y + t)^2 + \frac{1}{4}}.
	\end{equation*}

	\subsection{The case \(M=M_0\)}
	In this case, we can get that 
	\[I_{M_0}(\sigma,f)=0,\]for \(\sigma\in \Pi_{\rm{disc}}(M_1,t)\), by (\ref{tem}).

    Write the set of \(\sigma\) we discussed in subsection \ref{quotent} as \(\Sigma(M_{21},t)\). We have
    \begin{theorem}\label{spectr}
        For \(f\in\mathcal{H}(G(F_v))\), the invariant distribution \(I_M \left( \mathscr{I}_{P_1\cap M}^M(\sigma_\lambda), f \right)\) equals
\[ \rm{tr}({\pi}(f)),\]if \(M=G,\ \pi\in \Pi_{\rm{disc}}(G,t)\),\[-\sqrt{\frac{3}{2}}\int_{i\mathfrak{a}_{M_{21}}^*/i\mathfrak{a}_G^*}h_{\xi,\eta}(z)\rm{tr}(\mathscr{I}_{P_{21}}(\sigma_{\kappa,z},f))dz,\]if \(M=M_{21},\ \sigma\in\Sigma(M_{21},t)\).
    \end{theorem}

    Then according to the splitting formula in \cite{Art88(1)}:
    \begin{proposition}
        Let $\pi=\pi_1\otimes\pi_2$ and $X=(X_1,X_2)$ be as above.
Then for any function
\[
f=f_1f_2,\qquad
f_i\in\mathcal H_{\mathrm{ac}}\bigl(G(F_{S_i})\bigr),
\]
we have
\begin{align*}
    I_M(\pi,X,f)
=
\sum_{L_1,L_2\in\mathcal L(M)}
d_M^G(L_1,L_2)\,
\widehat I_M^{L_1}
\bigl(\pi_1,X_1,f_{1,L_1}\bigr)\,
\widehat I_M^{L_2}
\bigl(\pi_2,X_2,f_{2,L_2}\bigr).
\end{align*}
    \end{proposition}
We can write the invariant trace formula in global case.

	\subsection{The coefficient \(a^{M_1}_{\rm{disc}}(\sigma)\)}
	Let
	\[
	G_n = \mathrm{GL}(n), \quad E_n(\xi) = \rm{Ind}_{P_n}^{G_n}(\xi \otimes \cdots \otimes \xi), \quad n = 2, 3.
	\]
	It is easy to check that \[a^{M_{21}}_{\rm{disc}}(\sigma\otimes \eta)=a^{\GGL(2)}_{\rm{disc}}(\sigma),\qquad \sigma\in \Pi_{\rm{disc}}(\GGL(2)),\ \eta\in \Pi_{\rm{disc}}(\GGL(1)).\]
	
	\begin{theorem}\label{a_s}
		For \(n = 2, 3\),
		\[
		a_{\mathrm{disc}}^{G_n}(\pi) = 
		\begin{cases} 
			1, & \pi\in\Pi_{\rm{cusp}}(G_n)\  \text{or}\  \pi = \xi \circ \det, \\
			\frac{(-1)^{n-1}}{n^2}, & \pi = E_n(\xi),\\
			0, & \text{otherwise}.
		\end{cases}
		\]
		Furthermore, for \(M_0\) and \(M_{21}\),
		\[
		a_{\mathrm{disc}}^{M_0}(\chi_1 \otimes \chi_2 \otimes \chi_3) = 1, \qquad a_{\mathrm{disc}}^{M_{21}}(\tau \otimes \eta) = a_{\mathrm{disc}}^{\GGL(2)}(\tau).
		\]
	\end{theorem}
	\begin{proof}
		By the multiplicity one theorem for the cuspidal spectrum and the classification of the discrete spectrum by Moeglin-Waldspurger, the discrete spectrum of \(G_n\) consists of cuspidal representations and \(\xi \circ \det\), all with multiplicity 1.
		
		According to Arthur's discrete expansion, additional contributions arise from regular Weyl self-intertwinings of proper Levi subgroups. For \(\mathrm{GL}(3)\), and \(M_{21}\) satisfies
		\[
		W^{G_3}(\mathfrak{a}_{M_{21}})_{\mathrm{reg}} = \varnothing,
		\]
		Therefore, for \(n = 2, 3\), we only need to consider the diagonal torus \(T_n=\GGL(1)^n\). For the expansions used and the definitions of the coefficients, see [\cite{Art05}, (21.19), (22.5)].
		
		The regular elements of \(S_n\) in \(\mathfrak{a}_{T_n}^{G_n}\) are precisely the \(n\)-cycles, of which there are \((n - 1)!\), each satisfies
		\[
		\left| \det(1 - w \mid \mathfrak{a}_{T_n}^{G_n}) \right| = n.
		\]
		Only \(w\)-invariant inducing data contribute to the trace. For an $n$-cycle, such data are necessarily of the form \(\xi^{\otimes n}\). By [\cite{MW89}, \S II.1, (4)], the rank-one global self-intertwining operator corresponding to two identical characters is \(-\operatorname{Id}\). Combining this with the composition relations of step-by-step induction and intertwining operators, we obtain
		\[
		M(w,0)\big|_{E_n(\xi)} = (-1)^{n-1}\operatorname{Id}.
		\]
		Since $E_n(\xi)$ is irreducible, Arthur's expansion yields
		\[
		a_{\mathrm{disc}}^{G_n}(E_n(\xi)) = \frac{1}{n!} \, (n-1)! \, \frac{(-1)^{n-1}}{n} = \frac{(-1)^{n-1}}{n^2},
		\]
		which equals \(-\frac{1}{4}\) and \(\frac{1}{9}\) respectively, and there are no other proper Levi contributions. 
	\end{proof}
    Now we have 
    \begin{theorem}\label{main3}
        For any \(f\in \mathcal{H}(G(\mathbb{A}))\), the spectral side of invariant trace formula of \(\GGL(3)\) is 
        \begin{align*}
		I(f) = & \sum_{t>0}\sum_{M \in \mathcal{L}(T)} \sum_{M_1 \subset M}\frac{|W_0^{M_1}|}{|W_0^G|} \sum_{\sigma \in \Pi_{\mathrm{disc}}(M_1, t)} a_{\mathrm{disc}}^{M_1}(\sigma)\int_{i\mathfrak{a}_{M_1}^* / i\mathfrak{a}_M^*} r_{M_1}^M(\sigma, \lambda) I_M \left( \mathscr{I}_{P_1\cap M}^M(\sigma_\lambda), f \right) d\lambda.
	\end{align*}
    \end{theorem}

\end{document}